\pdfoutput=1
\RequirePackage{ifpdf}
\ifpdf 
\documentclass[pdftex]{sigma}
\else
\documentclass{sigma}
\fi

\usepackage{stmaryrd} 
\usepackage[scr=boondoxo, scrscaled=1.05]{mathalpha} 
\usepackage{fancyhdr}
\usepackage{enumitem}
\usepackage{xparse} 
\usepackage{tikz}
\usepackage{tikz-cd}
\usepackage{rotating} 

\let\Re\relax
\DeclareMathOperator{\Re}{Re}
\let\Im\relax
\DeclareMathOperator{\Im}{Im}

\usetikzlibrary{bbox}
\usetikzlibrary{fadings}

\pgfdeclareradialshading{radialedge}{\pgfpointorigin}{%
 color(0bp)=(transparent!0);
 color(20bp)=(transparent!0);
 color(22bp)=(transparent!10);
 color(24bp)=(transparent!90);
 color(25bp)=(transparent!100)
}
\pgfdeclarefading{radial edge}{\pgfuseshading{radialedge}}%

\newcommand{\holo}{\mathcal{H}}
\newcommand{\singexp}[2]{\mathcal{H}L^\infty_{#1, #2}}
\newcommand{\singexpalg}[1]{\singexp{#1}{\bullet}}
\newcommand{\dualsingexp}[1]{\widehat{\mathcal{H}}L^\infty_{#1}}

\newcommand{\maps}{\colon}

\newcommand{\Z}{\mathbb{Z}}
\newcommand{\R}{\mathbb{R}}
\newcommand{\C}{\mathbb{C}}

\newcommand{\volterra}{\mathcal{V}}
\newcommand{\hardpart}{\mathcal{V}_0}
\newcommand{\softpart}{\mathcal{V}_\star}

\newcommand{\solwhole}{f}
\newcommand{\solproto}{f_0}
\newcommand{\solptb}{f_\star}

\DeclareMathOperator{\var}{var}

\newcommand{\series}[1]{\tilde{#1}}
\newcommand{\fracderiv}[3]{\partial^{#1}_{#2, #3}}
\newcommand{\blankbox}{{\fboxsep 0pt \colorbox{lightgray}{\phantom{$h$}}}}
\newcommand{\van}{\mathfrak{m}}
\DeclareMathOperator{\Ai}{Ai}
\usetikzlibrary{matrix,shapes,arrows,decorations.pathmorphing}
\tikzset{commutative diagrams/arrow style=math font}
\tikzset{commutative diagrams/.cd,
mysymbol/.style={start anchor=center,end anchor=center,draw=none}}

\tikzset{
shift up/.style={
to path={([yshift=#1]\tikztostart.east)--([yshift=#1]\tikztotarget.west) \tikztonodes}
}
}
\DeclareMathAlphabet{\mathpzc}{OT1}{pzc}{m}{it}

\newcommand{\laplace}{\mathcal{L}}
\newcommand{\borel}{\mathcal{B}}
\newcommand{\aexp}{\text{\rm \ae}}

\definecolor{ietcoast}{RGB}{0, 150, 173}

\numberwithin{equation}{section}

\newtheorem{Theorem}{Theorem}[section]
\newtheorem{Corollary}[Theorem]{Corollary}
\newtheorem{Lemma}[Theorem]{Lemma}
\newtheorem{Proposition}[Theorem]{Proposition}
 { \theoremstyle{definition}
\newtheorem{Definition}[Theorem]{Definition}
\newtheorem{Remark}[Theorem]{Remark}
\newtheorem*{notation*}{Notation}
}

\usepackage{etoolbox}

\usetikzlibrary{math, 3d, calligraphy}

\definecolor{pwbeige}{RGB}{172, 146, 122}
\definecolor{pworange}{RGB}{245, 153, 30}

\tikzset{
  surf/.style={pwbeige!10},
  leaf/.style={thin, pwbeige!60},
  sing/.style={thick, pworange},
  ray/.style={thick}
}

\newcommand{\setfiberscale}[1]{
  \pgfmathsetmacro{\rfiber}{#1}
  \pgfmathsetmacro{\rshade}{1.5*\rfiber}
}
\setfiberscale{1.5}

\newcommand{\drawfiber}{
  \foreach \sgn in {-1, 1} {
    \calligraphy[scale=\sgn, copperplate, heavy, heavy line width=0.25mm]
      (\rfiber, \rfiber) -- (\rfiber, -\rfiber)
      (\rfiber, -\rfiber) -- (-\rfiber, -\rfiber);
  };
}

\pgfkeys{
  /laplace/.cd,
  drop/.store in=\fiberdrop,   drop/.value required,
  shade/.store in=\fibershade, shade/.default=true,
  degree/.store in=\fiberdeg,  degree/.value required,
  phase/.store in=\fiberphase, phase/.value required,
  dis/.store in=\fiberdis,     phase/.value required,
  dual/.store in=\fiberdual,   dual/.default=true
}

\newcommand{\fiberdescender}[1]{
  \calligraphy[copperplate, light, thin] (0, -\rfiber) -- (0, -#1);
}

\newcommand{\fibershading}[2]{
  \tikzmath{
    \shadestart = (-90 - #1) / \fiberdeg;
    \shadeend = (90 - #1) / \fiberdeg;
  }
  \fill[pwbeige!30, draw=pwbeige!60] (0, 0) ++(-#1:#2) -- ++(\shadestart:\rshade) arc (\shadestart:\shadeend:\rshade) -- cycle;
}

\newenvironment{fiber}[1][]{%
  \pgfkeys{/laplace/.cd, drop=0, shade=false, degree=1, phase=0, dis=0, #1}
  \begin{scope}
    \fill[surf] (-\rfiber, -\rfiber) rectangle (\rfiber, \rfiber);
    \clip (-\rfiber, -\rfiber) rectangle (\rfiber, \rfiber);
    \tikzmath{
      if \fibershade then {
        print \fibershading{\fiberphase}{\fiberdis};
      };
    }
}{
  \end{scope}
  \tikzmath{
    if \fiberdrop > 0 then {
      print \fiberdescender{\fiberdrop};
    };
  }
  \drawfiber
}

\newcommand{\genfibercontent}[1][]{
  \pgfkeys{/laplace/.cd, dual=false, #1}
  \pgfmathsetmacro{\hatch}{\rfiber/8};
  \pgfmathsetmacro{\firsthatch}{-\rfiber+\hatch}
  \pgfmathsetmacro{\nexthatch}{-\rfiber+2*\hatch}
  \pgfmathsetmacro{\hatchstop}{\rfiber-\hatch/2}
  \foreach \y in {\firsthatch, \nexthatch, ..., \hatchstop} {
    \tikzmath{
      coordinate \drawfrom;
      coordinate \drawto;
      if \fiberdual then {
        \drawfrom = (\y, -\rfiber);
        \drawto = (\y, \rfiber);
      } else {
        \drawfrom = (-\rfiber, \y);
        \drawto = (\rfiber, \y);
      };
    }
    \draw[leaf] (\drawfrom) -- (\drawto);
  };
  \fill circle (\dotsize);
}

\newcommand{\singfibercontent}{
  \foreach \sector in {30, 90, ..., 330} {
    \tikzmath{
      \farout = 45*floor((\sector + 59) / 45);
      \tfarout = 3*(\farout - \sector);
      \rfarout = \rfiber/sin(45+mod(\farout-45, 90));
    }
    \pgfmathsetmacro{\rmax}{sin(\tfarout) * pow(\rfarout, 3)}
    \foreach \r in {0.5, 1.0, ..., \rmax} {
      \pgfmathsetmacro{\win}{asin(\r/(2*\rmax))}
      \draw[leaf, rotate=\sector, domain=\win:{180-\win}, variable=\t] plot ({\t/3}:{(\r/sin(\t))^(1/3)});
    };
    
    \foreach \ang in {90, 210, 330} {
      \draw[sing] (\ang:{-sqrt(2)*\rfiber}) -- (\ang:{sqrt(2)*\rfiber});
    };
    \fill[sing] circle (\dotsize);
  };
}

\pgfmathsetmacro{\hatch}{2/3}
\tikzmath{
  \phase = 15;
  \dotsize = 0.04;
  \drop = 3.5;
  \xsing = 2;
  \zsing = 4*\hatch;
  \xgen = 6;
  \zgen = 8*\hatch;
  \rcut = 1/5;
}

\newcommand{\ray}[2]{
  \draw[ray] (#1, #2) -- (8, {#2 + (8-#1)*tan(\phase)});
}

\newcommand{\basicLaplace}{
\begin{tikzpicture}
\setfiberscale{2}

\begin{fiber}
  \genfibercontent
  \ray{-1}{0.4}
  \fill (-1, 0.4) circle (\dotsize);
  \node[anchor=north west, inner sep=2mm] at (-\rfiber, \rfiber) {$\zeta$};
\end{fiber}
\node[anchor=north, outer sep=1mm, font=\footnotesize]
  at (0, -\rfiber) {position space};

\begin{scope}[shift={(5.5, 0)}]
  \begin{fiber}[shade, phase=\phase, dis=0.6]
    \genfibercontent[dual]
    \node[anchor=north west, inner sep=2mm] at (-\rfiber, \rfiber) {$z$};
  \end{fiber}
  \node[anchor=north, outer sep=1mm, font=\footnotesize]
    at (0, -\rfiber) {frequency space};
\end{scope}
\end{tikzpicture}
} 

\newcommand{\phaseSpaceLaplace}{
\begin{tikzpicture}[z={(0.25cm, 0.2cm)}]
\begin{scope}[canvas is xz plane at y=0]
  \fill[surf]
    (0, 0) -- (0, \zsing-\rcut)
    -- ++(\xsing-\rcut, 0) arc (-90:90:\rcut) -- (0, {\zsing + \rcut})
    -- (0, 8) --(8, 8) -- (8, 0) -- cycle;
  
  \pgfmathsetmacro{\nexthatch}{2*\hatch}
  \pgfmathsetmacro{\singstop}{\zsing - \hatch/2}
  \foreach \z in {\hatch, \nexthatch, ..., \singstop} {
    \draw[leaf] (0, \z) -- (8, \z);
  };
  \pgfmathsetmacro{\singstart}{\zsing + \hatch}
  \pgfmathsetmacro{\singnext}{\zsing + 2*\hatch}
  \pgfmathsetmacro{\hatchstop}{8 - \hatch/2}
  \foreach \z in {\singstart, \singnext, ..., \hatchstop} {
    \draw[leaf] (0, \z) -- (8, \z);
  };
  
  \begin{scope}[sing]
    \draw (\xsing, \zsing) -- (8, \zsing);
    \draw (0, \zsing-\rcut) -- ++(\xsing-\rcut, 0) arc (-90:90:\rcut) -- (0, {\zsing + \rcut});
  \end{scope}
  
  \ray{\xgen}{\zgen}
  \ray{\xsing}{\zsing}
  
  \calligraphy[copperplate, heavy, heavy line width=0.25mm]
    (0, 0) -- (0, \zsing-\rcut)
    (0, \zsing+\rcut) -- (0, 8)
    (0, 8) -- (8, 8)
    (8, 8) -- (8, 0)
    (8, 0) -- (0, 0);
\end{scope}

\begin{scope}[canvas is xy plane at z=\zgen, shift={(\xgen, \drop)}]
  \begin{fiber}[drop=\drop, shade, phase=\phase]
    \genfibercontent[dual]
  \end{fiber}
\end{scope}

\begin{scope}[canvas is xy plane at z=\zsing, shift={(\xsing, \drop)}]
  \begin{fiber}[drop=\drop, shade, degree=3, phase=\phase]
    \singfibercontent
  \end{fiber}
\end{scope}

\fill (\xgen, 0, \zgen) circle (\dotsize);

\fill[sing] (\xsing, 0, \zsing) circle (\dotsize);

\begin{scope}[overlay]
  \node (surf label) at (9, 0, 4) {$B$};
  \path (surf label) ++(0, \drop, 0) node {$T^*B$};
  \node[outer sep=0.15mm, font=\footnotesize, align=left]
    (cut label) at (-1.25, 1.2, \zsing) {branch cut \\ at conical \\ singularity};
  \draw[->] (cut label) to[out=-75, in=180] (-0.25, 0, \zsing);
\end{scope}
\end{tikzpicture}
} 

\begin{document}

\allowdisplaybreaks

\newcommand{\arXivNumber}{2407.01412}

\renewcommand{\PaperNumber}{101}

\FirstPageHeading

\ShortArticleName{The Regularity of ODEs and Thimble Integrals with Respect to Borel Summation}

\ArticleName{The Regularity of ODEs and Thimble Integrals\\ with Respect to Borel Summation}

\Author{Veronica FANTINI~$^{\rm a}$ and Aaron FENYES~$^{\rm b}$}

\AuthorNameForHeading{V.~Fantini and A.~Fenyes}

\Address{$^{\rm a)}$~Laboratoire Math\'ematique d'Orsay, France}
\EmailD{\mail{veronica.fantini@universite-paris-saclay.fr}}
\URLaddressD{\url{https://sites.google.com/view/vfantini/home-page}}

\Address{$^{\rm b)}$~Studio Infinity, USA}
\EmailD{\mail{aaron.fenyes@fareycircles.ooo}}
\URLaddressD{\url{https://ooo.fareycircles.ooo}}

\ArticleDates{Received January 06, 2025, in final form November 05, 2025; Published online December 03, 2025}

\Abstract{Through Borel summation, one can often reconstruct an analytic solution of a~problem from its asymptotic expansion. We view the effectiveness of Borel summation as a regularity property of the solution, and we show that the solutions of certain differential equation and integration problems are regular in this sense. By taking a geometric perspective on the Laplace and Borel transforms, we also clarify why ``Borel regular'' solutions are associated with special points on the Borel plane. The particular classes of problems we look at are level~$1$ ODEs and exponential period integrals over one-dimensional Lefschetz thimbles. To expand the variety of examples available in the literature, we treat various examples of these problems in detail.}


\Classification{34M25; 34M30; 34M40; 40G10; 41A60}

\section{Introduction}

\subsection{The unreasonable effectiveness of Borel summation}\label{intro:summation}
You can often find a formal power series
\[ \series{\Phi} = \frac{c_0}{z^\tau} + \frac{c_1}{z^{\tau+1}} + \frac{c_2}{z^{\tau+2}} + \frac{c_3}{z^{\tau+3}} + \cdots , \]
with $\tau \in (0, 1]$, that looks or acts like a solution to a problem whose actual solutions are holomorphic functions of $z$. For example, if you want to understand how the solutions of the holomorphic ordinary differential equation (ODE)
\begin{equation}
\left[ \left[ \left(\frac{\partial}{\partial z}\right)^2 - 1 \right] + z^{-1} \frac{\partial}{\partial z} - \left(\frac{1}{3}\right)^2 z^{-2} \right] \Phi= 0 \label{eqn:bessel_rescaled_ex}
\end{equation}
behave near $z = \infty$, you might start by looking for formal {\em trans-monomial} solutions ${\rm e}^{-\alpha z} \series{\Phi}$, where $\series{\Phi}$ is a formal power series of the kind above.\footnote{Trans-mononials can be of different kinds: for example, fractional powers of $z$ can appear in the exponential. In this paper, however, we will only deal with trans-monomials of the form ${\rm e}^{-\alpha z} \series{\Phi}$ for $\alpha \in \C$.} Setting $\alpha = -1$ and $\tau = \frac{1}{2}$ gives a~well-behaved recurrence relation for $\series{\Phi}$, which produces the solution
\begin{equation}
{\rm e}^{z} \left[ \frac{(-1)!!}{z^{1/2}} + \frac{5}{72} \cdot \frac{1!!}{z^{3/2}} + \frac{385}{31104} \cdot \frac{3!!}{z^{5/2}} + \frac{17017}{6718464} \cdot \frac{5!!}{z^{7/2}} + \cdots \right] \label{series:bessel_ex}
\end{equation}
and its constant multiples (see \cite[equation~(10.40.1)]{dlmf}). As another example, you might rewrite the integral
\begin{equation}\label{int:bessel_ex}
\Phi = \int_{\mathcal{C}} \exp\bigl[-z \bigl(4 u^3 - 3 u\bigr)\bigr]\,{\rm d}u
\end{equation}
as
\[ {\rm e}^{z} \frac{1}{2\sqrt{3}} \int_{-\infty}^\infty {\rm e}^{-z t^2/2} \left[ 1 - \frac{\sqrt{3}}{9} t + \frac{5}{72} t^2 - \frac{4\sqrt{3}}{243} t^3 + \frac{385}{31104} t^4 - \frac{7\sqrt{3}}{2187} t^5 + \frac{17017}{6718464} t^6 - \cdots \right] {\rm d}t \]
using the substitution $\frac{1}{2} t^2 = 4u^3 - 3u + 1$. Na\"{i}vely integrating term by term, you again get the trans-monomial~\eqref{series:bessel_ex}, up to a constant.

Once you have the formal solution $\series{\Phi}$, you might try to get an actual solution by applying {\em Borel summation}, which turns a formal power series into a function asymptotic to it. Borel summation works in three steps.
\begin{enumerate}\itemsep=0pt
\item[(1)] Thinking of $z$ as a ``frequency variable'', we take the Borel transform (also known as formal inverse Laplace transform) of $\series{\Phi}$, producing a formal power series $\series{\phi}$ in a new ``position variable'' $\zeta$.
\item[(2)] If $\series{\phi}$ has a positive radius of convergence, we sum it to get a holomorphic function $\hat{\phi}$ on a~neighborhood of $\zeta = 0$. Then, by analytic continuation, we expand the domain of $\hat{\phi}$ to a~Riemann surface $B$ with a distinguished 1-form $\lambda$, which is the continuation of ${\rm d}\zeta$.
\item[(3)] 
If $\hat{\phi}$ grows slowly enough along an infinite ray $\zeta \in \alpha + {\rm e}^{{\rm i}\theta}(0, \infty)$, its Laplace transform~\smash{$\hat{\Phi} := \laplace_{\zeta, \alpha}^\theta \hat{\phi}$} turns out to be a holomorphic function of $z$, well-defined on some sector of the frequency plane. In this case, we say $\series{\Phi}$ is {\em Borel summable}, and we call $\hat{\Phi}$ its {\em Borel sum} at $\alpha$ in the direction $\theta$.
\end{enumerate}

\begin{figure}[h]
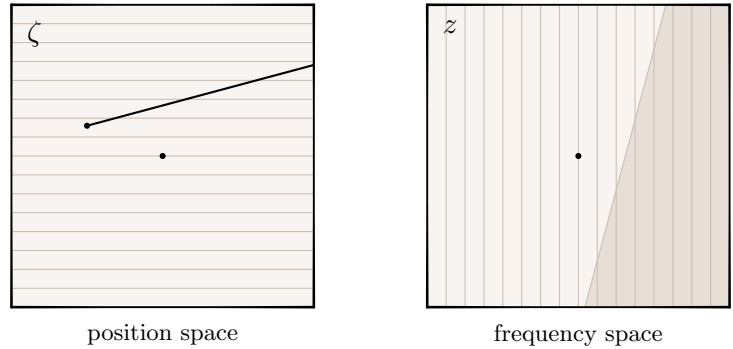
\centering
\basicLaplace
\caption{The Laplace transform \smash{$\laplace^\theta_{\zeta, \alpha}$} integrates a function along the ray $\zeta \in \alpha + {\rm e}^{{\rm i}\theta}(0, \infty)$ in the position domain, turning it into a function on some half-plane $\Re\bigl({\rm e}^{{\rm i}\theta} z\bigr) > \Lambda$ in the frequency domain.}
\end{figure}

The series $\series{\Phi}$ and its Borel sum $\hat{\Phi}$ have a special relationship, which is best described in the language of {\em Gevrey asymptoticity}.
\begin{Definition}
On an open, possibly bounded sector $\Omega$ around $\infty$, a holomorphic function $\Phi$ is {\em asymptotic} to a power series
\begin{equation}\label{series:asymp-definition}
a_0 + \frac{a_1}{z} + \frac{a_2}{z^2} + \frac{a_3}{z^3} + \frac{a_4}{z^4} + \cdots
\end{equation}
if, along each ray toward $\infty$,
\[ \left|\Phi - \left(a_0 + \frac{a_1}{z} + \frac{a_2}{z^2} + \cdots + \frac{a_n}{z^n} \right) \right| \in o(z^{-n}) \]
over all orders $n$. We will write
\smash{$\Phi \sim \sum_{n \geq 0} a_n z^{-n}$}
to denote this relationship.
\end{Definition}
This definition generalizes straightforwardly to ``fractionally shifted'' power series and trans-monomials.
\begin{Definition}
On an open, possibly bounded sector $\Omega$ around $\infty$, a holomorphic function $\Phi$ is {\em asymptotic} to a fractionally shifted trans-monomial
\[
{\rm e}^{\alpha z} z^{-\sigma} \left(a_0 + \frac{a_1}{z} + \frac{a_2}{z^2} + \frac{a_3}{z^3} + \frac{a_4}{z^4} + \cdots\right)
\]
if, along each ray toward $\infty$,
\[ \left|\Phi - {\rm e}^{\alpha z} z^{-\sigma} \left(a_0 + \frac{a_1}{z} + \frac{a_2}{z^2} + \cdots + \frac{a_n}{z^n} \right) \right| \in o\bigl({\rm e}^{-\alpha|z|} z^{-\sigma-n}\bigr) \]
over all orders $n$. We will again use $\sim$ to denote this relationship.
\end{Definition}

If a function is asymptotic to two power series of the form \eqref{series:asymp-definition}, those series must be equal; this can be deduced from the triangle inequality. The power series which a given function is asymptotic to on a given open sector, if it exists, is therefore an intrinsic property of the function: we call it the function's {\em asymptotic expansion} on that sector. The functions that have asymptotic expansions on a given sector form a ring~\cite[Appendix~A.4]{nikolaev2023existence}, and the map that sends such a function to its asymptotic expansion is a ring homomorphism. Since asymptotic expansions on overlapping sectors must agree, we can think of this homomorphism as being determined by a ray.
\begin{Definition}
For any direction $\vartheta$ and any trans-monomial space ${\rm e}^{\alpha z} z^{-\sigma} \C \big\llbracket z^{-1} \big\rrbracket$, the {\em asymptotic expansion} homomorphism
\[
 \aexp^\vartheta \colon\ \mathcal{O}(\C) \to {\rm e}^{\alpha z} z^{-\sigma} \C \big\llbracket z^{-1} \big\rrbracket\]
is the partially defined map that sends a holomorphic function on an open subset of $\C$ to its asymptotic expansion on a sector containing the end of the ray $z \in {\rm e}^{{\rm i}\vartheta}(0, \infty)$. To be in the domain of $\aexp^\vartheta$, a function just needs to have an asymptotic expansion on one such sector.
\end{Definition}
\begin{Definition}\label{def:unif-gevrey-asymp}
On an open sector $\Omega$ around $\infty$, a holomorphic function $\Phi$ is {\em uniformly $1$-Gevrey asymptotic} to a power series of the form \eqref{series:asymp-definition} if there is a constant $A \in (0, \infty)$ for which
\[ \left|\Phi - \left(a_0 + \frac{a_1}{z} + \frac{a_2}{z^2} + \cdots + \frac{a_{n-1}}{z^{n-1}} \right) \right| \lesssim \frac{A^n n!}{|z|^n} \]
over all orders $n$ and all $z \in \Omega$. We will use $\lesssim$ throughout the paper to mean ``bounded by a~constant multiple of'', with the same constant over all the specified variables.
\end{Definition}
To compare the Borel sum $\hat{\Phi}$ from step~(3) with the original series $\series{\Phi}$, let us take it one more step, sending it back into the world of formal power series by taking its asymptotic expansion.
\begin{enumerate}[start=4]\itemsep=0pt
\item[(4)] By construction, $\aexp^{-\theta} \hat{\Phi} = \series{\Phi}$. It turns out that $\hat{\Phi}$ is not only asymptotic to $\series{\Phi}$, but {\em {uniformly} $1$-Gevrey}-asymptotic~\cite[Corollary~5.23]{diverg-resurg-i}.
\end{enumerate}
Note that the Borel sum and the asymptotic expansion are taken along the complementary directions $\theta$ and $-\theta$, respectively.

The Borel summation process is summarized in the following diagram:
\begin{center}
\begin{tikzcd}
& \text{problem} & \\
\hat{\Phi} \arrow[ru, dotted, no head, tail] \arrow[rr, mapsto, "\aexp^{-\theta}", "\text{($1$-Gevrey)}"'] & & \series{\Phi} \arrow[lu, no head, tail, "\parbox{11mm}{\centering\scriptsize formally solves}"'] \arrow[dd, mapsto, "\borel"] \\
& & \\
\hat{\phi} \arrow[uu, mapsto, "\laplace_b^\theta"] & & \series{\phi} \arrow[ll, mapsto, "\text{sum}"].
\end{tikzcd}
\end{center}

You cannot be sure {\em a priori} that $\hat{\Phi}$ solves your original problem, even if you know that $\series{\Phi}$ is the asymptotic expansion of an actual solution. After all, $\series{\Phi}$ is asymptotic to many functions; Borel summation simply picks one of them. In many cases, however, Borel summation picks correctly, delivering an actual solution to your problem. Borel recognized this phenomenon as a~key motivation for the study of divergent series \cite{borel-memoire,borel-lecons}, and much has been said about it since (see Section~\ref{sec:historical-context}); the question of how it happens is the starting point for this paper.

\subsection{A new perspective: Borel regularity}

\subsubsection{Introducing Borel regularity}
We will look at the effectiveness of Borel summation from an analytic perspective, as a regularity property of problems and their holomorphic solutions. Suppose we have a holomorphic function~$\Phi$ which is a solution of some problem. If the asymptotic expansion $\series{\Phi} = \text{\ae}^{-\theta} \Phi$ is well defined, it should be a formal solution of the same problem; we can take this as a basic requirement for our concept of a formal solution:
\begin{center}
\begin{tikzcd}
& \text{problem} & \\
\Phi \arrow[ru, no head, tail, swap, "\text{solves}"' ] \arrow[rrd, mapsto, "\aexp^{-\theta}"]& & \\
\hat{\Phi} \arrow[rr, "\aexp^{-\theta}", "\text{($1$-Gevrey)}"'] & & \series{\Phi} \arrow[luu, no head, tail, swap, "\parbox{11mm}{\centering\scriptsize formally solves}"] \arrow[dd, mapsto, "\borel_\zeta"] \\
& & \\
\hat{\phi} \arrow[uu, mapsto, "\laplace_{\zeta, 0}^\theta"] & & \series{\phi} \arrow[ll, mapsto, "\text{sum}"]
\end{tikzcd}
\end{center}
If $\series{\Phi}$ is Borel summable, as described in Section~\ref{intro:summation}, its Borel sum $\hat{\Phi}$ is a new holomorphic function.
Since different functions can have the same asymptotic expansion, taking the Borel sum of the asymptotic expansion of a function must smooth away some details. We will therefore call this process {\em Borel regularization}. Explicitly, Borel regularization works in four steps:
\begin{enumerate}\itemsep=0pt
\item[(1)] Take the asymptotic expansion $\series{\Phi} = \aexp^{-\theta} \Phi$.
\item[(2)] Take the Borel transform $\series{\phi} = \borel_\zeta \series{\Phi}$.
\item[(3)] Take the sum $\hat{\phi}$ of $\series{\phi}$, and expand its domain to a Riemann surface $B$ with a distinguished $1$-form, as before.
\item[(4)] Take the Laplace transform \smash{$\hat{\Phi} = \laplace_{\zeta, 0}^\theta \hat{\phi}$}. This is possible, by definition, if $\series{\Phi}$ is Borel summable.\footnote{More generally, we could take the Laplace transform $\laplace^\theta_{\zeta, \alpha}$ starting at any point $\zeta = \alpha$ in $B$. The choice of starting point will play an important role in our study of ODEs and thimble integrals.}
\end{enumerate}
We will say a function is {\em Borel regularizable} if its asymptotic expansion is well defined and Borel summable, ensuring that we can carry out all the steps.

Defining a regularization process picks out a class of regular functions: the ones that are invariant under regularization. We will say a function is {\em Borel regular} if it is Borel regularizable and Borel regularization leaves it unchanged. In other words, $\Phi$ is Borel regular if $\hat{\Phi} = \Phi$.

\subsubsection{Borel regularity sometimes explains why Borel summation works}\label{borel-reg:explanatory-power}
The idea of Borel regularity can help us understand why Borel summation is so effective in some situations. Roughly speaking, Borel summation works well for problems that admit solutions in terms of the Laplace transform.

The central goal of this paper is to explain, from this perspective, why Borel summation works well for the two kinds of problems exemplified in Section~\ref{intro:summation}.
\begin{enumerate}\itemsep=0pt
\item[(1)] Solving a level~$1$ linear ODE \cite[Section 2.1]{EcalleIII}, \cite[Section~5.2.2.1]{diverg-resurg-iii}.\footnote{Calling a linear ODE {\em level~$1$} is a short way of saying that it has ``a single level, equal to $1$, at infinity''. That means a certain Newton polygon associated with its singularity at infinity has exactly one positive-slope side, which has slope~$1$. In general, the {\em levels} of a linear ODE at a given singularity are the slopes of the positive-slope sides of the associated Newton polygon~\cite[Proposition~3.3.14]{diverg-resurg--ii}.} Equation~\eqref{eqn:bessel_rescaled_ex} is an example. More generally, we will consider equations of the form $\mathcal{P} \Phi = 0$ given by a~differential operator of the form
\[ \mathcal{P} = P\left(\frac{\partial}{\partial z}\right) + \frac{1}{z} Q\left(\frac{\partial}{\partial z}\right) + \frac{1}{z^2} R\bigl(z^{-1}\bigr), \]
where $P$ is a monic degree-$d$ polynomial, $Q$ is a degree-$(d-1)$ polynomial that is non-zero at every root of $P$, and $R\bigl(z^{-1}\bigr)$ is holomorphic in some disk $|z| > A$ around $z = \infty$. We will restrict our attention to the case where the roots of $P$ are simple (see Section~\ref{borel_reg-ODE}), and to the solutions associated with roots $-\alpha$ where \smash{$\tau_\alpha:=\frac{Q(-\alpha)}{P'(-\alpha)}$} is real and positive (see Theorem~\ref{thm:exist_uniq_ODE}).
The form of $\mathcal{P}$ derives from \cite[equation~(2.2.3), p.~105]{EcalleIII}, which encompasses both linear and nonlinear ODEs of level $1$. Borel summation is more involved for the nonlinear ones, raising the question of whether our analysis generalizes.
\item[(2)] Evaluating a one-dimensional {\em thimble integral}. Integral~\eqref{int:bessel_ex} is an example. More generally, we will consider integrals of the form
\[
I = \int_{\mathcal{C}}{\rm e}^{-zf} \nu,
\]
where $X$ is a~one-dimensional complex manifold; $f \maps X \to \C$ is a holomorphic function with isolated, non-degenerate critical points; $\nu \in \Omega^1(X)$ is holomorphic $1$-form on $X$; and $\mathcal{C}$ is {\em Lefschetz thimble}---a type of contour described in Section~\ref{sec:why_borel_thimble}. Under these conditions, $I$ is a holomorphic function of $z$.
\end{enumerate}
These two problems are closely linked. By playing with derivatives of a thimble integral, you can often find a linear ODE that the integral satisfies. Conversely, for many classical ODEs, there are useful bases of thimble integral solutions.
\subsection{Goals and Results}
We present the results, proofs, and examples in this paper with a few central goals in mind. Our first goal is to discuss Borel summability in a way that carefully separates the analytic and formal sides of the theory. In particular, we highlight results that can be proven in purely analytic ways, like Theorem~\ref{thm:exist_uniq_ODE} for level~$1$ ODEs and Lemma~\ref{lem:thimble_proj_formula} for thimble integrals.

Our second goal is to emphasize how the problems we study become more regular when posed in the position domain (the Borel plane). When we take level~$1$ differential equations from the frequency domain to the position domain, we not only turn their irregular singularities into regular ones, but also turn them into integral equations, which provide more smoothness and can be solved more constructively. When we write thimble integrals as Laplace transforms, their irregularly shaped integration contours collapse to rays in the position domain.

Our third goal is to illustrate the results we discuss with detailed examples involving a~diversity of Borel regular functions. We treat several examples both as thimble integrals and as solutions of level~$1$ ODEs, so you can compare these parallel explanations of how Borel regularity arises. We include ODEs of both second degree and higher degree. Using the widely discussed example of the Airy function, we compare our treatment with other approaches and conventions found in the literature.

We state our main results in Sections~\ref{sec:why_borel_ODE}--\ref{sec:why_borel_thimble} and prove them in Section~\ref{sec:proof_main_results}. Our examples are given in Section~\ref{sec:examples}. Section~\ref{sec:org} describes the layout of the paper in more detail.
\subsubsection[Why does Borel summation work for solutions of level 1 ODEs?]{Why does Borel summation work for solutions of level 1 ODEs?}\label{sec:why_borel_ODE}

Consider a linear level~$1$ differential operator $\mathcal{P}$ of the form described in~Section~\ref{borel-reg:explanatory-power}. This operator will always have an irregular singularity at $z = \infty$. We can find holomorphic solutions of the equation $\mathcal{P}\Phi = 0$ by taking the Borel sums of formal trans-monomial solutions. However, this only works when there is a coincidence between the exponential factor ${\rm e}^{-\alpha z}$ in the formal solution, the base point $\zeta = \alpha$ we use for Borel summation, and the characteristic equation~${P(-\alpha) = 0}$. Why?

A first clue is the observation that every solution we get from Borel summation is, by definition, the Laplace transform \smash{$\Phi = \laplace_{\zeta,\alpha}^\theta \phi$} of a function $\phi$ on the position domain. In fact, the differential equation $\mathcal{P}\Phi = 0$ is equivalent to an integral equation $\hat{\mathcal{P}}_\alpha \phi = 0$.

Furthermore, when $-\alpha$ is a root of $P$, the integral equation $\hat{\mathcal{P}}_\alpha \phi = 0$ has a special solution $\psi_\alpha$ which is distinguished, up to scaling, by its power-law asymptotics at $\zeta = \alpha$~\cite[Theorem~4.6]{reg-sing-volterra}. This solution is regular enough to have a well-defined Laplace transform $\Psi_\alpha$, which is a special solution of $\mathcal{P}\Phi = 0$ distinguished by its own characteristic asymptotics.
\begin{Theorem}[Theorem~\ref{re:thm:exist_uniq_ODE}]\label{thm:exist_uniq_ODE}
Let $-\alpha$ be a root of $P$ where $\tau_\alpha$ is real and positive. Consider an open sector $\Omega_\alpha$ which has an opening angle of $\pi$ or less, has $\zeta = \alpha$ at its tip, and does not touch any other root of $P(-\zeta)$. The equation $\mathcal{P}\Phi = 0$ has a unique solution $\Psi_\alpha$ in the affine subspace~\smash{${\rm e}^{-\alpha z} \bigl[ z^{-\tau_\alpha} + \dualsingexp{-\tau_\alpha-1}\bigl(\widehat{\Omega}_\alpha^\bullet\bigr) \bigr] $}
of the space ${\rm e}^{-\alpha z} \dualsingexp{-\tau_\alpha}\bigl(\widehat{\Omega}_\alpha^\bullet\bigr)$ from Definition~{\rm\ref{def:H_hat}}.
\end{Theorem}
The special solution $\psi_\alpha$ is also regular enough to have a ``fractionally shifted'' power series expansion $\series{\psi}_\alpha$, whose formal Laplace transform $\series{\Psi}_\alpha$ is a formal solution of $\mathcal{P}\series{\Phi} = 0$. By construction, $\series{\Psi}_\alpha$ is Borel summable, with Borel sum $\Psi_\alpha$.
\begin{Theorem}[Theorem~\ref{re:thm:soln_is_Borel_sum}]\label{thm:soln_is_Borel_sum}
The analytic solution $\Psi_\alpha$ from Theorem~{\rm\ref{thm:exist_uniq_ODE}} is the Borel sum of a~formal trans-monomial solution $\series{\Psi}_\alpha \in {\rm e}^{-\alpha z} z^{-\tau_\alpha} \C \big\llbracket z^{-1} \big\rrbracket$ of the equation $\mathcal{P}\series{\Phi} = 0$.
\end{Theorem}

We can use this result to show that $\Psi_\alpha$ is Borel regular. Because of the fractional power~$z^{-\tau_\alpha}$, we use Lemma~\ref{lem:laplace-bridge} instead of a more standard result on the asymptotics of a Borel sum, such as \cite[Theorem~5.20]{diverg-resurg-i}.
\begin{Corollary}[Corollary~\ref{re:cor:soln_borel-reg}]\label{cor:soln_borel-reg}
The analytic solution $\Psi_\alpha$ from Theorem~{\rm\ref{thm:exist_uniq_ODE}} is Borel regular. This is because $\Psi_\alpha$ is asymptotic to the formal solution $\series{\Psi}_\alpha$ from Theorem~{\rm\ref{thm:soln_is_Borel_sum}}.
\end{Corollary}
In many cases, by choosing from among the solutions $\Psi_\alpha$, we can build a frame of Borel-regular solutions of the equation $\mathcal{P}\Phi = 0$.

From an analytic perspective, Theorem~\ref{thm:soln_is_Borel_sum} justifies the ansatz that Poincar\'e used to solve the equation $\mathcal{P}\series{\Phi} = 0$. On the other hand, from a formal perspective, one can start from the computation of the Poincar\'e solutions and the observation that their Borel transforms converge, and then use Theorem~\ref{thm:exist_uniq_ODE} to establish their Borel summability. We do this in Section~\ref{sec:new-summability-proof}.

As discussed in Section~\ref{sec:history_ODE}, the novel aspect of Theorem~\ref{thm:exist_uniq_ODE} is that it is stated and proved in the world of functional analysis, rather than the world of formal series. Theorem~\ref{thm:soln_is_Borel_sum} is a~new proof of a special case of an old result mentioned in Section~\ref{sec:Poincare method}: the multisummability of formal solutions of linear systems of ODEs~\cite{braaksma2006laplace,ramis1991series}. The new proof is interesting because it starts from the analytic solutions found in Theorem~\ref{thm:exist_uniq_ODE}, rather than the formal solutions found by Poincar\'e.

\subsubsection{Why does Borel summation work for thimble integrals?}\label{sec:why_borel_thimble}
In Section~\ref{borel-reg:explanatory-power}, we gave the general form of a thimble integral associated with a map $f \maps X \to \C$. The important thing about $\C$ here turns out to be its {\em translation surface} structure, which picks out an atlas of charts related by translations and provides a global sense of direction. The geometry of translation surfaces, described more fully in Section~\ref{sec:transl}, provides a natural setting for the Laplace transform. It will therefore be useful to place thimble integrals in this setting as well.

Given a translation chart $\zeta$ on a translation surface $B$, we write $\mathcal{J}_{\zeta, \alpha}^\theta$ to denote the ray in direction $\theta$ starting at $\zeta = \alpha$. Consider a map $f \maps X \to B$ with a critical value at $\zeta = \alpha$. Each connected component of the preimage $f^{-1}\bigl(\mathcal{J}_{\zeta, \alpha}^\theta\bigr)$ runs through one of the critical points that $f$ sends to $\zeta = \alpha$. When the complex manifold $X$ is one-dimensional, as it will be throughout this paper, a {\em Lefschetz thimble} is an oriented component of $f^{-1}\bigl(\mathcal{J}_{\zeta, \alpha}^\theta\bigr)$ that runs through a single, non-degenerate critical point of $f$. A one-dimensional Lefschetz thimble is determined up to orientation by the critical point $a$ that it runs through and the direction $\theta$ of the ray that it projects to, so we will often denote it with a symbol like $\mathcal{C}^\theta_a$. Because, by our definition, a~Lefschetz thimble can only run through a single critical point $a$, the ray that it projects to cannot touch any critical value except $f(a)$.
\begin{figure}[!ht]\centering
\includegraphics[width=6.5cm]{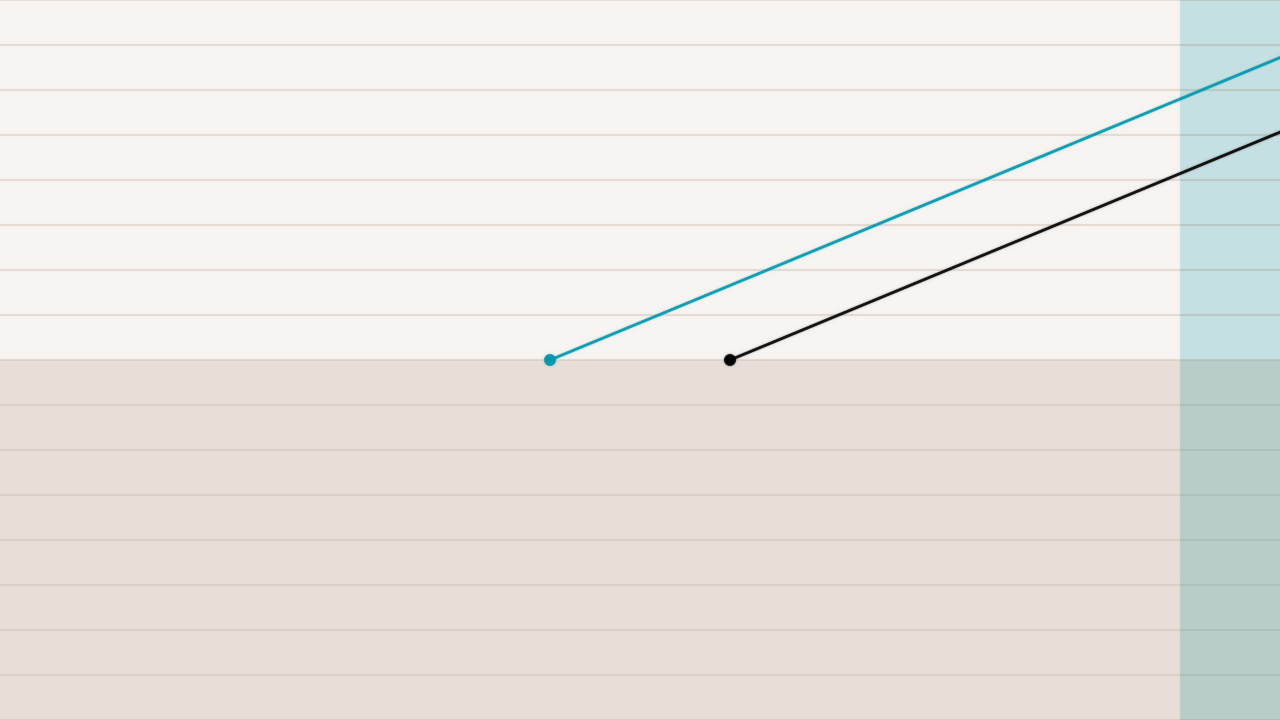}
\hspace{0.5cm}
\includegraphics[width=6.5cm]{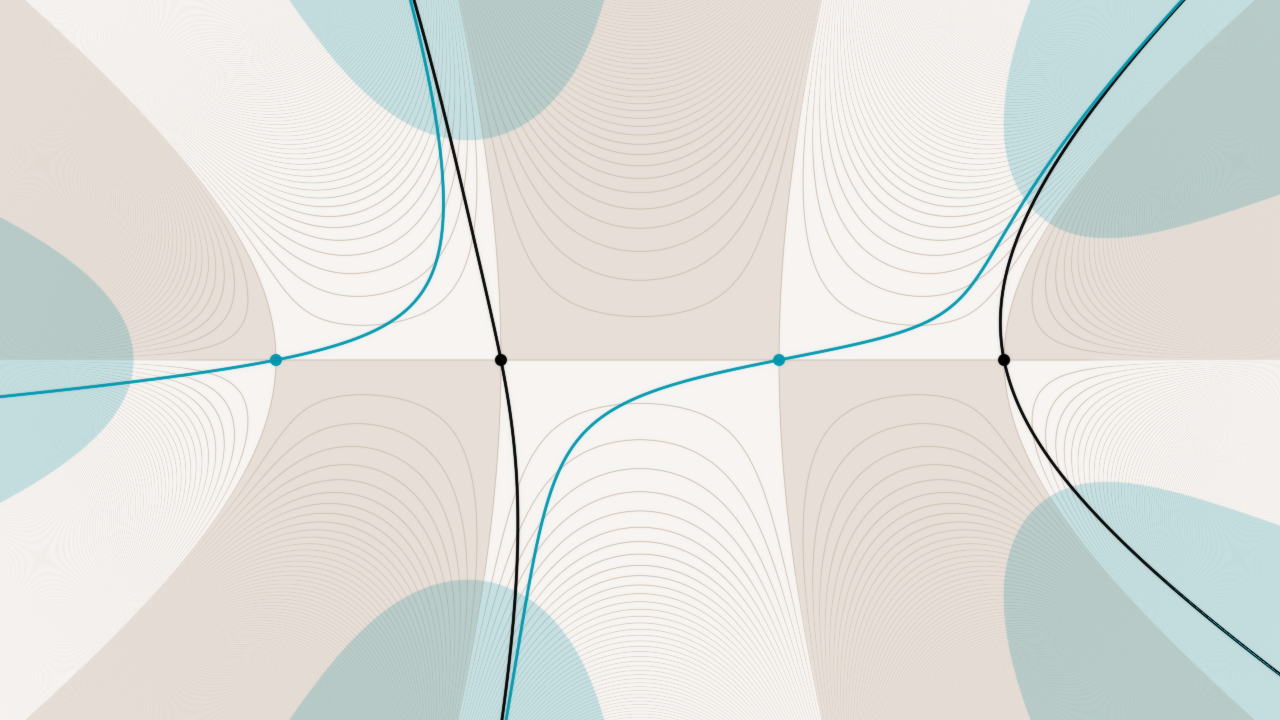}
\caption{On the left, we see the rays \smash{$\mathcal{J}^{\pi/8}_{\zeta, \pm 1}$} on the translation surface $B = \C$, where $\zeta$ is the standard coordinate. On the right, we see the complex manifold $X = \C$, colored according to a~holomorphic map $f \maps X \to B$ whose critical values are $\zeta = \pm 1$. The Lefschetz thimbles over the rays \smash{$\mathcal{J}^{\pi/8}_{\zeta, \pm 1}$} are shown.}
\end{figure}

We are now ready to define a one-dimensional thimble integral precisely, at our desired level of generality.
\begin{Definition}
Fix a one-dimensional complex manifold $X$; a translation surface $B$; and a~holomorphic map $f \maps X \to B$ with isolated, non-degenerate critical points. A one-dimensional {\em thimble integral} has the form
\[
 I_a = \int_{\mathcal{C}_a^\theta} {\rm e}^{-zf} \nu,
 \]
where $\nu$ is a holomorphic $1$-form on $X$, $a$ is a critical point of $f$, and $\mathcal{C}^\theta_a$ is the Lefschetz thimble that runs through $a$ in direction $\theta$ with a chosen orientation. Since $I_a$ depends on the frequency coordinate $z$, it is a holomorphic function on the frequency domain.
\end{Definition}
We will present a proof that one-dimensional thimble integrals are always Borel regular. The proof hinges on the following well-known formula, which can be used to express any one-dimensional thimble integral as a Laplace transform.
\begin{Lemma}[adapted from {\cite[Section~3.3]{pham}}]\label{lem:thimble_proj_formula}
A function $\iota_a$ with \smash{$I_a = \laplace_{\zeta, \alpha}^\theta \iota_a$} is given by the {\em thimble projection formula}
\[
 \iota_a = \frac{\partial}{\partial \zeta} \biggl( \int_{\mathcal{C}_a^\theta(\zeta)}\nu \biggr),
 \]
where $\mathcal{C}_a^\theta(\zeta)$ is the part of $\mathcal{C}_a^\theta$ that goes through $f^{-1}\bigl(\bigl[\alpha,\zeta {\rm e}^{{\rm i}\theta}\bigr]\bigr)$.
\end{Lemma}
Notice that the path $\mathcal{C}_a^\theta(\zeta)$ used in the lemma starts and ends in $f^{-1}(\zeta)$.

Once we know $\iota_a$, we can go through the steps of the Borel regularization process and show that we end up with the same function we started with.
\begin{Theorem}[Theorem \ref{thm:maxim-proof}]\label{thm:maxim}
If the integral defining $I_a$ is absolutely convergent, then $I_a$ is Borel regular. More explicitly,
\begin{enumerate}\itemsep=0pt
\item[$(1)$] 
The function $I_a$ has an asymptotic expansion $\series{I}_a := \aexp^{-\theta} I_a$, which lies in the space ${\rm e}^{-z \alpha} \smash{ z^{-\frac{1}{2}}} \allowbreak\times \C\big\llbracket z^{-1}\big\rrbracket$. Here, $\theta$ is the direction of the ray $\mathcal{J}^\theta_{\zeta, \alpha}$ that defines the thimble.
\item[$(2)$] 
The Borel transform $\series{\iota}_a := \borel_\zeta \series{I}_a$ converges near $\zeta = \alpha$.
\item[$(3)$] 
The sum of $\series{\iota}_a$ can be analytically continued along the ray $\mathcal{J}_{\zeta, \alpha}^\theta$. Its Laplace transform along that ray is well defined, and equal to $I_a$.
\end{enumerate}
\end{Theorem}
For one-dimensional thimbles, this result confirms the expectation that evaluating a thimble integral should be equivalent to taking the Borel sum of its asymptotic series~\cite[Section~2]{dunne-unsal}. Confirming this expectation in higher dimensions would be an interesting topic for future work, since our argument may not generalize straightforwardly.

As an alternative to our argument, Theorem~\ref{thm:maxim} could be deduced from Nevanlinna's theorem~\cite{nevanlinna}. However, that would require a careful asymptotic analysis of $I_a$: one would have to show that its asymptotic series is $1$-Gevrey and uniformly asymptotic to it on a large enough sector. We avoid this asymptotic analysis by working in the position domain and recasting $I_a$ as a Laplace transform integral.
\subsubsection{Other results}\label{sec:other_results}
\begin{figure}[t]
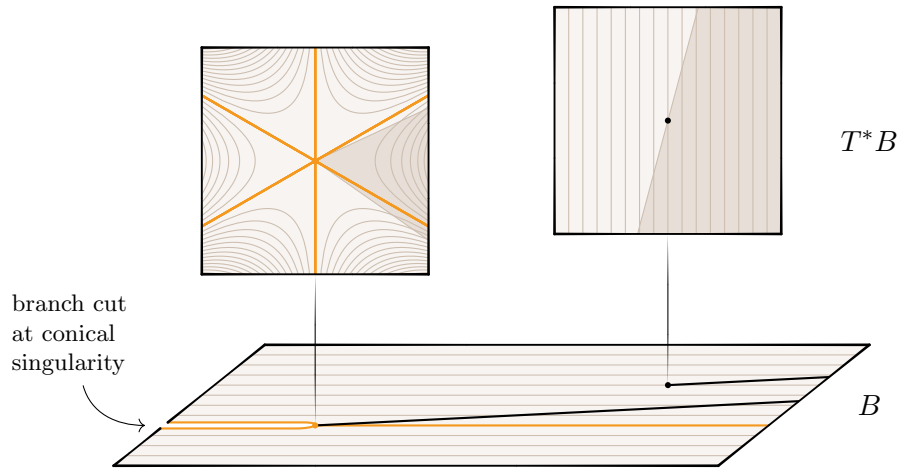
\centering
\phaseSpaceLaplace

\vspace{2mm}

\caption{The frequency coordinate $z$ on the cotangent spaces of an ordinary point and a~singular point. The singularity shown here has cone angle $6\pi$, like the singularities of the translation surface associated with the Airy function.}\label{fig:geometric-picture-intro}
\end{figure}

{\bf New perspectives on the Laplace transform.}
For Borel summation to work as described in Theorems~\ref{thm:soln_is_Borel_sum} and~\ref{thm:maxim}, one usually needs to take the Laplace and Borel transforms using a particular position coordinate and integration base point. To explain this aspect of Borel regularity, we introduce a geometric picture of the Laplace and Borel transforms, generalizing the position domain from the complex plane to a translation surface $B$. Each translation chart~$\zeta$ gives a different version of the Laplace transform, whose associated frequency domain is the cotangent space $T^*_{\zeta = 0} B$. The frequency coordinate $z$ is a canonical coordinate on $T^*_{\zeta = 0}B$, or on~\smash{$T^*_{\zeta = 0}B^{\otimes n}$} if $\zeta = 0$ is a singularity with cone angle $2\pi n$, as illustrated in Figure~\ref{fig:geometric-picture-intro}. The details of this perspective are discussed in Sections \ref{sec:geometry_laplace} and \ref{sec:geometry_borel}.

We depart slightly from traditional Laplace transform methods for solving ODEs. One traditionally relates ODEs on the frequency domain to ODEs on the position domain~\cite{braaksma2006laplace,laplace-tfm}, but we find it more natural to relate them to integral equations on the position domain. The good behavior of integral equations and their solutions is central to our proof that certain ODEs have Borel regular solutions.

{\bf New and old examples.}
We illustrate our main results with detailed treatments of several examples. Generalizing the classic example of the Airy equation (see Appendix~\ref{airy-appendix}), we find frames of Borel regular solutions for the Airy--Lucas and generalized Airy equations (see Sections~\ref{example_AL} and \ref{sec:gen-airy}). Noting that the Airy--Lucas equations reduce to special cases of the modified Bessel equation, we also find Borel regular solutions of the general modified Bessel equation (see Section~\ref{sec:mod-bessel-lift}). Similarly, noting that the Airy function is a thimble integral with a third-degree polynomial in the exponent, we show that a general third-degree thimble integral is Borel regular (see Section~\ref{sec:deg3}). Finally, as an example of a fourth-order ODE, we solve the equation describing the vibration of a triangular cantilever (see Section~\ref{sec:catilever}).

Each of our examples focuses on either a level~$1$ ODE whose solutions can be expressed as thimble integrals or a thimble integral that satisfies a level~$1$ ODE. This provides an opportunity to compare different explanations for why each example function is Borel regular.

In the literature on Borel summation, the Airy equation has been discussed many times, using different approaches and conventions. In Appendix~\ref{airy-comparison}, we explain how the treatments in~\cite[Section 2.2]{lectures-Marino}, \cite[Section 6.14]{diverg-resurg-i}, and \cite[Section 2.2]{kawai-takei} line up with ours.

{\bf Links to resurgence.}
The examples we study all involve {\em resurgent functions} in the position domain. The theory of resurgence, introduced by \'Ecalle~\cite{EcalleI,EcalleII,EcalleIII}, describes divergent series in the frequency domain in terms of the analytic properties of their Borel transforms. Resurgent analysis is an effective tool for computing ``subleading contributions'' in perturbation theory, and it occasionally leads to the resummation of the original divergent series. It has been applied extensively to linear ODEs \cite{loday1994stokes,diverg-resurg--ii}, and it has yielded results about nonlinear ODEs as well~\cite{schiappa-PI,costin-PI,costin_kruskal,diverg-resurg-iii}.

Geometric arguments show that when a thimble integral has a holomorphic Morse function in the exponent, as we have in at least some of our examples, the corresponding function on the position domain is always resurgent~\cite{Maxim_slide_ERC}, \cite[Section~6.2]{kontsevich2022analyticity}. Conjecturally, this property extends to more general thimble integrals---perhaps even infinite-dimensional ones \cite[Examples~5--6]{Maxim_slide_ERC}. The linear ODEs and $1$-dimensional thimble integrals that we study provide toy examples of how resurgent functions arise and behave. As a concrete example, in Section~\ref{resurgence-AL}, we use resurgence to find the Stokes constants of the Airy--Lucas equations.
\subsection{Organization of the paper}\label{sec:org}
In Section~\ref{sec:historical-context}, we contextualize our work by reviewing some classical results related to Borel regularity. These results come from the literature on asymptotics, ODEs, and integrals over Lefschetz thimbles.

In Section~\ref{sec:Laplace-Borel-general}, we introduce a geometric perspective on the Laplace and Borel transforms, which will be useful for stating and proving our results. To make this section more self-contained, we review basic definitions and properties in Section~\ref{laplace:ordinary} and the beginnings of Sections \ref{sec:geometry_borel} and~\ref{sec:borel-laplace-homom}.

In Section~\ref{sec:proof_main_results}, we restate and prove our main results. In Section~\ref{sec:examples}, we illustrate our results by working through detailed examples of how Borel regular functions arise, both as $1$-dimensional thimble integrals and as solutions of level~$1$ ODEs.

In Appendix~\ref{airy-appendix}, reprising our analysis of the Airy--Lucas equations in Section~\ref{sec:examples}, we focus on the Airy equation as a concrete special case. We then explain how our treatment lines up with other discussions of the Airy equation in the literature on Borel summation.

In Appendix~\ref{shifting}, we prove a fact about the relationship between integral and differential equations which is used in some of our examples.

\subsection{Notation}
In reasoning about Borel summation, it is important to keep track of whether we are working with holomorphic functions or formal series, and whether we are working on the position domain (the Borel plane) or the frequency domain (the $z$ plane). We have adopted a variable naming convention that makes these distinctions apparent at a glance.
\begin{notation*}
Functions and series on the frequency and position domains are named as follows.
\begin{center}
\begin{tabular}{c|c|c}
& \textbf{Analytic}: no tilde & \textbf{Formal}: tilde \\[1mm] \hline
\vphantom{\rule{0mm}{5mm}} \textbf{Frequency domain}: upper case & $\Phi$ & $\series{\Phi}$ \\[1mm] \hline
\vphantom{\rule{0mm}{5mm}} \textbf{Position domain}: lower case & $\phi$ & $\series{\phi}$ \\[1mm]
\end{tabular}
\end{center}
We will sometimes use a hat to emphasize that we got a holomorphic function $\hat{\phi}$ by taking the sum of formal series $\series{\phi}$.

The frequency and position coordinates $z$ and $\zeta$ are excepted from this convention. The frequency coordinates is lower-case, even though it lives on the frequency domain, and neither variable gets a tilde when it serves as a formal variable rather than a coordinate function.
\end{notation*}

Here are some examples of this notation in action. Since the Laplace transform $\laplace_{\zeta, 0}$ turns functions on the position domain into functions on the frequency domain, we might write ${\laplace_{\zeta,0} \phi = \Phi}$. On the other hand, the Borel transform $\borel_\zeta$ turns formal series in the frequency variable into formal series in the position variable, so we might write $\series{\phi} = \borel_\zeta \series{\Phi}$.

We will follow our convention from~\cite{reg-sing-volterra} for integrals along paths.
\begin{notation*}[{reproduced from \cite[Section~1.2.3]{reg-sing-volterra}}]
When we integrate a holomorphic function along a path, the integrand will always be either a 1-form, like $\varphi {\rm d}\zeta$, or a density, like $|\varphi {\rm d}\zeta|$~\cite[Section~1.8]{local-viewpoint}. On a simply connected domain, all integration paths are equivalent, so we will typically specify an integral by its start and end points. These points may be given directly, as in the expression $\int_b^a \varphi \,{\rm d}\zeta$, or described in terms of coordinates, as in the equivalent expression~${\int_{\zeta = \beta}^a \varphi\, {\rm d}\zeta}$. In the integral
\smash{$\int_{\zeta = \beta}^a (\zeta(a) - \zeta)^{-1/2} \varphi \,{\rm d}\zeta$},
notice that $\zeta(a)$ is a number, which stays constant throughout the integration, while $\zeta$ is a~function, whose value changes along the integration path.

Sometimes, in expressions like
\smash{$ \int_{\zeta(a') = \beta}^a w\bigl(a, a'\bigr) \varphi\bigl(a'\bigr) \,{\rm d}\zeta\bigl(a'\bigr)$},
it is useful to give a name $a'$ to the point moving along the integration path. In this case, all of the functions and differential forms whose values change along the integration path are explicitly evaluated at $a'$.
\end{notation*}
\section{Historical context}\label{sec:historical-context}
\subsection{Borel regularity as a good approximation condition}
Borel regular functions can be characterized as functions that are approximated well, asymptotically, by polynomials. Watson showed a century ago \cite[Part II, Section 9]{watson2} that a function $\Phi$ is Borel regular whenever its asymptotic expansion for large $|z|$ is uniformly $1$-Gevrey-asymptotic in an obtuse-angled sector at infinity (see Definition~\ref{def:unif-gevrey-asymp}).

Watson's theorem was soon improved by Nevanlinna \cite{nevanlinna} (see a modern proof in \cite[Theorem~B.15]{nikolaev2023exact} and a generalization to power series with fractional power exponents in \cite{delabaere--rosoamanana}), and improved again later by Sokal \cite{sokal1980improvement}. These improvements tell us that the obtuse-angled sector around $\infty$ in the statement above can be replaced with an open disk whose boundary touches infinity---that is, an half-plane which does not contain $0$, and may be displaced from $0$ by some distance $\Lambda$.
\begin{center}
\begin{tikzpicture}
\setfiberscale{2.5}
\renewcommand{\phase}{35}
\newcommand{\dis}{1.2}
\begin{fiber}[shade, phase=\phase, dis=\dis]
\genfibercontent[dual]
\draw (0, 0)--(-\phase:\rshade);
\draw[rotate=-\phase] (\dis, -0.1)--(\dis, 0.1);
\draw (0, 0)--(0.75, 0);
\draw[shorten <=0.3mm, shorten >=0.3mm] (0:0.5) arc (0:-\phase:0.5) node[midway, anchor=180-\phase/2] {$-\theta$};
\node[anchor=90-\phase, inner sep=2mm] at (-\phase:\dis) {$\Lambda$};
\node[anchor=north west, inner sep=2mm] at (-\rfiber, \rfiber) {$z$};
\end{fiber}
\end{tikzpicture}
\end{center}
When the half-plane extends along the $-\theta$ direction, the sum $\hat{\phi}$ of the Borel transform of $\aexp^{-\theta} \Phi$ has an absolutely convergent Laplace transform along the $\theta$ direction. In fact, we have $|\hat{\phi}| \lesssim {\rm e}^{\Lambda |\zeta|}$ uniformly over any constant-radius neighborhood of the ray $\zeta \in {\rm e}^{{\rm i}\theta}[0, \infty)$.
\begin{center}
\begin{tikzpicture}
\setfiberscale{2.5}
\pgfmathsetmacro{\rshade}{1.5*\rfiber}
\renewcommand{\phase}{35}
\newcommand{\radius}{0.85}
\begin{fiber}
\fill[pwbeige!30, draw=pwbeige!60, rotate=\phase] (\rshade, \radius)--(0, \radius) arc (90:270:\radius)--(0, -\radius)--(\rshade, -\radius)--cycle;
\genfibercontent
\draw (0, 0)--(0.75, 0);
\draw[ray] (0, 0)--(\phase:\rshade);
\draw[shorten <=0.3mm, shorten >=0.3mm] (0:0.5) arc (0:\phase:0.5) node[midway, anchor=180+\phase/2] {$\theta$};
\node[anchor=north west, inner sep=2mm] at (-\rfiber, \rfiber) {$\zeta$};
\end{fiber}
\end{tikzpicture}
\end{center}

The Watson--Nevanlinna--Sokal characterization of Borel regular functions is totally general, which means it cannot take advantage of any extra structure provided by the problem you are trying to solve. We take the opposite approach, showing that certain functions are Borel regular just because of the extra structure provided by the problems they solve.
\subsection[Solving level 1 ODEs]{Solving level 1 ODEs}\label{sec:history_ODE}
The study of irregular singular differential equations on complex domains has a long history, and interesting phenomena distinguish irregular singular equations from regular singular ones. Analytic solutions of a linear ODE with an irregular singularity exhibit the Stokes phenomenon: asymptotic behaviour that changes sharply from one direction to the next. Formally, this behaviour is captured in a {\em formal integral} solution: an expression
\[ \sum_{\alpha \in A} {\rm e}^{-\alpha z} z^{\tau_\alpha} \series{F}_\alpha \]
in which $A$ is a finite set of complex numbers, each $\tau_\alpha$ is a real number, and each $\series{F}_\alpha$ is a formal power series in $\C\big\llbracket z^{-1} \big\rrbracket$, where $z$ is a coordinate that puts the singularity at $z = \infty$. For each direction $\theta$, the term with the lowest value of $\Re\bigl({\rm e}^{{\rm i}\theta} \alpha\bigr)$ represents the dominant contribution to the solution's asymptotic behavior.

Each term of a formal integral is an example of a {\em trans-monomial}---the basic building block of a {\em trans-series}.\footnote{For the full definition of a trans-series, see~\cite{costin_transseries,costin_summability,van-der-hoeven2001complex,EcalleIII}.} The trans-monomials in the spaces ${\rm e}^{-\alpha z} z^\tau \C\big\llbracket z^{-1} \big\rrbracket$ are the only ones we will need to consider in this paper. Poincar\'e's method for solving linear level~$1$ ODEs, discussed in Section~\ref{sec:Poincare method}, produces trans-monomial solutions in these spaces. The Poincar\'e solutions serve as the terms of a general formal integral. Each Poincar\'e solution is the asymptotic expansion of an analytic solution, whose existence is guaranteed---non-constructively---by the main asymptotic existence theorem (MAET).

The Poincar\'e solutions and the MAET have inspired several other methods for solving linear level~$1$ ODEs. The ones we will discuss in this section are summarized in the following table. We classify solution methods by two features: whether the solutions they produce are analytic functions or formal series, and whether they solve the original equation on the frequency domain or an equivalent equation in the position domain.
\begin{center}
\begin{tabular}{l|l|l}
& \textbf{Analytic} & \textbf{Formal} \\ \hline
\textbf{Frequency} & & Poincar\'e: trans-monomial ansatz~\cite{int-irreg} \\ \hline
\textbf{Position} & \'{E}calle: resurgence & \'{E}calle: formal perturbation theory~\cite{EcalleIII,loday-Remy2011} \\
& Fixed-point iteration~\cite{reg-sing-volterra} \\
\end{tabular}
\end{center}

\subsubsection{The Poincar\'e method and Borel summation}\label{sec:Poincare method}
An ODE of the form described in Section~\ref{borel-reg:explanatory-power} always has a frame of solutions in the trans-monomial spaces ${\rm e}^{-\alpha z} z^{\tau_\alpha} \C\big\llbracket z^{-1} \big\rrbracket$ indexed by the roots $-\alpha$ of $P$, where \smash{$\tau_\alpha := \frac{Q(-\alpha)}{P'(-\alpha)}$}. These solutions can be found systematically, using an algorithm described by Poincar\'e~\cite{int-irreg}, \cite[Proposition~2.2.7, p.~111]{EcalleIII}. Poincar\'e observes that when $z$ is large, the constant-coefficient equation~$P\bigl(\frac{\partial}{\partial z}\bigr) \Phi = 0$ approximates the equation $\mathcal{P}\Phi = 0$ that we are trying to solve. The solutions of this approximate equation are the exponentials $\{{\rm e}^{-\alpha z} \mid P(-\alpha) = 0\}$. Poincar\'e guesses that each approximate solution ${\rm e}^{-\alpha z}$ can be turned into an exact solution ${\rm e}^{-\alpha z} \series{F}_\alpha$ using a formal correction factor $\series{F}_\alpha \in z^{-\tau_\alpha} \C\big\llbracket z^{-1} \big\rrbracket$. The correction factor can be found order by order, starting from any chosen constant term. We will refer to the resulting solutions \[{\rm e}^{-\alpha z} \series{F}_\alpha \in {\rm e}^{-\alpha z} z^{\tau_\alpha} \C\big\llbracket z^{-1} \big\rrbracket\] as Poincar\'e's formal solutions.

Although it is not constructive, the MAET guarantees the existence of a frame of analytic solutions asymptotic to Poincar\'e's formal ones~\cite[Chapter~14]{balser}. Furthermore, the Ramis index theorem tells us that the Borel transform sends Poincar\'e's formal solutions to convergent series on the position domain---the first step toward proving Borel summability~\cite{ramis_index}. This has motivated the development of summability methods, which promote formal solutions to analytic ones~\cite{diverg-resurg--ii,malgrange--fourier,malgrange1995sommation,malgrange92,ramis1991series}. These methods can be applied both within and beyond the world of linear level~$1$ ODEs. For instance, Costin has applied them to a broad class of nonlinear systems of ODEs~\cite{Costin:borel-sum-non-linear}.
\subsubsection{\'{E}calle's theory of resurgence}
\'{E}calle's theory of resurgence introduced a new perspective on formal solutions, mostly focused on the analysis of the position domain~\cite{EcalleIII,loday-Remy2011}. In a nutshell, resurgence provides information about the analytic continuation of the Borel transform of a divergent series. In particular, it reveals the locations of singularities in the position domain, and it quantifies the Stokes phenomenon in terms of ``Stokes constants'' describing analytic continuation around the singularities.

A series is resurgent if its Borel transform converges to an endlessly analytically continuable function on the position domain. Thus, resurgence creates a bridge between formal and analytic solutions of problems in the position domain.

The formal solutions of a level~$1$ ODE are always resurgent series~\cite[Proposition~2.2.1]{EcalleIII}. \'{E}calle proved this by using formal perturbation theory to solve a corresponding integral equation in the position domain. Using the ``formalism of singularities'' to understand how the singularities of the perturbative solution would propagate from one order of perturbation to the next, \'{E}calle showed that the perturbation series would converge to an endlessly analytically continuable function. Loday-Richaud and Remy went on to show that this function always has a well-defined Laplace transform, so the corresponding formal solution in the frequency domain is Borel summable~\cite{loday-Remy2011}. Malgrange got similar results by working with differential equations, rather than integral equations, on the position domain~\cite{malgrange--fourier}.

\'{E}calle's theory not only shed light on Borel summation, but also generalized it. \'{E}calle introduced generalized resummation methods that use ``well-behaved averages'' to handle infinite sequences of singularities along the summation axis~\cite{EcalleI,ecalle-dulac,menous_phdthesis,menous-mould}.
\subsubsection{Solving Riemann--Hilbert problems}
A level~$1$ ODE can be seen as the coordinate expression of a meromorphic connection on a~principal bundle---specifically, the kind of meromorphic connection called an {\em oper}~\cite{BD-opers}. Local solutions of the ODE correspond to local flat sections of the connection. We can therefore use the theory of meromorphic connections to study level~$1$ ODEs. In particular, given an ODE, we can set up a generalized Riemann--Hilbert problem whose solution would encode a frame of analytic solutions of the ODE~\cite[Section~3.2]{Dubrovin-Heun}.

To formulate the Riemann--Hilbert problem, we fix the asymptotic behaviour of the frame of solutions at $z = 0$ and $z = \infty$, and we also fix its discontinuities---the Stokes data prescribing how the frame jumps across each Stokes ray. A solution is then given by a piecewise analytic function which is holomorphic away from the Stokes rays and has the prescribed asymptotics and discontinuities.

In many cases, the solutions of Riemann--Hilbert problems can be found explicitly~\cite{BBS-RH,Tom-RH-1,Tom-RH-conifold,Dubrovin-tt_star,Dubrovin-Heun,GMN1}. These cases include some Riemann--Hilbert problems arising from level~$1$ ODEs. For example, in \cite[Section~6.2.1]{kontsevich2022analyticity}, Kontsevich and Soibelman discuss a Riemann--Hilbert problem for thimble integrals using the formalism of analytic wall-crossing structures. As mentioned in Section~\ref{borel-reg:explanatory-power}, every thimble integral satisfies some ODE. Kontsevich and Soibelman also explain why the asymptotics of these thimble integrals are resurgent series.
\subsubsection{Laplace transform methods and fixed point iteration}
All of the historical approaches we have discussed so far start by looking for formal solutions. In this paper, we will focus on a different kind of approach, where we start by looking for analytic solutions expressed as Laplace transforms. By carefully choosing the starting point of the Laplace transform, we get analytic solutions which are asymptotic to Poincar\'e's formal ones. These solutions will turn out to be Borel regular.

We use the Laplace transform to turn differential equations on the frequency domain into integral equations on the position domain. This runs parallel to \'{E}calle's use of the Borel transform. Building on the existence and uniqueness results we proved in~\cite{reg-sing-volterra}, we use Picard iteration to solve the integral equation in the position domain, and to show that the corresponding solution in the frequency domain is Borel regularizable.

The function spaces mentioned in Section~\ref{sec:laplace_analytic} play an important role in our approach. They also appear in the work of Braaksma, who uses them to solve systems of ODEs whose coefficients are expressed as Laplace transforms~\cite{braaksma2006laplace}. This overlap suggests an opportunity to combine the two approaches, extending our results to systems of equations, and Braaksma's to equations with more general coefficients.
\subsection{Thimble integrals}
Thimble integrals have been studied from different perspectives: in physics, they play an important technical role in quantum mechanics, where infinite-dimensional exponential integrals are supposed to give the expectation values of observable quantities \cite{dunne-unsal2,dunne-unsal,Fauvet_Menous_Queva,Tanizaki:2014tua}. In the algebraic geometric set-up, namely when, following the notation introduced in Section~\ref{borel-reg:explanatory-power}, $X$~is an $N$-dimensional algebraic variety over $\C$ and $f$ is a proper map $f\maps X\to\C$, thimble integrals are known as exponential period integrals~\cite{deligne2007singularites,Maxim_lectures,kontsevich-Fourier-exact,kontsevich2024holomorphic,fresan-notes}.
\subsubsection{Thimble integrals in physics}
In wave optics, a thimble integral can arise when we add up the secondary waves emanating from all the points along a wavefront~\cite{Fenyes-ihes-lecture}. For example, the Airy integral approximates the sum of the secondary waves coming from a wavefront with an inflection point~\cite{rainbows+glories}. The formalism of wave optics can also be applied to the scattering of quantum particles---for example, in the ``nuclear optical model''~\cite{nuclear-mie,nuclear-optical}. Knoll and Schaeffer's use of thimble integrals to describe nuclear scattering eigenstates~\cite{semiclass-scat} inspired Voros to find the energy eigenstates of the quartic oscillator by resumming the asymptotic series of thimble integrals~\cite{quartic-return}.

In the path integral picture of quantum mechanics, which draws inspiration from wave optics, infinite-dimensional thimble integrals are supposed to give the expectation values of observable quantities. In this context, the integral's formal asymptotic expansion is often better-defined than the integral itself, so physicists use Borel summation, resurgence, and related techniques to assign the integral a value. For instance, in the two-dimensional, $\mathcal{N}=2$ supersymmetric QFTs given by Landau--Ginzburg theories, Cecotti and Vafa showed that instantons can be computed using Picard--Lefschetz formulas. Indeed, the so-called wall-crossing phenomenon is equivalent to the Stokes phenomenon for thimble integrals~\cite{Cecotti:1992rm}. Similarly, in complex Chern--Simons theory, the path integral can be studied by decomposing it into thimble integrals and then arguing as in finite dimension~\cite{gukov-marino-purtrov-resurgence,Witten}. Other examples can be found in~\cite{Unsal--resurgence-gauge,dunne-unsal2,Berry1991,Berry_Howls,costin_kruskal,dunne-unsal,Fauvet_Menous_Queva,Garoufalidis--CS,GGM,GTM--CS,Howls,Howls97,pham1988resurgence,Tanizaki:2014tua}.
\subsubsection{Thimble integrals in geometry}
Thimble integrals can be used to explore the geometry of complex manifolds, and even to represent elements of a homology group. To see how this works, we first recall how the one-dimensional definition of a Lefschetz thimble arises from the general definition. Suppose that~${f \maps X \to B}$ is a~holomorphic map from an $N$-dimensional complex manifold $X$ to a translation surface $B$. Consider a translation chart $\zeta$ on $B$ and a non-degenerate critical point $a \in X$ that~$f$ sends to~${\zeta = \alpha}$. Near $a$, we can always find coordinates $t_1, \dots, t_N$ on $X$ with
$ f^*\zeta = \alpha + {\rm e}^{{\rm i} \theta} \bigl(t_1^2 + \cdots + t_N^2\bigr)$.
When the ray $\mathcal{J}_{\zeta, \alpha}^\theta$ avoids the other critical values of $f$, each such coordinate system defines a Lefschetz thimble: the oriented real-analytic submanifold where $t_1, \dots, t_N$ are real. There are many such coordinate systems, but they are all related by local holomorphic maps from $X$ to the complex orthogonal group $O_N(\C)$. Since $O_N(\C)$ has only two connected components, distinguished by the sign of the determinant, a Lefschetz thimble is determined up to homotopy and orientation by the critical point $a$ and the direction $\theta$. When $N = 1$, the homotopy freedom disappears, because $O_1(\C)$ is the discrete group $\{\pm 1\}$.

Lefschetz thimbles for $f$ with direction $\theta$ represent classes in the ``rapid decay homology'' group $H^\theta_N(X,f)$~\cite{pham}, \cite[Section~1.1]{fresan-notes}. Roughly speaking, this group describes homology relative to the region where $\Re\bigl({\rm e}^{-{\rm i}\theta} f\bigr)$ is large and positive. In this context, each thimble integral~$\int_{\mathcal{C}}{\rm e}^{-zf} \nu$
computes the pairing between a rapid decay homology class $\mathcal{C}$ and a ``twisted $1$-form'' ${\rm e}^{-zf} \nu$. When all the critical points of $f$ are non-degenerate, Lefschetz thimbles actually form a basis for $H^\theta_N(X,f)$. As $\theta$ varies, the groups $H^\theta_N(X,f)$ fit together into a local system over the circle of directions, which is singular at the direction of each segment connecting a pair of critical values. The monodromy of this local system can be found using the Picard--Lefschetz formula~\cite[Section~1]{Arnold}, \cite[2\`{e}me Partie, Section~3.3]{pham}.

Like solutions of level~$1$ ODEs, thimble integrals often exhibit the Stokes phenomenon. For one-dimensional thimble integrals, the Stokes constants can be interpreted geometrically as intersection numbers for pairs of thimbles~\cite{kontsevich2022analyticity}. More generally, finding the Stokes constants is equivalent to finding the monodromy of the local system $H^\theta_N(X,f)$ (see the recent development~\cite{kontsevich2024holomorphic}).
\subsubsection{Asymptotics of thimble integrals}
Thimble integrals are traditionally understood through their asymptotics, which can be found with the saddle point approximation~\cite{andersen2020resurgence,delabaere_dillinger_pham,delabaere-howls,Delabaere-Pham99,dingle1973asymptotic,Malgrange22,Pham83}. The asymptotic expansion of a thimble integral is typically divergent, but its Borel sum often matches the integral it came from. In other words, thimble integrals tend to be Borel regular. Theorem~\ref{thm:maxim} explains why this happens in the one-dimensional case.

To prove Theorem~\ref{thm:maxim}, we first shift our focus from the thimble integral to the corresponding analytic object in the position domain. The integral contains a hint about where to find this analytic object: the thimble itself, which is a real-analytic submanifold of a complex manifold. When the thimble is one-dimensional, we can use the well-known formula in Lemma~\ref{lem:thimble_proj_formula} to turn it into a function whose Laplace transform is the thimble integral. This formula generalizes to higher-dimensional thimbles~\cite{pham} and thimbles based at degenerate critical points~\cite[Section~1.2.2]{mistegard_phdthesis}.
\section{The Laplace and Borel transforms}\label{sec:Laplace-Borel-general}
\subsection{The geometry of the Laplace transform}\label{sec:geometry_laplace}
Classically, the Laplace transform turns functions on the position domain into functions on the frequency domain. In the study of Borel summation and resurgence, it is useful to see the position domain as a {\em translation surface} $B$, and the frequency domain as one of its cotangent spaces. Roughly speaking, the Laplace transform lifts holomorphic functions on $B$ to holomorphic functions on $T^*B$.
\subsubsection{Background on translation surfaces}\label{sec:transl}

{\bf A brief definition.}
A translation surface is a Riemann surface $B$ carrying a holomorphic $1$-form $\lambda$~\cite{zorich2006flat}. A {\em translation chart} is a local coordinate $\zeta$ with ${\rm d}\zeta = \lambda$. The standard metric on $\C$ pulls back along translation charts to a flat metric on $B$, with a conical singularity of angle~$2\pi n$ wherever $\lambda$ has a zero of order $n-1 > 0$.

We will call the zeros of $\lambda$ {\em branch points}. To explore the region around a branch point, it can be helpful to use a {\em translation parameter}: a function $\zeta$ which has ${\rm d}\zeta = \lambda$, but is not necessarily a local coordinate.

In all of our examples, $B$ will be a finite-type Riemann surface, and $\lambda$ will have a pole at each puncture. This level of generality allows for plenty of interesting behavior without letting~$B$ get too messy. Sections~2.4--2.5 of \cite{gupta2013meromorphic} give a sense of what $B$ can look like near a pole of $\lambda$.

{\bf A sense of direction.}
The translation structure gives $B$ a notion of direction as well as distance. Away from the branch points, we can talk about moving upward, rightward, or at any angle, just as we would on $\C$. At a branch point of cone angle $2\pi n$, we can also talk about moving upward, rightward, or at any angle in $\R/2\pi\Z$, but here there are $n$ directions that fit each description. To make this more concrete, note that around any point $b \in B$, there is a~unique translation parameter $\zeta_b$ that vanishes at $b$. This parameter is a translation chart when~$b$ is an ordinary point, and an $n$-fold branched covering of $\C$ when $b$ is a branch point of cone angle~$2\pi n$. In either case, $\zeta_b \in {\rm e}^{{\rm i}\theta} [0, \infty)$ is a ray or a set of rays leaving $b$ at angle $\theta \in \R/2\pi\Z$.

Near each branch point $b$, fix a coordinate $\omega_b$ with \smash{$\zeta_b = \frac{1}{n} \omega_b^n$}, where $2\pi n$ is the cone angle at~$b$. This lets us label each direction at $b$ with an ``extended angle'' in $\R/2\pi n\Z$. Of course, there are $n$ different choices for $\omega_b$.

{\bf The frequency coordinate.}\label{transl-freq}
Over the complement $B'$ of the branch points, the translation structure gives us a holomorphic map $z \maps T^*B' \to \C$. This map is an isomorphism on every fiber, trivializing $T^*B$ almost globally. Over a branch point $b$ of cone angle $2\pi n$, we get an analogous isomorphism $z \maps T^*_bB^{\otimes n} \to \C$. In both cases, we will call $z$ the {\em frequency coordinate} of $B$. Both cases are illustrated in Figure~\ref{fig:geometric-picture} below.
\begin{figure}[!ht]
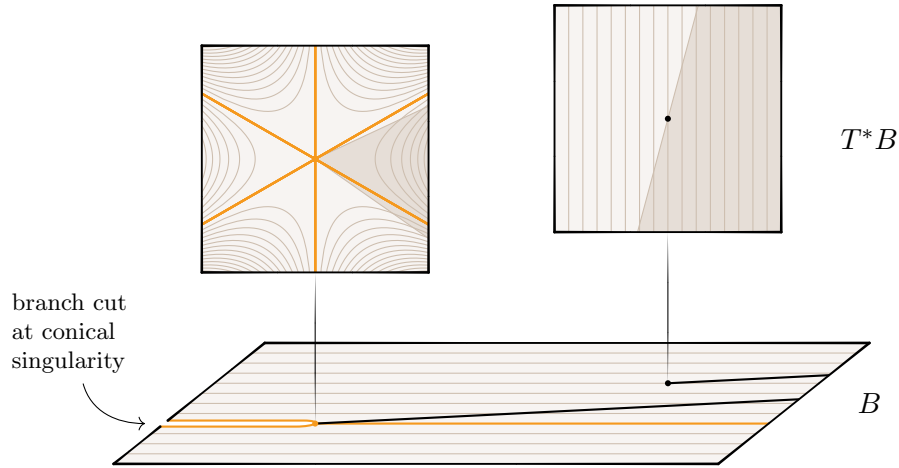
\centering
\phaseSpaceLaplace
\caption{The frequency coordinate $z$ on the cotangent spaces of an ordinary point and a~singular point. The singularity shown here has cone angle $6\pi$, like the singularities of the translation surface associated with the Airy function.}\label{fig:geometric-picture}
\end{figure}

\noindent
At an ordinary point, we can define $z$ simply as the map
\begin{equation*}
z \maps\ T^*_bB \to \C ,\qquad
\lambda\big|_b \mapsto 1.
\end{equation*}
To generalize $z$ to branch points, though, we need a more sophisticated definition. Recall that~${T^*_bB = \van_b / \van_b^2}$, where $\van_b$ is the ideal of holomorphic functions that vanish at $b$. Observing that~${(f + \van_b)^n}$ lies within $f^n + \van_b^{n+1}$ for any $f \in \van_b$, we can identify $T^*_bB^{\otimes n}$ with $\van_b^n / \van_b^{n+1}$ for~${n \ge 1}$. When $b$ is an ordinary point, the translation parameter $\zeta_b$ that vanishes at $b$ represents a nonzero element of $\van_b / \van_b^2$: the cotangent vector $\lambda|_b$. In general, $\zeta_b$ represents a nonzero element of $\van_b^n / \van_b^{n+1}$, where $2\pi n$ is the cone angle at $b$. We define $z$ as the isomorphism
\begin{equation*}
z \maps\ \van_b^n / \van_b^{n+1} \to \C,\qquad
\zeta_b + \van^{n+1} \mapsto 1.
\end{equation*}
When $b$ is a branch point, the coordinate $\omega_b$ we fixed in ``A sense of direction'' gives us an isomorphism
$
w_b \maps T^*_bB \to \C $, $
\omega_b + \van^2 \mapsto 1$
that makes the diagram
\begin{center}
\begin{tikzcd}
T^*_bB^{\otimes n} \arrow[r,"z"] & \C \\
T^*_bB \arrow[u,"\blankbox^n"] \arrow[r,"w"'] & \C \arrow[u,"\blankbox^n"']
\end{tikzcd}
\end{center}
commute. Here, $\blankbox^n$ represents the $n$th-power map.
\subsubsection{The Laplace transform over an ordinary point}\label{laplace:ordinary}
Pick a translation parameter $\zeta$ on $B$ and an extended angle $\theta \in \R$. The {\em Laplace transform}~$\laplace_{\zeta, \alpha}^\theta$ turns a local holomorphic function $\phi$ on $B$ into a local holomorphic function on~$T^*_{\zeta = 0} B$. When~${\zeta = 0}$ is an ordinary point, the Laplace transform is defined by the formula
\begin{equation}\label{laplace:int}
\laplace_{\zeta, \alpha}^\theta \phi = \int_{\mathcal{J}_{\zeta, \alpha}^\theta} {\rm e}^{-z\zeta} \phi \,{\rm d}\zeta,
\end{equation}
where $z$ is the frequency function and $\mathcal{J}_{\zeta, \alpha}^\theta$ is the ray $\zeta \in \alpha + {\rm e}^{{\rm i}\theta} [0, \infty)$. To make sense of this formula, we ask for the following conditions:
\begin{itemize}\itemsep=0pt
\item The starting point $\zeta = \alpha$ is in the domain of $\zeta$. Once we have this, we can continue $\zeta$ along the whole ray $\mathcal{J}_{\zeta, \alpha}^\theta$.
\item The ray $\mathcal{J}_{\zeta, \alpha}^\theta$ avoids the branch points after leaving $\zeta = \alpha$.
\item The integral converges. We ensure this by putting conditions on $\phi$ and $z$.
\begin{itemize}\itemsep=0pt
\item With respect to the flat metric, $\phi$ is uniformly of exponential type $\Lambda$ along the ray~$\mathcal{J}_{\zeta, \alpha}^\theta$, and is locally integrable throughout the ray.\footnote{Recall that a function $\phi$ is of exponential type $\Lambda$ if for every $\varepsilon>0$, there is a constant $A_\varepsilon$ (which may depends on $\varepsilon$) such that $|\phi|\le A_\varepsilon {\rm e}^{(\Lambda+\varepsilon)|\zeta|}$. We instead require a uniform constant $A$ such that $|\phi| \le A {\rm e}^{\Lambda|\zeta|}$.}
\item The value of $z$ satisfies the inequality $\Re\bigl({\rm e}^{{\rm i}\theta} z\bigr) > \Lambda$, which cuts out a half-plane $H^\theta_\Lambda$ in $T^*_{\zeta = 0} B$.
\end{itemize}
\end{itemize}
For any $\sigma > -1$, the conditions on $\phi$ are satisfied by all of the functions in the spaces $\singexp{\sigma}{\Lambda}(\Omega_\alpha)$ introduced in \cite[Definition~2.3]{reg-sing-volterra}. Here, the domain $\Omega_\alpha$ must contain the ray \smash{$\mathcal{J}_{\zeta, \alpha}^\theta$}, and the norm~$\|\cdot\|_{\sigma, \Lambda}$ is taken with respect to $\zeta = \alpha$.

For convenience, we recall the definition of $\singexp{\sigma}{\Lambda}(\Omega)$. Let $\holo(\Omega)$ denote the space of holomorphic functions on the domain $\Omega$.
\begin{Definition}\label{def:HL_inf}
For any $\sigma, \Lambda \in \R$, let
$\singexp{\sigma}{\Lambda}(\Omega) \subset \holo(\Omega) $
be the subspace consisting of functions~$f$ with $|f| \lesssim |\zeta|^\sigma {\rm e}^{\Lambda|\zeta|}$ over $\Omega$, equipped with the norm
\[ \|f\|_{\sigma,\Lambda} := \sup_\Omega |\zeta|^{-\sigma} {\rm e}^{-\Lambda|\zeta|} |f|. \]
\end{Definition}
\subsubsection{The Laplace transform over a branch point}
When $\zeta = 0$ is a branch point, we can still use formula~\eqref{laplace:int} to define $\laplace_{\zeta, \alpha}^\theta \phi$ on $T_{\zeta = 0}^*B$, as long as we take care of a few subtleties. Thanks to the labelling choices we made in Section~\ref{sec:transl}, the extended angle $\theta \in \R$ still picks out a ray $\mathcal{J}_{\zeta, \alpha}^\theta$. The function $z$ is defined on $T_b^*B^{\otimes n}$, where~$2\pi n$ is cone angle at $\zeta = 0$, so we pull it back to \smash{$T^*_{\zeta = 0} B$} along the $n$th-power map. This amounts to substituting $w_b^n$ for $z$ in formula~\eqref{laplace:int}. The inequality $\Re\bigl({\rm e}^{{\rm i}\theta} z\bigr) > \Lambda$ cuts out a half-plane in~$T_b^*B^{\otimes n}$, which pulls back to $n$ sector-like regions in $T_b^*B$ of angle $\pi/n$. We only define $\laplace_{\zeta, \alpha}^\theta \phi$ on one of them: the one centered around the ray \smash{$w_b \in {\rm e}^{-{\rm i}\theta/n}[0, \infty)$}.
\subsubsection{Change of translation chart}\label{sec:change-translation}
Suppose $\zeta$ is a translation chart on $B$, and $\zeta = \alpha$ is an ordinary point. Let $\zeta_\alpha$ be the translation chart with $\zeta = \alpha + \zeta_\alpha$. The Laplace transforms $\laplace_{\zeta_\alpha, 0}$ and $\laplace_{\zeta, \alpha}$, which both turn functions on $B$ into functions on $T^*_{\zeta = \alpha}B$, are related in the following way.
\begin{Lemma}\label{translation}
If the Laplace transform of $\varphi$ is well defined, then
 \begin{equation}
 \label{change-chart}
 {\rm e}^{-z\alpha} \laplace_{\zeta_\alpha, 0} \varphi = \laplace_{\zeta, \alpha} \varphi.
\end{equation}
\end{Lemma}
\begin{proof}
With a change of variable in the integral that defines the Laplace transform, we see that
\begin{align*}
\laplace_{\zeta, \alpha} \varphi & = \int_{\mathcal{J}_{\zeta,\alpha}} {\rm e}^{-z \zeta} \varphi \,{\rm d}\zeta = \int_{\mathcal{J}_{\zeta_\alpha,0}} {\rm e}^{-z(\alpha + \zeta_\alpha)} \varphi \,{\rm d}\zeta_\alpha = {\rm e}^{-z\alpha } \int_{\mathcal{J}_{\zeta_\alpha,0}} {\rm e}^{-z\zeta_\alpha} \varphi\, {\rm d}\zeta_\alpha \\
& = {\rm e}^{-z\alpha } \laplace_{\zeta_\alpha, 0} \varphi.\tag*{\qed}
\end{align*} \renewcommand{\qed}{}
\end{proof}

We now consider a rescaling of the translation structure of $B$, expanding displacements by a~factor of $\mu \in (0, \infty)$. The coordinate $\xi = \mu\zeta$ is a chart for the new translation structure. The corresponding frequency coordinate $x \maps T^*B \to B$ is given by ${\rm d}\xi \mapsto 1$, so $x = \mu^{-1} z$.
\begin{Lemma}
Let $\varphi\in\singexp{\sigma}{\Lambda}(\Omega)$ with $\sigma>-1$ and for some $\Lambda>0$, then
$\laplace_{\xi, 0} \varphi = \mu \laplace_{\zeta, 0} \varphi$.
\end{Lemma}
\begin{proof}
 From the computation
 \begin{equation*}
\laplace_{\xi, 0} \varphi = \int_{\mathcal{J}_{\xi, 0}} {\rm e}^{-x\xi} \varphi \,{\rm d}\xi = \int_{\mathcal{J}_{\zeta, 0}} {\rm e}^{-z \zeta} \varphi \mu \,{\rm d}\zeta = \mu \laplace_{\zeta, 0} \varphi,
\end{equation*}
we get the desired result.
\end{proof}

\subsection{Analysis of the Laplace transform}\label{sec:laplace_analytic}
\subsubsection{Regularity and decay properties}\label{sec:reg-decay}
\begin{figure}[!ht]\centering
\begin{tikzpicture}
\setfiberscale{2}
\pgfmathsetmacro{\rshade}{1.5*\rfiber}
\newcommand{\dis}{1}
\renewcommand{\phase}{0}
\newcommand{\spread}{30}

\begin{fiber}
\fill[pwbeige!30, draw=pwbeige!60] (0, 0)--(\phase-\spread:1.5*\rshade)--(\phase+\spread:1.5*\rshade)--cycle;
\genfibercontent
\node[anchor=-11, inner sep=2mm] at (0,0) {$\zeta = \alpha$};
\node[anchor=north west, inner sep=2mm] at (-\rfiber, \rfiber) {$\zeta$};
\end{fiber}

\begin{scope}[shift={(5.5, 0)}]
\begin{fiber}
\fill[pwbeige!30, draw=pwbeige!60] (-\phase-\spread:\dis) +(-\phase-\spread-90:\rshade)--+(0, 0) arc (-\phase-\spread:-\phase+\spread:\dis)--++(-\phase+\spread+90:\rshade)-- (\rshade, \rshade)--(\rshade, -\rshade)--(0, -\rshade)--cycle;
\genfibercontent[dual]
\node[anchor=north west, inner sep=2mm] at (-\rfiber, \rfiber) {$z$};
\end{fiber}
\end{scope}
\end{tikzpicture}
\caption{A sector $\Omega_\alpha$ in the position domain, and the corresponding union of half-planes $\widehat{\Omega}_\alpha^\Lambda$ in the frequency domain.}\label{fig:sectorial_domain-pos-fre}
\end{figure}
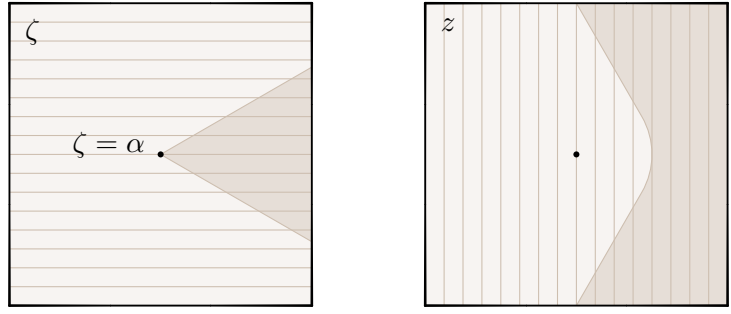

Suppose $\Omega_\alpha$ is an open sector with $\zeta = \alpha$ at its tip, and an opening angle of $\pi$ or less. For any~${\Lambda \in \R}$, let $\widehat{\Omega}_\alpha^\Lambda$ be the union of the half-planes $\Re\bigl({\rm e}^{{\rm i}\theta} z\bigr) > \Lambda$ over all angles $\theta$ in the opening of~$\Omega_\alpha$, as illustrated in Figure~\ref{fig:sectorial_domain-pos-fre}. The domain \smash{$\widehat{\Omega}_\alpha^\Lambda$} is defined wherever $z$ is well defined: in~$T_b^*B$ over each ordinary point $b \in B$, and in $T_b^*B^{\otimes n}$ over each branch point $b \in B$ with cone angle~$2\pi n$.
To~describe the image of the Laplace transform $\laplace_{\zeta,\alpha}$, we introduce frequency-domain counterparts for the function spaces $\singexp{\sigma}{\Lambda}(\Omega_\alpha)$ from Definition~\ref{def:HL_inf}.

\begin{Definition}\label{def:H_hat}
Let \smash{$\dualsingexp{\sigma}\bigl(\widehat{\Omega}_\alpha^\Lambda\bigr)$} be the space of holomorphic functions $\Phi$ on $\widehat{\Omega}_\alpha^\Lambda$ with \smash{$|\Phi| \lesssim \Delta^\sigma$}, where $\Delta$ is the function that measures distance to the boundary of \smash{$\widehat{\Omega}_\alpha^\Lambda$}. The norm $\|\Phi\|_{\sigma, \Lambda} =\allowbreak\smash{ \sup_{\widehat{\Omega}_\alpha^\Lambda} \Delta^{-\sigma} |\Phi|}$ turns \smash{$\dualsingexp{\sigma}\bigl(\widehat{\Omega}_\alpha^\Lambda\bigr)$} into a Banach space.
\end{Definition}
To make statements about \smash{$\dualsingexp{\sigma}\bigl(\widehat{\Omega}_\alpha^\Lambda\bigr)$} that hold when $\Lambda$ is large enough, it is helpful to have a space that encompasses \smash{$\dualsingexp{\sigma}\bigl(\widehat{\Omega}_\alpha^\Lambda\bigr)$} for all $\Lambda$. Since restriction gives an inclusion \[
\dualsingexp{\sigma}\bigl(\widehat{\Omega}_\alpha^\Lambda\bigr) \hookrightarrow \dualsingexp{\sigma}\bigl(\widehat{\Omega}_\alpha^{\Lambda'}\bigr)
\]
 whenever $\Lambda < \Lambda'$, we can define a colimit space \smash{$\dualsingexp{\sigma}\bigl(\widehat{\Omega}_\alpha^\bullet\bigr)$} that all of the spaces~\smash{$\dualsingexp{\sigma}\bigl(\widehat{\Omega}_\alpha^\Lambda\bigr)$} include into.

\begin{Proposition}[following~\cite{sternin1995borel}]\label{prop:laplace-cont}
Let $\Omega_\alpha$ be an open sector with $\zeta = \alpha$ at its tip, and an opening angle of $\pi$ or less. For any $\sigma > 0$ and $\Lambda \ge 0$, and any angle $\theta$ in the opening of $\Omega_\alpha$, the Laplace transform \smash{$\laplace_{\zeta_\alpha, 0}^\theta$} is a continuous map \smash{$\singexp{\sigma-1}{\Lambda}(\Omega_\alpha) \to \dualsingexp{-\sigma}\bigl(\widehat{\Omega}_\alpha^\Lambda\bigr)$}, with a norm of at most $\Gamma(\sigma)$.
\end{Proposition}
\begin{proof}
Given some $\phi \in \singexp{\sigma-1}{\Lambda}(\Omega_\alpha)$, we compute
\begin{align*}
\big|\laplace_{\zeta_\alpha, 0}^\theta \phi\big| & = \left| \int_{\mathcal{J}_{\zeta_\alpha, 0}^\theta} {\rm e}^{-z\zeta_\alpha} \phi\, {\rm d}\zeta_\alpha \right| \le \int_{\mathcal{J}_{\zeta_\alpha,0}^\theta} {\rm e}^{-\Re(z\zeta_\alpha)} |\zeta_\alpha|^{\sigma-1} {\rm e}^{\Lambda |\zeta_\alpha|} \|\phi\|_{\sigma-1, \Lambda} \,|{\rm d}\zeta_\alpha| \\
& \le \int_{\mathcal{J}_{\zeta_\alpha,0}^\theta} {\rm e}^{(\Lambda - c_{z, \theta}|z|)|\zeta_\alpha|} |\zeta_\alpha|^{\sigma-1} \|\phi\|_{\sigma-1, \Lambda} \,|{\rm d}\zeta_\alpha|,
\end{align*}
where $c_{z, \theta} = \cos(\arg(z) + \theta)$. When $|z|$ is large and $c_{z, \theta}$ is positive, the integrand shrinks exponentially as $|\zeta_\alpha|$ grows. This shows that for each angle $\theta$ in the opening of $\Omega_\alpha$, the integral defining \smash{$\laplace_{\zeta_\alpha, 0}^\theta \phi$} converges in some region of \smash{$\widehat{\Omega}_\alpha^\Lambda$}. It also shows that for different angles $\theta$, the functions $\laplace_{\zeta_\alpha, 0}^\theta \phi$ match where their domains overlap. We can therefore glue these functions together into one big Laplace transform of $\phi$, defined at large values of $|z|$ across the whole opening angle of $\widehat{\Omega}_\alpha^\Lambda$.

We can now simplify the calculation of \smash{$\laplace_{\zeta_\alpha, 0}^\theta \phi$} by looking at~$\arg(z)$ and using the closest angle~$\theta$ in the opening of $\Omega_\alpha$. This keeps $\Lambda - c_{z, \theta}|z|$ equal to $\Delta$, the distance to the boundary of~\smash{$\widehat{\Omega}_\alpha^\Lambda$}. It follows that
\begin{equation*}
|\laplace_{\zeta_\alpha, 0}^\theta \phi| \le \int_{\mathcal{J}_{\zeta_\alpha, 0}^{\arg(z)}} {\rm e}^{-\Delta|\zeta_\alpha|} |\zeta_\alpha|^{\sigma-1} \|\phi\|_{\sigma-1, \Lambda}\, |{\rm d}\zeta_\alpha| = \int_0^\infty {\rm e}^{-\Delta t} t^{\sigma-1} \|\phi\|_{\sigma-1, \Lambda} \,{\rm d}t.
\end{equation*}
The integral on the last line converges throughout $\widehat{\Omega}_\alpha^\Lambda$. We can evaluate it by observing that it is a Laplace transform integral, with the roles of position and frequency played by $t$ and $\Delta$ respectively
\smash{$\big|\laplace_{\zeta_\alpha, 0}^\theta \phi\big| \le \Gamma(\sigma) \Delta^{-\sigma} \|\phi\|_{\sigma-1, \Lambda}$}.
In terms of the metric on $\dualsingexp{-\sigma}\bigl(\widehat{\Omega}_\alpha^\Lambda\bigr)$ defined above, this bound says that
\smash{$\big\|\laplace_{\zeta_\alpha, 0}^\theta \phi\big\|_{-\sigma, \Lambda} \le \Gamma(\sigma) \|\phi\|_{\sigma-1, \Lambda}$},
which is what we wanted to show.
\end{proof}

\begin{Proposition}\label{prop:inverse_laplace_analytic}
Let $\Omega_\alpha$ be an open sector with tip $\zeta = \alpha$ and opening angle $\Theta \le \pi$. Choose some $\sigma > 0$ and $\Lambda \ge 0$, and let $\Lambda' = \Lambda / \cos(\Theta)$. Under the conditions of Proposition~{\rm\ref{prop:laplace-cont}}, for any $\varepsilon \in \big(0, \frac{\Theta}{2}\big)$, the Laplace transform
\[ \laplace_{\zeta_\alpha, 0}^\theta \maps\ \singexp{\sigma-1}{\Lambda}(\Omega_\alpha) \to \dualsingexp{-\sigma}\bigl(\widehat{\Omega}_\alpha^\Lambda\bigr) \]
has a continuous left inverse
\[ \bigl(\laplace_{\zeta_\alpha, 0}^\theta\bigr)^{-1} \maps\ \dualsingexp{-\sigma}\bigl(\widehat{\Omega}_\alpha^\Lambda\bigr) \to \singexp{\sigma-1}{\Lambda'}(\Omega_\alpha^\varepsilon), \]
with a norm of at most
\smash{$ \frac{\Gamma(1-\sigma)}{\pi \sin(\varepsilon/2)}$}.
Here, $\Omega_\alpha^\varepsilon \subset \Omega_\alpha$ is the open sector created by cutting a sector of angle $\varepsilon$ off each edge of $\Omega_\alpha$.
\end{Proposition}

\begin{proof}
The boundary of $\widehat{\Omega}_\alpha^\Lambda$ consists of two rays connected by a circular arc. Extend the rays until they meet, forming the boundary of a sector with opening angle $\pi - \Theta$. Expand that sector by angle $\varepsilon/2$ on both sides. The edges of the expanded sector form a path $\mathcal{C}$ within $\widehat{\Omega}_\alpha^\Lambda$, as shown in Figure~\ref{fig:ref-path}.
\begin{figure}[!ht]\centering
\begin{tikzpicture}
\setfiberscale{2.5}
\pgfmathsetmacro{\rshade}{1.5*\rfiber}
\newcommand{\dis}{1.25}
\renewcommand{\phase}{0}
\newcommand{\spread}{55.5}
\pgfmathsetmacro{\peakdis}{\dis / cos(\spread)}

\newcommand{\pathspread}{35}
\newcommand{\angmarkrad}{1}
\newcommand{\distmarkrad}{2.5}
\pgfmathsetmacro{\bdrydist}{\distmarkrad*sin(\pathspread-\spread)}

\begin{fiber}
\fill[pwbeige!30, draw=pwbeige!60, thick] (-\phase-\spread:\dis) +(-\phase-\spread-90:\rshade)--+(0, 0) arc (-\phase-\spread:-\phase+\spread:\dis)-- ++(-\phase+\spread+90:\rshade)--(\rshade, \rshade)--(\rshade, -\rshade)--(0, -\rshade)--cycle;

\genfibercontent[dual]

\begin{scope}[pwbeige!80!black]
\draw (\peakdis, 0) +(-\phase-\spread-90:1.5*\rshade)--+(0, 0)--+(-\phase+\spread+90:1.5*\rshade);
\draw[shorten <=0.3mm, shorten >=0.4mm] (\peakdis, 0) ++(\phase+\spread+90:\angmarkrad) arc[start angle=\phase+\spread+90, delta angle=\pathspread-\spread, radius=\angmarkrad] node [midway, circle, inner sep=1, anchor=-90+0.5*\spread+0.5*\pathspread)] {$\frac{\varepsilon}{2}$};
\draw[shorten >=0.3mm] (\peakdis, 0) ++(\phase+\pathspread+90:\distmarkrad)--++(\phase+\spread:\bdrydist) node[midway, circle, inner sep=0, anchor=-90+\spread] {$t \sin(\frac{\varepsilon}{2})$};
\end{scope}

\draw[very thick] (\peakdis, 0) +(-\phase-\pathspread-90:\rshade)--+(0, 0)--+(-\phase+\pathspread+90:\rshade);
\draw (\peakdis, 0) ++(\phase+\pathspread+90:0.5*\distmarkrad) node[anchor=180+\pathspread] {$t$};
\draw (\peakdis, 0) ++(\phase+\pathspread+90:\distmarkrad) +(\pathspread:1mm)--+(180+\pathspread:1mm);

\node[anchor=north west, inner sep=2mm] at (-\rfiber, \rfiber) {$z$};
\end{fiber}
\end{tikzpicture}
\caption{The thick black line is the path $\mathcal{C}$.}\label{fig:ref-path}
\end{figure}
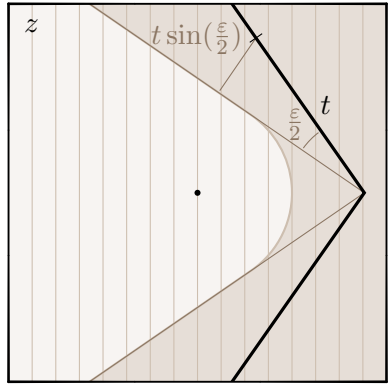

\noindent
The path $\mathcal{C}$ is made of two rays, which meet at an angle of $\pi - \Theta + \varepsilon$. Parameterize it using the arc length parameter $t$ which is zero where the rays meet. Along $\mathcal{C}$, the distance $\Delta$ to the boundary of $\widehat{\Omega}_\alpha^\Lambda$ satisfies the estimate
$\Delta \ge \sin(\varepsilon/2) |t|$,
as shown in Figure~\ref{fig:ref-path}. On $\Omega_\alpha^\varepsilon \times \mathcal{C}$, we have
\begin{equation}\label{eqn:inverse-laplace-exp-bound}
\Re(z\zeta_\alpha) \le |\zeta_\alpha| \bigl(\Lambda' - \sin(\varepsilon/2) |t|\bigr).
\end{equation}
To see why this bound holds, fix a value of the coordinate $\zeta_\alpha$ representing a point in $\Omega_\alpha^\varepsilon$. Consider the half-plane $\Re\bigl(z \frac{\zeta_\alpha}{|\zeta_\alpha|}\bigr) > \Lambda'$, illustrated in Figure~\ref{fig:bound-real-part} for an extreme value of $\zeta_\alpha$.
\begin{figure}[!ht]\centering
\begin{tikzpicture}
\setfiberscale{2.5}
\pgfmathsetmacro{\rshade}{1.5*\rfiber}
\newcommand{\dis}{1.25}
\renewcommand{\phase}{0}
\newcommand{\spread}{55.5}
\pgfmathsetmacro{\peakdis}{\dis / cos(\spread)}

\newcommand{\pathspread}{35}
\newcommand{\angmarkrad}{1}
\newcommand{\distmarkrad}{2.5}
\pgfmathsetmacro{\bdrydist}{\distmarkrad*sin(\pathspread-\spread)}

\pgfmathsetmacro{\sectorshave}{2*(\spread - \pathspread)}
\pgfmathsetmacro{\shavedspread}{\spread - \sectorshave}

\begin{fiber}
\fill[pwbeige!30, draw=pwbeige!60] (0, 0)--(\phase-\shavedspread:1.5*\rshade)--(\phase+\shavedspread:1.5*\rshade)--cycle;

\begin{scope}[pwbeige!80!black]
\draw (\phase-\spread:1.5*\rshade)--(0, 0)--(\phase+\spread:1.5*\rshade);
\draw[shorten <=0.3mm, shorten >=0.3mm] (\phase+\shavedspread:\angmarkrad) arc[start angle=\phase+\shavedspread, delta angle=\sectorshave, radius=\angmarkrad] node [midway, circle, inner sep=1, anchor=-180+0.5*\spread+0.5*\pathspread)] {$\varepsilon$};
\draw[shorten <=0.3mm, shorten >=0.3mm] (\phase-\shavedspread:\angmarkrad) arc[start angle=\phase-\shavedspread, delta angle=-\sectorshave, radius=\angmarkrad] node [midway, circle, inner sep=1, anchor=-180-0.5*\spread-0.5*\pathspread)] {$\varepsilon$};
\end{scope}

\fill (\phase-\shavedspread:1.7) circle (1.25*\dotsize);

\genfibercontent
\node[anchor=-11, inner sep=2mm] at (0,0) {$\zeta = \alpha$};
\node[anchor=north west, inner sep=2mm] at (-\rfiber, \rfiber) {$\zeta$};
\end{fiber}

\begin{scope}[shift={(6.5, 0)}] 

\begin{fiber}[shade, phase=\phase+2*\pathspread-\spread, dis=\peakdis]

\genfibercontent[dual]

\begin{scope}[pwbeige!80!black]
\draw (\peakdis, 0) +(-\phase-\spread-90:1.5*\rshade)--+(0, 0)--+(-\phase+\spread+90:1.5*\rshade);
\draw (\peakdis, 0)--+(\phase+2*\pathspread-\spread+90:1.5*\rshade);
\draw[shorten <=0.3mm, shorten >=0.4mm] (\peakdis, 0) ++(\phase+\spread+90:\angmarkrad) arc[start angle=\phase+\spread+90, delta angle=\pathspread-\spread, radius=\angmarkrad] node [midway, circle, inner sep=1, anchor=-90+0.5*\spread+0.5*\pathspread)] {$\frac{\varepsilon}{2}$};
\draw[shorten <=0.4mm, shorten >=0.3mm] (\peakdis, 0) ++(\phase+\pathspread+90:\angmarkrad) arc[start angle=\phase+\pathspread+90, delta angle=\pathspread-\spread, radius=\angmarkrad] node [midway, circle, inner sep=1, anchor=-90+1.5*\pathspread-0.5*\spread)] {$\frac{\varepsilon}{2}$};
\draw[shorten >=0.3mm] (\peakdis, 0) ++(\phase+\pathspread+90:\distmarkrad)--++(\phase+2*\pathspread-\spread-180:\bdrydist) coordinate[midway] (dist-anchor);
\end{scope}

\draw[very thick] (\peakdis, 0) +(-\phase-\pathspread-90:\rshade)--+(0, 0)--+(-\phase+\pathspread+90:\rshade);
\draw (\peakdis, 0) ++(\phase+\pathspread+90:0.7*\distmarkrad) node[circle, inner sep=0.5, fill=pwbeige!10] {$t$};
\draw (\peakdis, 0) ++(\phase+\pathspread+90:\distmarkrad) +(\pathspread:1mm)--+(180+\pathspread:1mm);

\node[anchor=north west, inner sep=2mm] at (-\rfiber, \rfiber) {$z$};
\end{fiber}

\begin{scope}[pwbeige!80!black]
\draw[out=-55, in=90] (1.7, 3.2) node[anchor=-35, inner sep=1] {$\Re\bigl(z\frac{\zeta_\alpha}{|\zeta_\alpha|}\bigr) > \Lambda'$} to (2.1, 1.9);
\draw[out=-15, in=105, shorten >=0.5mm, -stealth] (0.6, 2.9) node[anchor=east, inner sep=2] {$t \sin(\frac{\varepsilon}{2})$} to (dist-anchor);
\end{scope}

\end{scope} 
\end{tikzpicture}
\caption{For this figure, we have chosen a value of $\zeta_\alpha$ representing a point on the boundary of the sector $\Omega_\alpha^\varepsilon$. This is a worst-case scenario for inequality~\eqref{eqn:inverse-laplace-exp-bound}, since our assumptions require~$\zeta_\alpha$ to represent a point inside $\Omega_\alpha^\varepsilon$.}\label{fig:bound-real-part}
\end{figure}
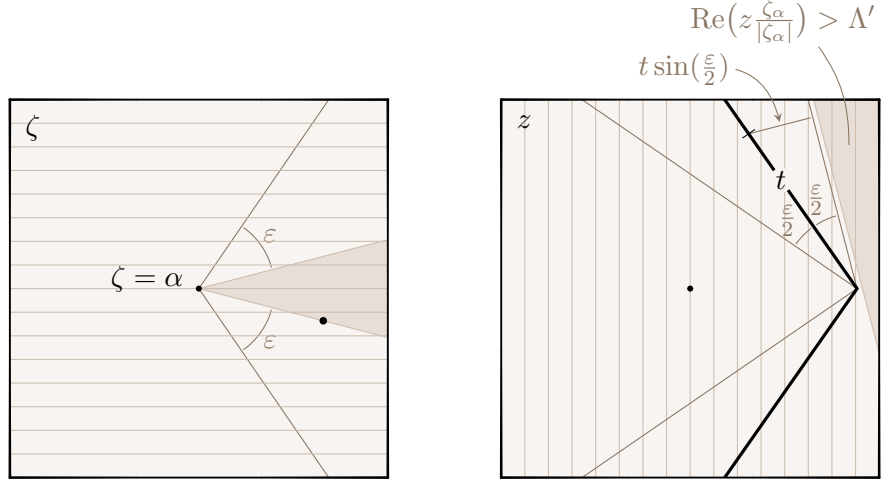
By construction, the curve $\mathcal{C}$ lies entirely outside this half-plane, which means that $\smash{\Re\bigl(z \frac{\zeta_\alpha}{|\zeta_\alpha|}\bigr) }\!\le \Lambda'$ throughout~$\mathcal{C}$. In fact, each of the rays that make up $\mathcal{C}$ makes an angle of at least $\varepsilon/2$ with the boundary line \smash{$\Re\bigl(z \frac{\zeta_\alpha}{|\zeta_\alpha|}\bigr) = \Lambda'$}. Therefore, as we move outward along~$\mathcal{C}$, our distance to the boundary of the half-plane \smash{$\Re\bigl(z \frac{\zeta_\alpha}{|\zeta_\alpha|}\bigr) > \Lambda'$} increases at a rate of at least $\sin(\frac{\varepsilon}{2})$ with respect to the arc length parameter $t$.

The inverse Laplace transform is given by the formula
\[ \bigl(\laplace_{\zeta_\alpha, 0}^\theta\bigr)^{-1} \Phi = \frac{1}{2 \pi {\rm i}} \int_{\mathcal{C}} {\rm e}^{z\zeta_\alpha} \Phi\,{\rm d}z. \]
When $\Phi$ is in $\dualsingexp{-\sigma}\bigl(\widehat{\Omega}_\alpha^\Lambda\bigr)$, we have the bound
\begin{align*}
\big|\bigl(\laplace_{\zeta_\alpha, 0}^\theta\bigr)^{-1} \Phi\big| & \le \frac{1}{2 \pi} \int_{\mathcal{C}} {\rm e}^{\Re(z\zeta_\alpha)} \Delta^{-\sigma} \|\Phi\|_{\sigma, \Lambda}\,{\rm d}z \\
& \le \frac{1}{2 \pi} \int_{-\infty}^\infty {\rm e}^{|\zeta_\alpha| \left(\Lambda' - \sin(\varepsilon/2) |t|\right)} (\sin(\varepsilon/2) |t|)^{-\sigma} \|\Phi\|_{\sigma, \Lambda} \,{\rm d}t \\
& = {\rm e}^{|\zeta_\alpha| \Lambda'} \|\Phi\|_{\sigma, \Lambda} \frac{1}{2 \pi} \int_{-\infty}^\infty {\rm e}^{-|\zeta_\alpha| \sin(\varepsilon/2) |t|} (\sin(\varepsilon/2) |t|)^{-\sigma}\, {\rm d}t,
\end{align*}
which we can rewrite as
\begin{align*}
\big|\bigl(\laplace_{\zeta_\alpha, 0}^\theta\bigr)^{-1} \Phi\big| & \le {\rm e}^{|\zeta_\alpha| \Lambda'} \|\Phi\|_{\sigma, \Lambda} \frac{1}{2 \pi} \int_{-\infty}^\infty {\rm e}^{-|\zeta_\alpha| |s|} |s|^{-\sigma} \frac{{\rm d}s}{\sin(\varepsilon/2)} \\
& = {\rm e}^{|\zeta_\alpha| \Lambda'} \|\Phi\|_{\sigma, \Lambda} \frac{1}{\pi \sin(\varepsilon/2)} \int_0^\infty {\rm e}^{-|\zeta_\alpha|s} s^{-\sigma}\,{\rm d}s
\end{align*}
using the new parameter $s = \sin(\varepsilon/2) t$. Recognizing the integral in the last line as a Laplace transform integral, with the roles of position and frequency played by $s$ and $|\zeta_\alpha|$ respectively, we have
\[ \big|\bigl(\laplace_{\zeta_\alpha, 0}^\theta\bigr)^{-1} \Phi\big| \le {\rm e}^{|\zeta_\alpha| \lambda'} \|\Phi\|_{\sigma, \Lambda} \frac{\Gamma(1-\sigma)}{\pi \sin(\varepsilon/2)} |\zeta_\alpha|^{\sigma-1}, \]
which is what we wanted to show.
\end{proof}

\subsection{The geometry of the Borel transform}\label{sec:geometry_borel}
The Laplace transform $\laplace_{\zeta, 0}$ acts in an especially simple way on powers of the coordinate $\zeta$
\[ \laplace_{\zeta, 0}\left[\frac{\zeta^n}{n!}\right] = z^{-n-1}. \]
Here, we are thinking of $z$ as the frequency coordinate on $T^*_{\zeta = 0}B$, as described in Section~\ref{transl-freq}. We can get a function on $T^*_{\zeta = \alpha}B$ instead by taking the Laplace transform with respect to the coordinate $\zeta_\alpha$ defined by $\zeta = \zeta_\alpha + \alpha$
\[ \laplace_{\zeta_\alpha, 0}\left[\frac{\zeta_\alpha^n}{n!}\right] = z^{-n-1}. \]
On each cotangent space, we can define a formal inverse of the Laplace transform by turning negative powers of $z$ back into powers of the appropriate translation coordinate. This formal inverse is called the {\em Borel transform}. To be more precise, the Borel transform $\borel_\zeta$ on $T^*_{\zeta = 0}B$ is the inverse of $\laplace_{\zeta,0}$ on monomials
\begin{center}
\begin{tikzcd}[every arrow/.append style={shift left}]
 \big\{z^{-1}, z^{-2}, z^{-3}, z^{-4}, \dots \big\} \arrow{d}{\borel_{\zeta}} \\ \big\{1, \zeta, \frac{1}{2!} \zeta^2, \frac{1}{3!} \zeta^3, \dots\big\} \arrow{u}{\laplace_{\zeta, 0}}
\end{tikzcd}
\end{center}
and it extends to formal power series by countable linearity
\begin{equation*}
\borel_\zeta \maps\ z^{-1} \C \big\llbracket z^{-1} \big\rrbracket \to \C \llbracket \zeta \rrbracket ,\qquad
\sum_{n \ge 0} a_n z^{-n-1} \mapsto \sum_{n \ge 0} a_n \frac{\zeta^n}{n!}.
\end{equation*}
This definition extends straightforwardly to fractional powers of $z$. Observing that
$\laplace_{\zeta,0}[\zeta^\sigma]=\Gamma(\sigma+1)z^{-\sigma-1}$
for every $\sigma \in \R \setminus \Z_{\leq 0}$, we define
\smash{$
\borel_\zeta\bigl[z^{-\sigma-1}\bigr] := \frac{\zeta^{\sigma}}{\Gamma(\sigma+1)}$}.
Then, for any $\sigma \in \R \setminus \Z_{\leq 0}$, we can extend by countable linearity to the space $z^{-\sigma}\C\big\llbracket z^{-1}\big\rrbracket$ of ``fractionally shifted'' formal~power series.

On a different fiber of $T_{\zeta=\alpha}^*B$, the Borel transform $\borel_{\zeta_\alpha}$ in the translated coordinate $\zeta_\alpha$ is defined on monomials as the inverse of the Laplace transform $\laplace_{\zeta_\alpha,0}$
\begin{center}
\begin{tikzcd}[every arrow/.append style={shift left}]
 \big\{z^{-1}, z^{-2}, z^{-3}, z^{-4}, \dots \big\} \arrow{d}{\borel_{\zeta_\alpha}} \\ \big\{1, \zeta_\alpha, \frac{1}{2!} \zeta_\alpha^2, \frac{1}{3!} \zeta_\alpha^3, \dots\big\} \arrow{u}{\laplace_{\zeta_\alpha, 0}}
\end{tikzcd}
\end{center}
It extends countable linearity to a map $\borel_{\zeta_\alpha} \maps z^{-1} \C \big\llbracket z^{-1} \big\rrbracket \to \C \llbracket \zeta_\alpha \rrbracket$.
\subsubsection{Action on trans-monomials}\label{sec:action_transseries}
We can extend the definition of the Borel transform to these trans-monomials by recognizing that $\laplace_{\zeta, \alpha}$ sends ordinary monomials to trans-monomials, and then taking advantage of the relationship between $\laplace_{\zeta, \alpha}$ and $\laplace_{\zeta_\alpha, 0}$ given by identity~\eqref{change-chart}.
\begin{Definition}
On each ordinary monomial within a trans-monomial, $\borel_\zeta$ acts as the formal inverse of $\laplace_{\zeta,\alpha}$,
\[
\laplace_{\zeta,\alpha}\borel_\zeta\bigl[{\rm e}^{-z\alpha} z^{-n-1} \bigr]={\rm e}^{-z\alpha} z^{-n-1}.
\]
The action of $\borel_\zeta$ extends to all of ${\rm e}^{-z\alpha} \C\big\llbracket z^{-1} \big\rrbracket$ by countable linearity.
\end{Definition}
From this definition, and identity~\eqref{change-chart}, we deduce that
\begin{gather*}
{\rm e}^{-z\alpha}\laplace_{\zeta_\alpha,0}\borel_\zeta\bigl[{\rm e}^{-z\alpha} z^{-n-1}\bigr] = {\rm e}^{-z\alpha} z^{-n-1} ,\qquad
\laplace_{\zeta_\alpha,0}\borel_\zeta\bigl[{\rm e}^{-z\alpha} z^{-n-1}\bigr] = z^{-n-1} ,\\
\borel_{\zeta}\bigl[{\rm e}^{-z\alpha}z^{-n-1}\bigr] = \frac{\zeta_\alpha^n}{n!}.
\end{gather*}
In other words, the diagram
\begin{center}
\begin{tikzcd}[every arrow/.append style={shift left}]
\big\{z^{-1}, z^{-2}, z^{-3}, z^{-4}, \dots \big\} \arrow[dd,"\borel_\zeta"'] \arrow[rr,"{\rm e}^{-z\alpha}"]& & {\rm e}^{-z\alpha} \big\{z^{-1}, z^{-2}, z^{-3}, z^{-4}, \dots \big\} \arrow[dd, "\borel_\zeta"] \\
& & \\
\big\{1, \zeta, \zeta^2, \zeta^3, \dots\big\} \arrow[rr, "\mathsf{T}_{-\alpha}^*"'] & & \big\{1, \zeta_\alpha, \zeta_\alpha^2, \zeta_\alpha^3, \dots\big\}\arrow[lluu,"\laplace_{\zeta_\alpha, 0}" description]\arrow{uu}{\laplace_{\zeta, \alpha}}
\end{tikzcd}
\end{center}
commutes, where $\mathsf{T}_{-\alpha}$ denotes translation by $-\alpha$. Notice that $\laplace_{\zeta, 0}$ and $\laplace_{\zeta, \alpha}$ produce functions on the fiber $T^*_{\zeta = 0}B$, while $\laplace_{\zeta_\alpha, 0}$ produces functions on the fiber $T^*_{\zeta = \alpha}B$, so this argument depends on the fact that $z$ is defined almost globally $T^*B$.

We can extend the Borel transform to the fractionally shifted trans-monomial spaces ${\rm e}^{-z\alpha} z^\tau\allowbreak\times \C\big\llbracket z^{-1} \big\rrbracket$, where $\tau \in \R \setminus \Z$, by applying the same argument to the definition at the end of Section~\ref{sec:geometry_borel}. We conclude, in particular, that
\smash{$\borel_{\zeta}\bigl[{\rm e}^{-z\alpha}z^{-\tau-1}\bigr] = \frac{\zeta_\alpha^\tau}{\Gamma(\tau+1)}$}
for any $\tau \in \R\setminus\Z$.

\subsubsection{Change of translation chart}\label{transl-borel}
We will show that the Borel transform is compatible with the change of translation chart for the Laplace transform in Section~\ref{sec:change-translation}. To be more precise, we will show that the diagram{\samepage
\begin{center}
\begin{tikzcd}
{\rm e}^{-z\alpha}\big\{z^{-1}, z^{-2}, z^{-3}, z^{-4}, \dots \big\} \arrow[r,"{\rm e}^{z\alpha}"] & \big\{z^{-1}, z^{-2}, z^{-3}, z^{-4}, \dots \big\}\arrow[d, "\borel_{\zeta_\alpha}"] \\
\big\{1, \zeta, \zeta^2, \zeta^3, \dots\big\}\arrow[u,"\laplace_{\zeta,\alpha}"]\arrow[r,"\text{change chart}"'] & \big\{1, \zeta_\alpha, \zeta_\alpha^2, \zeta_\alpha^3, \dots\big\}
\end{tikzcd}
\end{center}
commutes, where the functions of $z$ are functions on the fiber \smash{$T^*_{\zeta=\alpha}B$}.}

\begin{proof}
We want to show that
$\borel_{\zeta_\alpha}[{\rm e}^{z\alpha} \laplace_{\zeta,\alpha}[\zeta^n]]=\zeta^n$.
Recall that $\borel_{\zeta_\alpha}$
is the formal inverse of $\laplace_{\zeta_\alpha,0}$ on the cotangent fibre over $\alpha$. Thus, taking the Borel transform on both side of the identity~\eqref{change-chart} we find
 \begin{align*}
 \borel_{\zeta_\alpha}[{\rm e}^{z\alpha} \laplace_{\zeta,\alpha}[\zeta^n]]&=\borel_{\zeta_\alpha}\laplace_{\zeta_\alpha,0}\bigl[\zeta^n\bigr]=\borel_{\zeta_\alpha}\laplace_{\zeta_\alpha,0}\left[\sum_{k=0}^n{n\choose k}\zeta_\alpha^k \alpha^{n-k}\right]\\
 &=\sum_{k=0}^n{n\choose k} \alpha^{n-k}\borel_{\zeta_\alpha}\laplace_{\zeta_\alpha,0}\bigl[\zeta_\alpha^k \bigr]=\sum_{k=0}^n{n\choose k} \alpha^{n-k}\zeta_\alpha^k=\zeta^n.\tag*{\qed}
 \end{align*}\renewcommand{\qed}{}
\end{proof}

\subsubsection{Action on translations in the frequency domain}
Translations in the frequency domain do not play any role in our work, but they do appear in other contexts---for example, in Laplace transform methods for difference equations. Since~$\borel_\zeta$ is the formal inverse of $\laplace_{\zeta,0}$, we can use the properties of the Laplace transform to deduce how~$\borel_\zeta$ acts on translations in the frequency domain. It turns out that
$\borel_\zeta \mathsf{T}_{-c}^* \series{\Phi} = {\rm e}^{-c\zeta }\borel_\zeta \series{\Phi}$
for all~${\series{\Phi} \in z^{-1} \C\big\llbracket z^{-1} \big\rrbracket}$.
Indeed, the following diagram
\begin{center}
\begin{tikzcd}[every arrow/.append style={shift left}]
\big\{z^{-1}, z^{-2}, z^{-3}, z^{-4}, \dots \big\} \arrow{dd}{\borel_\zeta}\arrow[rr,"\mathsf{T}_{-c}^*"]& & \big\{(z+c)^{-1}, (z+c)^{-2}, (z+c)^{-3}, (z+c)^{-4}, \dots \big\} \arrow{dd}{\borel_\zeta} \\
& & \\
\big\{1, \zeta, \zeta^2, \zeta^3, \dots\big\} \arrow[rr, "{\rm e}^{-c\zeta}"'] \arrow{uu}{\laplace_{\zeta,0}} & & {\rm e}^{-c\zeta}\big\{1, \zeta, \zeta^2, \zeta^3, \dots\big\}\arrow{uu}{\laplace_{\zeta,0}}
\end{tikzcd}
\end{center}
commutes, where the functions of the variable $z$ belong to the fiber $T^*_{\zeta=0}B$. If $z$ is a function on the fiber over $\zeta=\alpha$, then it is enough to replace $\borel_\zeta$ with $\borel_{\zeta_\alpha}$ and to use $\laplace_{\zeta_\alpha,0}$ as the inverse.
\subsection{The Borel and Laplace transforms as algebra homomorphisms}\label{sec:borel-laplace-homom}
{\bf The Laplace transform as algebra homomorphism.}
The Laplace transform $\laplace_{\zeta,\alpha}$ is defined on the space of holomorphic functions in the position domain which are integrable at $\zeta=\alpha$ and exponentially bounded at infinity. This subspace of holomorphic functions on $B$ has a $\C$-algebra structure, with product given by the \textit{convolution product}.
\begin{Definition}\label{def:convolution}
Choose a translation coordinate $\zeta$ on $B$. Let $\Omega_\alpha \subset B$ be an open set with the point $\zeta = \alpha$ on its boundary, and let $\Omega_\beta \subset B$ be an open set with $\zeta = \beta$ on its boundary. For any $\phi \in \singexpalg{\sigma}(\Omega_\alpha)$ and $\psi \in \singexpalg{\sigma'}(\Omega_\beta)$ with $\sigma,\sigma'\in (-1,\infty)$, the convolution product $\phi \ast_{\zeta, \alpha, \beta} \psi$ is the function on $B$ defined by the integral
\[
\phi \ast_{\zeta, \alpha, \beta} \psi := \int_{\mathcal{S}_{\zeta, \alpha, \beta}} \phi(a) \psi\bigl(a'\bigr)\, {\rm d}\zeta(a),
\]
where $\mathcal{S}_{\zeta, \alpha, \beta}$ is the path $[0, 1] \to B^2$ given by
\[
\zeta_\alpha(a) = t\zeta_{\alpha+\beta}, \qquad
\zeta_\beta\bigl(a'\bigr) = (1-t)\zeta_{\alpha+\beta}
\]
in terms of the standard coordinate $t$ on $[0, 1]$.
\end{Definition}
\begin{notation*}
We denote $\ast_{\zeta,\alpha,\alpha}$ by $\ast_{\zeta,\alpha}$ for short. The definition of $\ast_{\zeta, 0, 0}$ given by Definition~\ref{def:convolution} agrees with the standard definition of convolution \cite[Definition~5.12]{diverg-resurg-i}, so we denote this operation simply by $\ast_\zeta$.
\end{notation*}
\begin{Lemma}
Take $\phi\in\singexp{\sigma}{\Lambda}(\Omega_\alpha)$ and $\psi\in\singexp{\sigma'}{\Lambda'}(\Omega_\beta)$, with $\sigma,\sigma'\in (-1,\infty)$. The convolution~$\phi\ast_{\zeta,\alpha,\beta}\psi$ belongs to $\singexp{\sigma''}{\Lambda''}(\Omega_{\alpha+\beta})$, where $\Lambda'' = \Lambda + \Lambda'$ and $\sigma'' = 1+\sigma+\sigma'$.
\end{Lemma}
\begin{proof}
From Definition~\ref{def:convolution}, we get the bound
\begin{align*}
\left|\phi\ast_{\zeta,\alpha,\beta}\psi\right|&\le\int_{\mathcal{S}_{\zeta,\alpha,\beta}} |\phi(a)| |\psi\bigl(a'\bigr)| \,|{\rm d}\zeta(a)| = \int_0^1|\phi(a)| |\psi\bigl(a'\bigr)| |\zeta_{\alpha+\beta}| \,{\rm d}t,
\end{align*}
using the fact that ${\rm d}\zeta(a)={\rm d}\zeta_\alpha(a)={\rm d}(t\zeta_{\alpha+\beta})=\zeta_{\alpha+\beta} {\rm d}t$ in the second step. We can then bound the right-hand side in terms of the norms of $\phi$ and $\psi$
\begin{gather*}
|\zeta_{\alpha+\beta}|\int_0^1 |\phi(a)| |\psi\bigl(a'\bigr)| \,{\rm d}t
\\
\qquad\le \|\phi\|_{\sigma,\Lambda} \|\psi\|_{\sigma',\Lambda'} |\zeta_{\alpha+\beta}|\int_0^1 |\zeta_\alpha(a)|^{\sigma}{\rm e}^{\Lambda |\zeta_{\alpha}(a)|} |\zeta_\beta\bigl(a'\bigr)|^{\sigma'}{\rm e}^{\Lambda' |\zeta_{\beta}\bigl(a'\bigr)|} \,{\rm d}t\\
\qquad= \|\phi\|_{\sigma,\Lambda} \|\psi\|_{\sigma',\Lambda'} |\zeta_{\alpha+\beta}|^{1+\sigma+\sigma'}\int_0^1 |t|^{\sigma}{\rm e}^{\Lambda t |\zeta_{\alpha+\beta}|} |1-t|^{\sigma'}{\rm e}^{\Lambda' (1-t)|\zeta_{\alpha+\beta}|} \,{\rm d}t\\
\qquad\lesssim |\zeta_{\alpha+\beta}| ^{1+\sigma+\sigma'} {\rm e}^{(\Lambda+\Lambda')|\zeta_{\alpha+\beta}|}\int_0^1 t^{\sigma} (1-t)^{\sigma'} \,{\rm d}t\\
\qquad\lesssim |\zeta_{\alpha+\beta}| ^{1+\sigma+\sigma'} {\rm e}^{(\Lambda+\Lambda')|\zeta_{\alpha+\beta}|}.\tag*{\qed}
\end{gather*}\renewcommand{\qed}{}
\end{proof}

The product of two Laplace transforms based at different points is the Laplace transform of a~convolution.
\begin{Proposition}
For any $\phi\in\singexp{\sigma}{\Lambda}(\Omega_\alpha)$ and $\psi\in\singexp{\sigma'}{\Lambda'}(\Omega_\beta)$ with $\sigma,\sigma'\in (-1,\infty)$, we have
\begin{equation}\label{eqn:gen-laplace-convol}
\laplace^\theta_{\zeta,\alpha+\beta}[\phi \ast_{\zeta,\alpha,\beta} \psi] = \laplace^\theta_{\zeta,\alpha}[\phi] \laplace^\theta_{\zeta,\beta}[\psi].
\end{equation}
\end{Proposition}
More generally, this proposition holds for any holomorphic functions $\phi$ and $\psi$ for which the convolution and all of the Laplace transforms in equation~\eqref{eqn:gen-laplace-convol} are well defined.
\begin{proof}
For simplicity, we will write out the proof for $\theta = 0$. The proof for general $\theta$ is analogous.

The left-hand side of equation~\eqref{eqn:gen-laplace-convol} expands to
\begin{align*}
\laplace_{\zeta,\alpha,\beta}[\phi \ast_{\zeta,\alpha, \beta} \psi] &= \int_{\mathcal{J}_{\zeta, \alpha+\beta}} \,{\rm d}\zeta {\rm e}^{-z\zeta}\int_{\mathcal{S}_{\zeta,\alpha,\beta}} \phi(a) \psi\bigl(a'\bigr) \,{\rm d}\zeta(a),
\end{align*}
where $\mathcal{S}_{\zeta, \alpha, \beta}$ is the path described in Definition~\ref{def:convolution}. On this path, the equation
\begin{equation}\label{eq:slice-plane}
\zeta_\alpha(a) + \zeta_\beta\bigl(a'\bigr) = \zeta_{\alpha+\beta}
\end{equation}
holds, and the coordinates $\zeta_\alpha(a)$ and $\zeta_\beta\bigl(a'\bigr)$ have the same argument as $\zeta_{\alpha+\beta}$. In particular, over the ray $\mathcal{J}_{\zeta, \alpha+\beta}$, the coordinates $\zeta_\alpha(a)$, $\zeta_\beta\bigl(a'\bigr)$ and $\zeta_{\alpha+\beta}$ are all real and non-negative.\footnote{For general $\theta$, we have the condition that $\zeta_\alpha(a)$, $\zeta_\beta\bigl(a'\bigr)$ and $\zeta_{\alpha+\beta}$ all lie on the ray ${\rm e}^{{\rm i}\theta}[0, \infty)$.} The point where $\mathcal{S}_{\zeta, \alpha, \beta}$ starts has $\zeta_\alpha(a) = 0$, and the point where $\mathcal{S}_{\zeta, \alpha, \beta}$ ends has $\zeta_\beta\bigl(a'\bigr) = 0$.

The convolution and Laplace integrals can be combined into an integral over a real two-dimensional wedge $\mathcal{N}_{\zeta, \alpha, \beta}$ within the complex three-dimensional space $B^3$
\[ \laplace_{\zeta,\alpha+\beta}[{\phi}\ast_{\zeta, \alpha,\beta}{\psi}] = \int_{\mathcal{N}_{\zeta, \alpha, \beta}} {\rm e}^{-z\zeta} \phi(a) \psi\bigl(a'\bigr) \,{\rm d}\zeta(a) \wedge {\rm d}\zeta. \]
The coordinates $\zeta_\alpha(a)$, $\zeta_\beta\bigl(a'\bigr)$, $\zeta_{\alpha+\beta}$ parameterize $B^3$, and the real wedge $\mathcal{N}_{\zeta, \alpha, \beta}$ is cut out by equation~\eqref{eq:slice-plane} and the conditions
$
a \in \mathcal{J}_{\zeta, \alpha}$,
$a' \in \mathcal{J}_{\zeta, \beta}$,
$\zeta_{\alpha+\beta} \in [0, \infty)$.
Equation~\eqref{eq:slice-plane} can be rewritten as $\zeta(a) + \zeta\bigl(a'\bigr) = \zeta$, which implies the identities
\begin{gather*}
{\rm d}\zeta(a) + {\rm d}\zeta\bigl(a'\bigr) = {\rm d}\zeta ,\\
{\rm d}\zeta(a) \wedge \bigl[{\rm d}\zeta(a) + {\rm d}\zeta\bigl(a'\bigr)\bigr] = {\rm d}\zeta(a) \wedge {\rm d}\zeta, \\
{\rm d}\zeta(a) \wedge {\rm d}\zeta\bigl(a'\bigr) = {\rm d}\zeta(a) \wedge {\rm d}\zeta.
\end{gather*}
We can substitute the last identity into the integral over $\mathcal{N}_{\zeta, \alpha, \beta}$ to get
\begin{align*}
\laplace_{\zeta,\alpha+\beta}[{\phi}\ast_{\zeta, \alpha,\beta}{\psi}] & = \int_{\mathcal{N}_{\zeta, \alpha, \beta}} {\rm e}^{-z[\zeta(a) + \zeta(a')]} \phi(a) \psi\bigl(a'\bigr) \,{\rm d}\zeta(a) \wedge {\rm d}\zeta\bigl(a'\bigr) \\
& = \int_{\mathcal{N}_{\zeta, \alpha, \beta}} \bigl[ {\rm e}^{-z\zeta(a)} \phi(a)\, {\rm d}\zeta(a) \bigr] \wedge \bigl[ {\rm e}^{z\zeta(a')} \psi\bigl(a'\bigr) \,{\rm d}\zeta\bigl(a'\bigr) \bigr].
\end{align*}
Thinking of $\mathcal{N}_{\zeta, \alpha, \beta}$ as the graph of $\zeta_\alpha(a) + \zeta_\beta\bigl(a'\bigr)$ over the quadrant $\mathcal{J}_{\zeta, \alpha} \times \mathcal{J}_{\zeta, \beta}$, we can use Fubini's theorem to express the integral as the product
\[ \laplace_{\zeta,\alpha+\beta}[{\phi}\ast_{\zeta, \alpha,\beta}{\psi}] = \bigg[ \int_{\mathcal{J}_{\zeta, \alpha}} {\rm e}^{-z\zeta} \phi\, {\rm d}\zeta \bigg] \bigg[ \int_{\mathcal{J}_{\zeta,\beta}} {\rm e}^{-z\zeta} \psi\, {\rm d}\zeta \bigg]. \tag*{\qed}
\]\renewcommand{\qed}{}
\end{proof}

\begin{Corollary}\label{prop:convolution_iso_laplace}
For any two functions $\phi \in \singexpalg{\sigma}(\Omega)$ and $\psi \in \singexpalg{\sigma'}(\Omega)$ with $\sigma, \sigma' \in (-1, \infty)$, we have
$
\laplace_{\zeta,0}[\phi\ast_{\zeta}\psi]=\laplace_{\zeta,0}[\phi] \laplace_{\zeta,0}[\psi]$.
\end{Corollary}

On a domain $\Omega$ which is star-shaped around $\zeta = 0$, the convolution product $\ast_{\zeta}$ induces a~$\C$-algebra structure on $\singexp{\sigma}{\Lambda}(\Omega)$. The Laplace transform $\laplace_{\zeta, 0}$ is special, because it is an algebra isomorphism with respect to $\ast_{\zeta}$.
\begin{Remark}
When the Taylor expansions of $\phi$ and $\psi$ around $\zeta=\alpha$ converge in a disk of radius $R$, the Taylor expansion of the convolution $\phi\ast_{\zeta,\alpha}\psi$ around $\zeta=\alpha$ also converges in a disk of radius $R$. This is proven for $\alpha = 0$ in \cite[Lemma~5.14]{diverg-resurg-i}.
\end{Remark}

{\bf The Borel transform as algebra isomorphism.}
The Borel transform $\borel_\zeta$ is defined on the space of formal power series $z^{-1}\C\big\llbracket z^{-1}\big\rrbracket$, which is a $\C$-algebra, with product given by the Cauchy product of formal power series. In addition, introducing the formal unit $\delta$, the Borel transform can be extended to $\C\big\llbracket z^{-1}\big\rrbracket$ (including constants), as $\borel_\zeta(1)=:\delta$. Therefore, the Borel transform is a $\C$-linear map
$\borel_\zeta\maps\C\big\llbracket z^{-1}\big\rrbracket\to\C\delta + \C\llbracket\zeta\rrbracket$.
In fact, since $\borel_\zeta$ is defined as the formal inverse of $\laplace_{\zeta,0}$ we can introduce the \textit{formal Laplace transform} as the map from $\C\llbracket \zeta\rrbracket$ to $z^{-1}\C\big\llbracket z^{-1}\big\rrbracket$ which acts as $\laplace_{\zeta,0}$ on monomials $\zeta^k$ and then it extends by countable linearity. As a result, we find that the Borel transform is a $\C$-linear invertible map whose inverse is the formal Laplace transform.

In addition, on both $\C\llbracket\zeta\rrbracket$ and $\C\big\llbracket z^{-1}\big\rrbracket$ there is the natural product given by the Cauchy product of formal power series. However, the Borel transform does not extend to a $\C$-algebra isomorphism, because $\borel_\zeta\bigl(\series{\Phi}\series{\Psi}\bigr)$ is different from the Cauchy product $\bigl(\borel_\zeta \series{\Phi}\bigr)\bigl(\borel_\zeta \series{\Psi}\bigr)$.

To see the Borel transform as an algebra isomorphism, we must use a different product on~$\C\llbracket\zeta\rrbracket$: the convolution product for formal power series. Let $\zeta$ be a coordinate on $B$. For any integers $m$, $n$,
\begin{equation}\label{eq: convolution_formal}
\frac{\zeta^m}{m!} \ast_\zeta \frac{\zeta^n}{n!} = \frac{\zeta^{m+n+1}}{(m+n+1)!} .
\end{equation}
\begin{Definition}\label{def:convolution_formal}
We extend the convolution product $\ast_\zeta$ by countable linearity to the whole ring of formal power series $\C\llbracket\zeta\rrbracket$
\begin{equation*}
 \left( \sum_{n=0}^\infty a_n \frac{\zeta^n}{n!}\right)\ast_\zeta\left(\sum_{n=0}^\infty b_n\frac{\zeta^n}{n!}\right)=\sum_{n=0}^\infty\sum_{j=0}^na_jb_{n-j}\left[\frac{\zeta^j}{j!}\ast_\zeta\frac{\zeta^{n-j}}{(n-j)!}\right]=\sum_{n=0}^\infty\sum_{j=0}^na_jb_{n-j}\frac{\zeta^{n+1}}{(n+1)!}.
\end{equation*}
Furthermore, we define the convolution product of $\delta$ with monomials
$
\delta\ast_\zeta\frac{\zeta^n}{n!}:=\frac{\zeta^n}{n!}$, $
\frac{\zeta^n}{n!}\ast_\zeta\delta:=\frac{\zeta^n}{n!}$,
in other words, $\delta$ is the unit for the convolution.
\end{Definition}
Notice that from equation~\eqref{eq: convolution_formal} and the definition of the convolution with $\delta$, we deduce the convolution product is commutative with unit $\delta$.

We can show that $\borel_\zeta$ is an algebra homomorphism from $\C\big\llbracket z^{-1}\big\rrbracket$ endowed with the Cauchy product to $\C\delta+\C\llbracket\zeta\rrbracket$ endowed with the convolution product.
\begin{Proposition}
Let $z$ be a coordinate on the cotangent fiber over $\zeta=0$, and $\series{\Phi},\series{\Psi}\in\C\big\llbracket z^{-1}\big\rrbracket$. The Borel transform $\borel_\zeta$ is an algebra homomorphism, namely
\smash{$\borel_\zeta\bigl(\series{\Phi}\series{\Psi}\bigr)=\borel_\zeta\series{\Phi} \ast_\zeta \borel_\zeta\series{\Psi}$}.
\end{Proposition}
\begin{proof}
 We assume $\series{\Phi}=\sum_{n=0}^\infty a_nz^{-n}$ and $\series{\Psi}=\sum_{n=0}^\infty b_nz^{-n}$. Then,
 \begin{align*}
\borel_\zeta\bigl(\series{\Phi}\series{\Psi}\bigr)={}&\borel_\zeta\Bigg(\sum_{n=0}^\infty\sum_{j=0}^na_jb_{n-j}z^{-n}\Bigg)=a_0b_0\delta +\sum_{n=1}^\infty\sum_{j=0}^na_jb_{n-j}\frac{\zeta^{n-1}}{(n-1)!}\\
={}&a_0b_0\delta+a_0\sum_{n=1}^\infty b_{n}\frac{\zeta^{n-1}}{(n-1)!}+b_0\sum_{n=1}^\infty a_n\frac{\zeta^{n-1}}{(n-1)!} \\
& + \sum_{n=1}^\infty\sum_{j=1}^{n-1}a_jb_{n-j}\left[\frac{\zeta^{j-1}}{(j-1)!} \ast_\zeta \frac{\zeta^{n-j-1}}{(n-j-1)!}\right]\\
={}&\left(a_0\delta+\sum_{n=1}^\infty a_n\frac{\zeta^{n-1}}{(n-1)!}\right)\ast_\zeta\left(b_0\delta+\sum_{n=1}^\infty b_n\frac{\zeta^{n-1}}{(n-1)!}\right)\\
={}&\borel_\zeta\series{\Phi} \ast_\zeta \borel_\zeta\series{\Psi}.\tag*{\qed}
\end{align*}\renewcommand{\qed}{}
\end{proof}

Similarly, $\borel_{\zeta_\alpha}$ is an algebra homomorphism from $\C\big\llbracket z^{-1}\big\rrbracket$ endowed with the Cauchy product to~${\C\delta+\C\llbracket\zeta_\alpha\rrbracket}$ endowed with the convolution product $\ast_{\zeta_\alpha}$, where $z$ is a coordinate on $T^*_{\zeta=\alpha}B$.

\begin{Remark}
If $\series{\phi}$, $\series{\psi}$ are convergent formal power series in $\C\{\zeta\}$, their convolution product is convergent too, and its sum is the convolution of the sums of $\series{\phi}$, $\series{\psi}$, respectively, as defined in Definition~\ref{def:convolution}. We will discuss it in Section~\ref{sec:Borel-gevrey}.
\end{Remark}

\subsection{Bridging the gap between formal and analytic objects}\label{sec:bridging}

\subsubsection{Formal to analytic}\label{sec:Borel-gevrey}
If we want to get analytic information from the Borel transform---for example, by doing Borel summation---we should focus on the series whose Borel transforms converge. These series, called {\em $1$-Gevrey} series, form a subspace $\C \big\llbracket z^{-1} \big\rrbracket_1$ within $\C \big\llbracket z^{-1} \big\rrbracket$. They are characterized by the factorial growth of their coefficients.
\begin{Definition}
A series
\smash{$\sum_{n \ge -1} a_n z^{-n-1}$}
is {\em $1$-Gevrey} if there is some $A > 0$ with $|a_n| \lesssim A^n n!$ over all $n \ge 0$.
\end{Definition}
\begin{Proposition}\label{prop:gevrey_to_convergent}
A series
\smash{$ \series{\Phi} = \sum_{n \ge -1} a_n z^{-n-1}$}
is $1$-Gevrey if and only if its Borel transform~${\series{\phi} = \borel_\zeta \series{\Phi}}$ is convergent.
\end{Proposition}
\begin{proof}
From the expression
\[ \series{\phi} = a_{-1} \delta + \sum_{n \ge 0} a_n \frac{\zeta^n}{n!}, \]
we see that $\series{\phi}$ converges if and only if there is some $A > 0$ with $|a_n / n!| \lesssim A^n$ over all $n \ge 0$. The constant $A$ is the radius of convergence of $\series{\phi}$.
\end{proof}

The space of $1$-Gevrey series is closed under the Cauchy product. In addition,
\[\borel_\zeta\maps\ \C\big\llbracket z^{-1}\big\rrbracket_1\to\C\delta + \C\{\zeta\}\]
extends to an algebra homomorphism, where $\C\{\zeta\}$ is endowed with the convolution product defined in Definition~\ref{def:convolution_formal}. In fact, we deduce that the Borel transform induces an algebra homomorphism on $1$-Gevrey series. Similarly, if $z$ is coordinate on the fiber over $\zeta=\alpha$, the Borel transform $\borel_{\zeta_\alpha}$ induces an isomorphism between $\C\big\llbracket z^{-1}\big\rrbracket_1$ endowed with the Cauchy product, and~${\C\delta+\C\{\zeta_\alpha\}}$ endowed with the convolution $\ast_{\zeta_\alpha}$.
\subsubsection{Analytic to formal}
Consider an open set $\Omega \subset \C$ that touches but does not contain $\zeta = 0$.\footnote{We say that a subset of a topological space {\em touches} the points in its closure. This is a way of expressing the nearness relation associated with the topology~\cite[Chapter~5, Definition~2.11]{joshi1983gen-top}.} The spaces $\singexp{\sigma}{\Lambda}(\Omega)$ from~\cite{reg-sing-volterra} are convenient for analytic discussions of the Laplace transform, as we saw in Section~\ref{sec:reg-decay}. To study Borel summation, however, we also need to take shifted Taylor expansions around~${\zeta = 0}$, turning holomorphic functions $\phi \in \singexp{\sigma}{\Lambda}(\Omega)$ into convergent power series~$\smash{\series{\phi} \in \zeta^\rho \C\{\zeta\}}$ with $\rho \ge \sigma$. The following discussion, leading up to Lemma~\ref{lem:shifted_holo_closed}, will help us stay within a~subspace of $\singexp{\sigma}{\Lambda}(\Omega)$ where the shifted Taylor series exists.

If a function on $\Omega$ has an ordinary Taylor expansion in $\C\{\zeta\}$, summing the series extends the function to a neighborhood of $\zeta = 0$. The subspace of $\singexp{\sigma}{\Lambda}(\Omega)$ where the ordinary Taylor series exists can thus be seen as the intersection
$\mathcal{O}_{\zeta = \alpha} \cap \singexp{\sigma}{\Lambda}(\Omega)$,
where $\mathcal{O}_{\zeta = \alpha}$ is the algebra of germs of holomorphic functions at $\zeta$. More generally, for any $\rho \ge \sigma$, the subspace of $\singexp{\sigma}{\Lambda}(\Omega)$ where the $\rho$-shifted Taylor series exists can be seen as the intersection~${ \zeta^\rho \mathcal{O}_{\zeta = \alpha} \cap \singexp{\sigma}{\Lambda}(\Omega)}$.

To study this subspace, it will be helpful to work on the universal covering \smash{$\pi \maps \widetilde{\C^\times} \to \C^\times$}. On a simply connected open set \smash{$\Upsilon \subset \widetilde{\C^\times}$} that touches $\zeta = 0$, a function has an ordinary Taylor expansion around $\zeta = 0$ if and only if it descends and extends to a holomorphic function on a~neighborhood of $\zeta = 0$ in $\C$. Thus, once again, we can see the subspace of $\singexp{\sigma}{\Lambda}(\Upsilon)$ where the $\rho$-shifted Taylor series exists as the intersection
$ \zeta^\rho \mathcal{O}_{\zeta = 0} \cap \singexp{\sigma}{\Lambda}(\Upsilon)$.
\begin{Lemma}\label{lem:shifted_holo_closed}
Consider a simply connected open set \smash{$\Upsilon \subset \widetilde{\C^\times}$} which touches $\zeta = 0$ and is star-shaped around $\zeta = 0$. Suppose that $\pi \maps \Upsilon \to \C^\times$ overlaps itself at the edges: it is two-to-one over a simply connected open set $\Upsilon_\Delta \subset \C^\times$ which is star-shaped around $\zeta = 0$, and it is one-to-one everywhere else.

\begin{figure}[!ht]\centering\vspace{-15mm}

\hspace*{-19mm}\begin{tikzpicture}
 \newcommand{\spill}{4}
 \fill[ietcoast!20, bezier bounding box, path fading=radial edge]
 (-\spill, -\spill) (\spill, \spill)
 (190:1.5) arc (190:41:1.5)--(41:\spill) arc (41:-48:\spill)--(-48:1.3) arc (-48:-170:1.3)--cycle;
 \fill[ietcoast!33] (0, 0)--(-170:1.3) arc (-170:-210:1.3)--cycle;
 \fill[ietcoast!33!black] circle (0.7mm) node[anchor=west, black, outer sep=1mm] {$\zeta = \alpha$};
 \node[ietcoast!33!black] at (-190:0.85) {$\Upsilon_\Delta$};
 \foreach \p in {(107.5:1.95), (46.5:2.46), (-56:2.24)} {
 \fill[ietcoast!33!black] \p circle (0.7mm);
 }
\end{tikzpicture}\vspace{-14mm}

\caption{The domain $\Upsilon$ at $\zeta=\alpha$.}\label{fig:domain_Upsilon}
\end{figure}
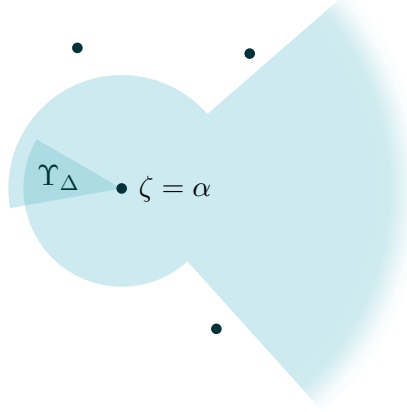

\noindent
For any $\sigma \in \R$, $\Lambda \in \R$, and non-integer $\rho \ge \sigma$, the intersection~${\zeta^\rho \mathcal{O}_{\zeta = 0} \cap \singexp{\sigma}{\Lambda}(\Upsilon)}$
is a~closed subspace of $\singexp{\sigma}{\Lambda}(\Upsilon)$.
\end{Lemma}

\begin{proof}
Each point $b \in \Upsilon_\Delta$ has two preimages $b_+, b_- \in \Upsilon$, where $b_+$ is $2\pi$ counter-clockwise of~$b_-$ around $\zeta = \alpha$. The {\em variation} operator $\var \maps \singexp{\sigma}{\Lambda}(\Upsilon) \to \singexp{\sigma}{\Lambda}(\Upsilon_\Delta)$, defined by\footnote{The target of the variation operator is a space of functions on $\Upsilon_\Delta$ because the variation is defined only on the overlap. Our notation follows the convention of \cite{diverg-resurg-i}, described before Definition~6.46.}
$\var \varphi \big|_b = \varphi(b_+) - \varphi(b_-) $
is bounded. A function in \smash{$\singexp{\sigma}{\Lambda}(\Upsilon)$} descends and extends to a function on a~neighborhood of~${\zeta = 0}$ in $\C$ if and only if it lies in the kernel of $\var - \bigl({\rm e}^{2 \pi {\rm i} \rho} - 1\bigr)$. This operator is bounded, so its kernel is closed.
\end{proof}

Taylor expansion and summation form a two-way link between holomorphic functions and formal power series. For regular enough functions, the formal Laplace transform preserves a~one-way vestige of this link: it turns Taylor expansions on the position domain into asymptotic expansions on the frequency domain.
\begin{Lemma}\label{lem:laplace-bridge}
Let $\Omega$ be an open sector with $\zeta = 0$ at its tip, and an opening angle of $\pi$ or less. Take some $\rho \in (-1,\infty)$. If \smash{$\varphi \in \singexp{\rho}{\Lambda}(\Omega)$} is the sum of a power series $\series{\varphi} \in \zeta^\rho \C\{\zeta\}$, then its Laplace transform \smash{$\Phi = \laplace_{\zeta, 0}^\theta \varphi$} is asymptotic to the formal Laplace transform \smash{$\series{\Phi} = \laplace_{\zeta, 0}^\theta \series{\varphi}$} on the domain \smash{$\widehat{\Omega}_0^\Lambda$} defined in Section~{\rm\ref{sec:reg-decay}}.
\end{Lemma}
\begin{Remark}
This result is similar to \cite[Theorem~5.20]{diverg-resurg-i}. It is weaker in that it does not establish Gevrey asymptoticity, but stronger in that it can be applied to functions with fractionally shifted Taylor expansions.
\end{Remark}
\begin{proof}[Proof of Lemma~\ref{lem:laplace-bridge}]
By adding points within the disk around $\zeta = 0$ where $\series{\varphi}$ converges, expand $\Omega$ to an open set $\Upsilon \subset \widetilde{\C^\times}$ of the kind used in Lemma~\ref{lem:shifted_holo_closed}, Figure~\ref{fig:domain_Upsilon}. Write $\series{\varphi}$ out in the form
\[ \series{\varphi} = a_0 \frac{\zeta^\rho}{\Gamma(\rho+1)} + a_1 \frac{\zeta^{\rho+1}}{\Gamma(\rho+2)} + a_2 \frac{\zeta^{\rho+2}}{\Gamma(\rho+3)} + a_3 \frac{\zeta^{\rho+3}}{\Gamma(\rho+4)} + \cdots, \]
observing that
\[ \series{\Phi} = \frac{a_0}{z^{\rho+1}} + \frac{a_1}{z^{\rho+2}} + \frac{a_2}{z^{\rho+3}} + \frac{a_3}{z^{\rho+4}} + \cdots. \]
For each $n \in \Z_{\ge 0}$, define the tail sum $\varepsilon_n \in \singexp{\rho+n}{\Lambda}(\Upsilon)$ by the equation
\[ \varphi = a_0 \frac{\zeta^\rho}{\Gamma(\rho+1)} + a_1 \frac{\zeta^{\rho+1}}{\Gamma(\rho+2)} + \cdots + a_{n-1} \frac{\zeta^{\rho+n-1}}{\Gamma(\rho+n)} + \varepsilon_n. \]
Taking the Laplace transform of both sides, we learn that
\[ \Phi = \frac{a_0}{z^{\rho+1}} + \frac{a_1}{z^{\rho+2}} + \cdots + \frac{a_{n-1}}{z^{\rho+n}} + E_n, \]
where $E_n = \laplace_{\zeta, 0}^\theta \varepsilon_n$. To show that $\Phi$ is asymptotic to $\series{\Phi}$, we just need to prove that $|z^{\rho+n} E_n| \to 0$ as $|z| \to \infty$ along every ray.

Because of the conditions we put on $\Omega$, we can apply Proposition~\ref{prop:laplace-cont}, which tells us that~$E_n$ lies in \smash{$\dualsingexp{-\rho-n-1}\bigl(\widehat{\Omega}_0^\Lambda\bigr)$}. This means that $|E_n| \lesssim \Delta^{-\rho-n-1}$, where $\Delta$ is the distance to the boundary of \smash{$\widehat{\Omega}_0^\Lambda$}. Observe that \smash{$\widehat{\Omega}_0^\Lambda$} is a union of half-planes, the distance to the boundary of each half-plane increases linearly along any ray that intersects the half-plane, and $\Delta$ is the infimum of these distances. From this, we can deduce that $|z|/\Delta$ goes to a constant limit as~${z \to \infty}$ along any ray that intersects $\widehat{\Omega}_0^\Lambda$.\footnote{The constant limit is not uniform over all rays, so we are showing that $\Phi$ is asymptotic, but not necessarily uniformly asymptotic, to $\series{\Phi}$.} We then reason that
\begin{align*}
|z^{\rho+n} E_n| & \lesssim |z|^{\rho+n} \Delta^{-\rho-n-1} = (|z|/\Delta)^{\rho+n} \Delta^{-1} \lesssim \Delta^{-1}
\end{align*}
along any ray that intersects $\widehat{\Omega}_0^\Lambda$. This, as discussed above, proves the desired result.
\end{proof}

\subsection{Action on fractional derivatives and fractional integrals}\label{sec:frac-diff-laplace}
In this subsection, we will consider the action of $\borel_\zeta$ and $\laplace_{\zeta,0}$ with respect to the action of derivatives in both the frequency and position domains, as well as fractional integrals and derivatives~$\fracderiv{\lambda}{\zeta}{0}$ in the position domain. This will be useful in Section~\ref{sec:examples}.

In our treatment, we will restrict the Borel transform to the space of $1$-Gevrey series and the Laplace transform to $\singexp{\sigma}{\Lambda}(\Omega_\alpha)$ with $\sigma>-1$ and $\Lambda>0$, for some domain $\Omega_\alpha$ as in Figure~\ref{fig:sectorial_domain-pos-fre}.
\begin{Definition}
For $\nu \in (-\infty, 1)$, the \textit{fractional integral} $\partial^{\nu-1}_{\zeta, \alpha}$ is defined by
\[ \partial^{\nu-1}_{\zeta, \alpha} \phi := \frac{1}{\Gamma(1-\nu)} \int_{\mathcal{S}_{\zeta,\alpha}} (\zeta-\zeta(a))^{-\nu} \phi(a) \,{\rm d}\zeta(a), \]
where $\mathcal{S}_{\zeta,\alpha}$ is the line segment from $\zeta(a)=\alpha$ to $\zeta(a)=\zeta$.
\end{Definition}
The fractional integral obeys the semigroup law \cite[Section~1.3]{mladenov2014advanced}
\smash{$
\fracderiv{\lambda}{\zeta}{\alpha} \fracderiv{\mu}{\zeta}{\alpha} = \fracderiv{\lambda+\mu}{\zeta}{\alpha} $}, $ \lambda, \mu \in (-\infty, 0)$,
and agrees with ordinary repeated integration when $\nu$ is an integer~\cite[equation~(35)]{mladenov2014advanced}.

For $\mu \in (0, 1)$ and integers $n \ge 0$, fractional derivatives \smash{$\fracderiv{n+\mu}{\zeta}{0}$} are defined by composing \smash{$\fracderiv{\mu-1}{\zeta}{0}$} with powers of \smash{$\frac{\partial}{\partial \zeta}$}. However, \smash{$\fracderiv{\mu-1}{\zeta}{0}$} and \smash{$\frac{\partial}{\partial \zeta}$} do not commute~\cite[equation~(54)]{mladenov2014advanced}. Various ordering conventions give various definitions of \smash{$\fracderiv{n+\mu}{\zeta}{0} \phi$}, which differ by operators that act on the germ of~$\phi$ at zero (see~\cite[Section 1.3]{mladenov2014advanced} and \cite{podlubny}). We will use the {\em Riemann--Liouville} convention.

\begin{Definition}\label{definition:frac_driv}
For $\mu \in (0, 1)$ and integers $n \ge 0$, the {\em Riemann--Liouville fractional derivative}~$\fracderiv{n+\mu}{\zeta}{0}$ is defined by
\[ \fracderiv{n+\mu}{\zeta}{0} := \left(\frac{\partial}{\partial \zeta}\right)^{n+1} \fracderiv{\mu-1}{\zeta}{0}. \]
\end{Definition}
The symbol $\fracderiv{\mu}{\zeta}{0}$ is now defined for any $\mu \in \R \smallsetminus \{0, 1, 2, 3, \dots\}$. It denotes a fractional integral when $\mu$ is negative, and a fractional derivative when $\mu$ is a positive non-integer.

The Riemann--Liouville fractional derivative is a left inverse of the fractional integral, in the sense that $\fracderiv{\lambda}{\zeta}{ 0} \fracderiv{-\lambda}{\zeta}{0}$ is the identity for all $\lambda \in (0, \infty)$. This extends the semigroup law
\begin{align*}
\fracderiv{\lambda}{\zeta}{0} \fracderiv{\mu}{\zeta}{0} = \fracderiv{\lambda+\mu}{\zeta}{0} ,\qquad \lambda \in \R \smallsetminus \{0, 1, 2, 3, \dots\},\quad\mu \in (-\infty, 0).
\end{align*}
\begin{Proposition}\label{prop:L-int-op}
For any $\phi\in\singexp{\sigma}{\Lambda}(\Omega)$ with $\sigma>-1$ and $\Lambda>0$, we have
\begin{enumerate}[label=(\roman*)]\itemsep=0pt
\item[$(1)$] 
for fractional derivatives:
$\laplace_{\zeta,0} \bigl[\fracderiv{\mu}{\zeta}{0} \phi\bigr] = z^{\mu} \laplace_{\zeta,0} \phi$ for every $\mu\in(0,1)$;
\item[$(2)$] 
for fractional integrals:
$\laplace_{\zeta,0} \bigl[\fracderiv{\lambda}{\zeta}{0} \phi\bigr] = z^{\lambda} \laplace_{\zeta,0} \phi$ for every $\lambda\in(-\infty,0)$;
\item[$(3)$] 
for whole derivatives:
$\laplace_{\zeta,0} \bigl[\bigl(\frac{\partial}{\partial\zeta}\bigr)^n \phi\bigr] = z^{n} \laplace_{\zeta,0} \phi$ for all $n\in\{1, 2, 3, \dots\}$, under the condition that \smash{$\phi^{(0)}, \phi^{(1)}, \dots, \phi^{(n)}$} all vanish at $\zeta = 0$;
\item[$(4)$] 
for monomial multiplication:
$\laplace_{\zeta,0} \bigl[\zeta^n \phi\bigr] = \bigl(-\frac{\partial}{\partial z}\bigr)^n \laplace_{\zeta,0} \phi$ for every $n\in\{0,1,2,3,\dots\}$.
\end{enumerate}
\end{Proposition}
\begin{proof}
Results (1) and (2) follow from the properties of the Laplace transform with respect to the convolution product. For every $\tau \in (-\infty, 0) \cup (0,1)$,
\begin{align*}
z^\tau\laplace_{\zeta,0}\phi&=\laplace_{\zeta,0}\left[\frac{\zeta^{-\tau-1}}{\Gamma(-\tau)}\right] \laplace_{\zeta,0}\phi=\laplace_{\zeta,0}\left[\frac{\zeta^{-\tau-1}}{\Gamma(-\tau)} \ast_\zeta \phi\right] =\laplace_{\zeta,0} \int_{\mathcal{S}_{\zeta,0}}\frac{\bigl(\zeta-\zeta(a)\bigr)^{-\tau-1}}{\Gamma(-\tau)} \phi(a)\, {\rm d}\zeta(a)\\
&=\laplace_{\zeta,0} \int_{\mathcal{S}_{\zeta,0} }\frac{\bigl(\zeta-\zeta(a)\bigr)^{-\tau-1}}{\Gamma(-\tau)} \phi(a) \, {\rm d}\zeta(a)=\laplace_{\zeta,0}\bigl[\fracderiv{\tau}{\zeta}{0}\phi\bigr],
\end{align*}
where $\mathcal{S}_{\zeta,0}$ is the line segment from $\zeta(a)=0$ to $\zeta(a) = \zeta$. Result~(3) is proven by repeated integration by parts, with the condition on $\phi^{(0)}, \phi^{(1)}, \dots \phi^{(n)}$ ensuring that the boundary terms vanish
\begin{align*}
 z^n\laplace_{\zeta,0}\phi&=\int_0^\infty {\rm e}^{-z\zeta} z^n \phi {\rm d}\zeta=(-1)^n\int_0^\infty \partial_\zeta^n \bigl[{\rm e}^{-z\zeta}\bigr] \phi \, {\rm d}\zeta=\laplace_{\zeta,0}\bigl[\partial_\zeta^n \phi\bigr].
\end{align*}
Result~(4) is also proven by repeated integration by parts, with the assumption that $\zeta^n\phi\in\singexp{\sigma}{\Lambda}(\Omega)$ ensuring that the integral converges.
\end{proof}

\begin{Lemma}\label{lem:frac-deriv-Borel}
For any non-integer $\mu \in (0, \infty)$ and any integer $k \ge 0$,
\[
\partial^\mu_{\zeta,0} \borel_\zeta \bigl[z^{-(k+1)}\bigr] = \borel_\zeta \bigl[z^\mu z^{-(k+1)}\bigr].\]
\end{Lemma}
\begin{proof}
Because every operation in the statement maps monomials to monomials, we do not need to distinguish carefully between formal and analytic arguments. Using a mostly analytic argument, we will show that for any $\mu \in (0, 1)$ and any integer $n \ge 0$, the claim holds with~${\mu = n + \alpha}$. First, evaluate
\begin{align*}
\partial^{\mu-1}_{\zeta,0} \borel_\zeta \bigl[z^{-(k+1)}\bigr] & = \frac{1}{\Gamma(1-\mu)} \int_0^\zeta (\zeta-\zeta(a))^{-\mu} \frac{\zeta(a)^k}{\Gamma(k+1)} \,{\rm d}\zeta(a) \\
& = \frac{1}{\Gamma(1-\mu) \Gamma(k+1)} \int_0^1 (\zeta-\zeta t)^{-\mu} (\zeta t)^k \zeta \,{\rm d}t \\
& = \frac{\zeta^{k-(\mu-1)}}{\Gamma(1-\mu) \Gamma(k+1)} \int_0^1 (1-t)^{-\mu} t^k \,{\rm d}t = \frac{\zeta^{k-(\mu-1)}}{\Gamma\bigl(k-(\mu-1)+1\bigr)}
\end{align*}
by reducing the integral to Euler's beta function~\cite[identity (5.12.1)]{dlmf}. This establishes that
\begin{equation}\label{frac-diff-laplace}
\left(\frac{\partial}{\partial \zeta}\right)^{n+1} \partial^{\mu-1}_{\zeta,0} \borel_\zeta \bigl[z^{-(k+1)}\bigr] = \frac{\zeta^{k-(n+\mu)}}{\Gamma(k-(n+\mu)+1)}
\end{equation}
for $n = -1$. If \eqref{frac-diff-laplace} holds for $n = m$, it also holds for $n = m+1$, because
\begin{align*}
\frac{\partial}{\partial \zeta} \left(\frac{\partial}{\partial \zeta}\right)^{m+1} \partial^{\mu-1}_{\zeta,0} \borel_\zeta \bigl[z^{-(k+1)}\bigr] & = \frac{\partial}{\partial \zeta} \left( \frac{\zeta^{k-(m+\mu)}}{\bigl(k-(m+\mu)\bigr) \Gamma\bigl(k-(m+\mu)\bigr)} \right)\\
& = \frac{\zeta^{k-(m+1+\mu)}}{\Gamma(k-(m+\mu))}.
\end{align*}
Hence, \eqref{frac-diff-laplace} holds for all $n \ge -1$, and the desired result quickly follows. The condition $\mu \in (0, 1)$ saves us from the trouble we would run into if $k-(m+\mu)$ were in $\Z_{\le 0}$. This is how we avoid the initial value corrections that appear in ordinary derivatives of Borel transforms.
\end{proof}

We can now prove the properties of the Borel transform analogous to the one of the Laplace transform from Proposition~\ref{prop:L-int-op}.
\begin{Proposition}\label{prop:frac-der-int-borel}
For any $1$-Gevrey series $\series{\Phi}\in\C\big\llbracket z^{-1}\big\rrbracket_1$, we have
\begin{enumerate}[label=(\roman*)]\itemsep=0pt
\item[$(1)$] 
for fractional derivatives:
\smash{$\borel_\zeta\bigl[ z^\mu \series{\Phi}\bigr]=\fracderiv{\mu}{\zeta}{0} \borel_\zeta\series{\Phi}$} for every $\mu\in(0, \infty)$;
\item[$(2)$] 
for fractional integrals:
\smash{$\borel_\zeta\bigl[z^\lambda \series{\Phi}\bigr]=\fracderiv{\lambda}{\zeta}{0} \borel_\zeta\series{\Phi}$} for every $\lambda\in(-\infty, 0)$;
\item[$(3)$] 
for whole derivatives:
\smash{$\borel_\zeta\bigl[\partial_z^{k} \series{\Phi}\bigr]=(-\zeta)^k\borel_\zeta\series{\Phi}$} for all $k\in\{0, 1, 2, 3, \dots\}$, under the condition that \smash{$\series{\Phi}\in z^{-k-1}\big\llbracket z^{-1}\big\rrbracket_1$}.\footnote{This property holds for all series in $\C\big\llbracket z^{-1}\big\rrbracket$, even ones that are not $1$-Gevrey.}
\end{enumerate}
\end{Proposition}
\begin{proof}
Result~(1) follows from Lemma~\ref{lem:frac-deriv-Borel}.
To prove result~(2), first notice that for $\lambda \in \R_{< 0}$,
\[\fracderiv{\lambda}{\zeta}{0}\zeta^{k}=\zeta^{k-\lambda}\frac{k!}{\Gamma(k-\lambda+1)}.\]
Then, given a series $\series{\Phi} = \sum_{n\geq 0} a_n z^{-n}$, we compute
\begin{align*}
\borel_\zeta \bigl[z^\lambda \series{\Phi}\bigr] & = \sum_{n\geq 0}a_n\frac{\zeta^{n-\lambda}}{\Gamma(n-\lambda+1)} = \sum_{n\geq 0}a_n\frac{\zeta^{n-\lambda}}{n!}\frac{n!}{\Gamma(n-\lambda+1)} = \sum_{n\geq 0} \frac{a_n}{n!} \fracderiv{\lambda}{\zeta}{0}\zeta^n = \fracderiv{\lambda}{\zeta}{0}\series{\Phi}
\end{align*}
using the fact that $\borel_\zeta\series{\Phi}$ is convergent in the last step.

Finally, result~(3) follows from a simple computation
\begin{align*}
\borel_\zeta\bigl[\partial_z^k\series{\Phi}\bigr]&=\borel_\zeta\left[\sum_{n=k+1}^\infty a_n \frac{(n+k-1)!}{(n-1)!} (-1)^{k} z^{-n-k}\right]\\
&=(-1)^k \sum_{n=k+1}^\infty a_n \frac{(n+k-1)!}{(n-1)!}\frac{\zeta^{n+k-1}}{(n+k-1)!}=(-\zeta)^k\sum_{n=k+1}^\infty a_n \frac{\zeta^{n-1}}{(n-1)!}\\
&=(-\zeta)^k\borel_\zeta \series{\Phi}.\tag*{\qed}
\end{align*} \renewcommand{\qed}{}
\end{proof}

\section{Proof of main results}\label{sec:proof_main_results}
\subsection{Borel regularity for ODEs}\label{borel_reg-ODE}
In Section~\ref{sec:regularity-results}, we show that linear level~$1$ ODEs of a certain form always have Borel regular solutions, indexed by the characteristic roots of their constant-coefficient parts. These solutions turn out to be asymptotic to the formal series solutions found by Poincar\'e.

In Section~\ref{sec:new-summability-proof}, we will look at the same class of ODEs from a more formal perspective, starting from the assumption that Poincar\'e's solutions exist and their Borel transforms converge. Using the existence and uniqueness results that characterize the Borel regular solutions, we give a new proof that Poincar\'e's solutions are Borel summable.
\subsubsection{Regularity results}\label{sec:regularity-results}
Let us recall the setting from Section~\ref{borel-reg:explanatory-power}. Let $\mathcal{P}$ be a linear differential operator of the form
\begin{equation}\label{eqn:operator-P}
\mathcal{P} = P\left(\frac{\partial}{\partial z}\right) + \frac{1}{z} Q\left(\frac{\partial}{\partial z}\right) + \frac{1}{z^2} R\bigl(z^{-1}\bigr),
\end{equation}
where
\begin{enumerate}\itemsep=0pt
\item[(1)] $P$ is a monic degree-$d$ polynomial whose roots are all simple;
\item[(2)] $Q$ is a degree-$(d-1)$ polynomial that is non-zero at every root of $P$;
\item[(3)] $R\bigl(z^{-1}\bigr)$ is holomorphic in some disk $|z| > A$ around $z = \infty$. In particular, the power series~\smash{$R\bigl(z^{-1}\bigr) = \sum_{j=0}^\infty R_j z^{-j}$}
converges in the region $|z| > A$.
\end{enumerate}
For each root $-\alpha$ of $P$, let $\hat{\mathcal{P}}_{\alpha}$ be the Volterra operator
\begin{equation}\label{eq:hat_P}
\hat{\mathcal{P}}_\alpha:=P(-\zeta)+\partial_{\zeta,\alpha}^{-1}\circ Q(-\zeta)+\partial_{\zeta,\alpha}^{-2}\circ R\bigl(\partial_{\zeta,\alpha}^{-1}\bigr),
\end{equation}
and define the constant \smash{$\tau_\alpha := \frac{Q(-\alpha)}{P'(-\alpha)}$}. From \cite[Lemma~4.3]{reg-sing-volterra}, we know that if $\psi_\alpha$ satisfies the equation \smash{$\hat{\mathcal{P}}_\alpha\psi_\alpha=0$}, then its Laplace transform $\Psi_\alpha:=\laplace_{\zeta,\alpha}^{\theta}\psi_\alpha$ satisfies the equation $\mathcal{P}\Psi_\alpha=0$, as long as the Laplace transform is well defined.

We can now pick out some special solutions of the equation $\mathcal{P}\Psi = 0$, which are associated with roots of $P$. We will see in Corollary~\ref{re:cor:soln_borel-reg} of Theorem~\ref{re:thm:soln_is_Borel_sum} that these solutions are Borel regular.
\begin{Theorem}[restatement of Theorem~\ref{thm:exist_uniq_ODE}]\label{re:thm:exist_uniq_ODE}
Choose a root $-\alpha$ of $P$ where $\tau_\alpha$ is real and positive. Choose an open sector $\Omega_\alpha$ which has an opening angle of $\pi$ or less, has $\zeta = \alpha$ at its tip, and does not touch any other root of $P(-\zeta)$. The equation $\mathcal{P}\Psi = 0$ has a unique solution~$\Psi_\alpha$ in the affine subspace
\smash{${\rm e}^{-\alpha z} \bigl[ z^{-\tau_\alpha} + \dualsingexp{-\tau_\alpha-1}\bigl(\widehat{\Omega}_\alpha^\bullet\bigr) \bigr] $}
of the space ${\rm e}^{-\alpha z} \dualsingexp{-\tau_\alpha}\bigl(\widehat{\Omega}_\alpha^\bullet\bigr)$ in Definition~{\rm\ref{def:H_hat}}.
\end{Theorem}
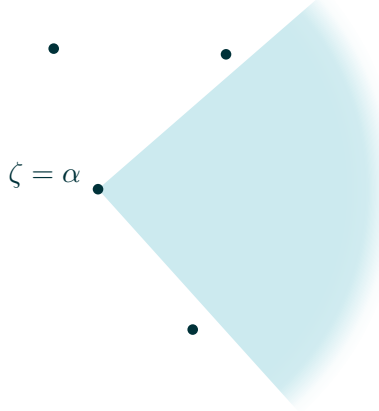
\begin{figure}[!ht]\centering\vspace{-18mm}
\hspace*{-25mm}\begin{tikzpicture}
\newcommand{\spill}{4}
\fill[ietcoast!20, bezier bounding box, path fading=radial edge]
 (-\spill, -\spill) (\spill, \spill)
 (0, 0)--(41:\spill) arc (41:-48:\spill)--cycle;
\fill[ietcoast!33!black] circle (0.7mm) node[anchor=-15, outer sep=1mm] {$\zeta = \alpha$};
\foreach \p in {(107.5:1.95), (46.5:2.46), (-56:2.24)} {
 \fill[ietcoast!33!black] \p circle (0.7mm);
}
\end{tikzpicture}\vspace{-12mm}

\caption{The sectorial domain $\Omega_\alpha$.}\label{fig:sectorial_domain--with roots P}
\end{figure}
\begin{proof}
After rescaling by a constant, \cite[Theorem~4.6]{reg-sing-volterra} tells us that the equation $\hat{\mathcal{P}}_\alpha \psi_\alpha = 0$ has a unique solution $\psi_\alpha$ in the affine subspace
\smash{$ \frac{\zeta^{\tau_\alpha-1}}{\Gamma(\tau_\alpha)} + \singexpalg{\tau_\alpha}(\Omega_\alpha) $}
of the space \smash{$\singexpalg{\tau_\alpha-1}(\Omega_\alpha)$}. The same is true on any smaller sector created by shaving a sector off each edge of $\Omega_\alpha$. By Theorem~\ref{translation} and the results of Section~\ref{sec:reg-decay}, the Laplace transform~\smash{$\laplace^\theta_{\zeta, \alpha}$} is a left-invertible map~\smash{$\singexpalg{\tau_\alpha}(\Omega_\alpha) \to {\rm e}^{-\alpha z} \dualsingexp{-\tau_\alpha-1}\bigl(\widehat{\Omega}_\alpha^\bullet\bigr)$}.
With $\Psi_\alpha = \laplace^\theta_{\zeta, \alpha} \psi_\alpha$, the desired result follows.
\end{proof}

\begin{Theorem}[restatement of Theorem~\ref{thm:soln_is_Borel_sum}]\label{re:thm:soln_is_Borel_sum}
The analytic solution $\Psi_\alpha$ from Theorem~{\rm\ref{re:thm:exist_uniq_ODE}} is the Borel sum of a formal trans-monomial solution $\series{\Psi}_\alpha \in {\rm e}^{-\alpha z} z^{-\tau_\alpha} \C \big\llbracket z^{-1} \big\rrbracket$ of the equation $\mathcal{P}\series{\Psi} = 0$.
\end{Theorem}
\begin{proof}
Take an open sector $\Omega_\alpha$ of the kind used in Theorem~\ref{re:thm:exist_uniq_ODE}, and expand it to an open set~$\smash{\Upsilon \subset \widetilde{\C^\times}}$ of the kind used in Lemma~\ref{lem:shifted_holo_closed}. Our arguments in \cite{reg-sing-volterra} took place on $\C$, but they work just as well on a covering space like \smash{$\widetilde{\C^\times}$}. We can therefore use \cite[Theorem~4.6]{reg-sing-volterra} to show that $\hat{\mathcal{P}}_\alpha \psi_\alpha = 0$ has a unique solution $\psi_\alpha$ in the affine subspace
\smash{$\frac{\zeta_\alpha^{\tau_\alpha-1}}{\Gamma(\tau_\alpha)} + \singexpalg{\tau_\alpha}(\Upsilon) $}
of the space $\singexpalg{\tau_\alpha-1}(\Upsilon)$. Restricting this solution to $\Omega_\alpha$ gives the analogous solution from the proof of Theorem~\ref{re:thm:exist_uniq_ODE}, as we can see from the uniqueness part of \cite[Theorem~4.6]{reg-sing-volterra}. To express~${\Psi_\alpha = \laplace_{\zeta, \alpha} \psi_\alpha}$ as the Borel sum of a trans-monomial, we first need to express $\psi_\alpha$ as the sum of a convergent power series. We will do this by showing that $\psi_\alpha$ lies in $\zeta_\alpha^{\tau_\alpha - 1} \mathcal{O}(\pi \Upsilon)$.

For convenience, let $p = P(-\zeta)$ and $q = Q(-\zeta)$. To prove \cite[Theorem~4.6]{reg-sing-volterra}, we rewrite the equation $\hat{\mathcal{P}}_\alpha \solwhole = 0$ as a regular singular Volterra equation $\solwhole = \volterra^\alpha \solwhole$ that satisfies the conditions of the underlying existence and uniqueness result~\cite[Theorem~1.6]{reg-sing-volterra}. The operator $\volterra^\alpha$ is the sum of the ``prototypical'' part
$ \hardpart^\alpha = -\frac{1}{p} \circ \partial^{-1}_{\zeta, \alpha} \circ q $
and the perturbation
\smash{$ \softpart^\alpha = \frac{1}{p} \circ \partial^{-2}_{\zeta, \alpha} \circ R\bigl(\partial^{-1}_{\zeta, \alpha}\bigr)$}.

Let us run through the proof of Theorem~1.6, specialized to the case we consider in Theorem~4.6. First, picking an arbitrary point $b \in \Omega_\alpha$, we show that the {\em prototype solution}
\[ \solproto(a) = \frac{1}{p(a)} \exp\left(-\int_{b}^{a}\frac{q}{p} \,{\rm d}\zeta_\alpha\right) \]
lies in the space $\singexpalg{\tau_\alpha - 1}(\Upsilon)$ and satisfies the equation $\solproto = \hardpart^\alpha \solproto$. We then look for a perturbation $\solptb$ that makes $\solwhole = \solproto + \solptb$ a solution of the Volterra equation we are trying to solve. This is equivalent to solving the inhomogeneous equation
\begin{equation}\label{eqn:inhomog-volterra}
\solptb = \softpart^\alpha \solproto + \volterra^\alpha \solptb.
\end{equation}
The central idea of the proof is to show that $\softpart^\alpha$ maps $\solproto$, and in fact all of $\singexpalg{\tau_\alpha - 1}(\Upsilon)$, into~$\singexpalg{\tau_\alpha}(\Upsilon)$, and that $\volterra^\alpha$ is a contraction of $\singexp{\tau_\alpha}{\Lambda}(\Upsilon)$ when $\Lambda$ is large enough. It follows, by the contraction mapping theorem, that equation~\eqref{eqn:inhomog-volterra} has a unique solution $\solptb$ in~$\singexp{\tau_\alpha}{\Lambda}(\Upsilon)$. More explicitly, we can solve equation~\eqref{eqn:inhomog-volterra} by fixed-point iteration. Defining
\begin{gather*}
\solptb^{(0)} = \softpart^\alpha \solproto ,\qquad
\solptb^{(1)} = \softpart^\alpha \solproto + \volterra^\alpha \solptb^{(0)} ,\qquad
\solptb^{(2)} = \softpart^\alpha \solproto + \volterra^\alpha \solptb^{(1)}, \\
\solptb^{(3)} = \softpart^\alpha \solproto + \volterra^\alpha \solptb^{(2)}, \qquad
\dots,
\end{gather*}
we get a sequence of functions that converges in $\singexp{\tau_\alpha}{\Lambda}(\Upsilon)$ to a solution $f_*$.

In \cite[Section~3.2.1]{reg-sing-volterra}, we rewrite the prototype solution in the form
\[ \solproto(a) = \frac{1}{p(a)} \left(\frac{\zeta_\alpha(a)}{\zeta_\alpha(b)}\right)^{\tau_\alpha} \exp\left[-\int_b^a \left( \frac{q}{p} + \frac{\tau_\alpha}{\zeta_\alpha} \right) {\rm d}\zeta_\alpha\right], \]
which shows that it represents a germ in $\zeta_\alpha^{\tau_\alpha-1} \mathcal{O}_{\zeta = \alpha}$. By working with Taylor series, we can show that $\softpart^\alpha$ maps
$\zeta_\alpha^{\tau_\alpha-1} \mathcal{O}_{\zeta=\alpha} \to \zeta^{\tau_\alpha} \mathcal{O}_{\zeta=\alpha} $
and $\volterra^\alpha$ maps $\zeta_\alpha^{\tau_\alpha} \mathcal{O}_{\zeta=\alpha}$ to itself. It follows that the sequence \smash{$\solptb^{(0)}, \solptb^{(1)}, \solptb^{(2)}, \dots$} lies in $\zeta_\alpha^{\tau_\alpha} \mathcal{O}_{\zeta=\alpha} \cap \singexp{\tau_\alpha}{\Lambda}(\Upsilon)$, which by Lemma~\ref{lem:shifted_holo_closed} is a closed subspace of \smash{$\singexp{\tau_\alpha}{\Lambda}(\Upsilon)$}. The limit $f_*$, which satisfies equation~\eqref{eqn:inhomog-volterra}, therefore represents a germ in $\zeta_\alpha^{\tau_\alpha} \mathcal{O}_{\zeta = \alpha}$. Recalling that the prototype solution represents a germ in $\zeta_\alpha^{\tau_\alpha-1} \mathcal{O}_{\zeta = \alpha}$, we can deduce that $\psi_\alpha$ represents a germ in $\zeta_\alpha^{\tau_\alpha-1} \mathcal{O}_{\zeta = \alpha}$ too.

We now cross over into the formal world, expanding $\psi_\alpha$ as a shifted Taylor series $\series{\psi}_\alpha \in \zeta_\alpha^{\tau_\alpha-1} \C\{\zeta\}$. By definition, the formal Laplace transform \smash{$\series{\Psi}_\alpha = \laplace_{\zeta, \alpha} \series{\psi}_\alpha$} lies in ${\rm e}^{-\alpha z} z^{-\tau_\alpha} \C \big\llbracket z^{-1} \big\rrbracket$.

When we take the Borel sum of \smash{$\series{\Psi}_\alpha$}, we start by retracing the construction steps described above. The Borel transform inverts the formal Laplace transform, bringing us back to $\series{\psi}_\alpha$, and summation then inverts the Taylor expansion, bringing us back to $\psi_\alpha$. We finish the Borel summation by taking the Laplace transform $\laplace^\theta_{\zeta, \alpha} \psi_\alpha$, giving $\Psi_\alpha$ by definition.
\end{proof}

\begin{Corollary}[restatement of Corollary~\ref{cor:soln_borel-reg}]\label{re:cor:soln_borel-reg}
The analytic solution $\Psi_\alpha$ from Theorem~{\rm\ref{re:thm:exist_uniq_ODE}} is Borel regular. This is because $\Psi_\alpha$ is asymptotic to the formal solution $\series{\Psi}_\alpha$ from Theorem~{\rm\ref{re:thm:soln_is_Borel_sum}}.
\end{Corollary}
\begin{proof}
In the proof of Theorem~\ref{re:thm:soln_is_Borel_sum}, we showed that the position-domain solution $\psi_\alpha \in \singexpalg{\tau_\alpha-1}(\Upsilon)$ is the sum of a shifted Taylor series $\series{\psi}_\alpha \in\smash{ \zeta^{\tau_\alpha-1} \C\{\zeta\}}$. Since $\Psi_\alpha$ is the Laplace transform of $\psi_\alpha$, and \smash{$\series{\Psi}_\alpha$} is the formal Laplace transform of \smash{$\series{\psi}_\alpha$}, Lemma~\ref{lem:laplace-bridge} tells us that \smash{$\Psi_\alpha$} is asymptotic to $\series{\Psi}_\alpha$. Taking the Borel sum of $\series{\Psi}_\alpha$ gives us back $\Psi_\alpha$, by construction, so $\Psi_\alpha$ is Borel regular.
\end{proof}

\subsubsection{A new proof of the Borel summability of the Poincar\'e solutions}\label{sec:new-summability-proof}
In~\cite{int-irreg}, Poincar\'e described a formal way to solve the equation $\mathcal{P}\series{\Phi} = 0$ for an operator of the form~\eqref{eqn:operator-P}. Choose a root $-\alpha$ of $P$ and look for a formal solution $\series{\Psi}_\alpha$ in the space \smash{${\rm e}^{-z\alpha} z^{-\tau_\alpha}\C\big\llbracket z^{-1} \big\rrbracket$}, recalling that \smash{$\tau_\alpha := \frac{Q(-\alpha)}{P'(-\alpha)}$}. Once the leading coefficient of $\series{\Psi}_\alpha$ is chosen, the rest of the coefficients are determined, and can be found order by order. The Ramis index theorem guarantees that $\series{\Psi}_\alpha$ will be $1$-Gevrey~\cite{ramis_index}.

We have come this far in our discussion of Poincar\'e's method through formal reasoning in the frequency domain. From here, there are various ways to show that $\series{\Psi}_\alpha$ is Borel summable. We do it by using \cite[Theorem~4.6]{reg-sing-volterra} to show that the relevant integral equation in the position domain has a unique analytic solution.
\begin{Theorem}\label{thm:summability_ODE}
Choose a root $-\alpha$ of $P$ where $\tau_\alpha$ is real and positive. Choose an open sector~$\Omega_\alpha$ which has an opening angle of $\pi$ or less, has $\zeta = \alpha$ at its tip, and does not touch any other root of~$P(-\zeta)$. If a $1$-Gevrey trans-monomial $\series{\Psi}_\alpha \in {\rm e}^{-z\alpha} z^{-\tau_\alpha}\C\big\llbracket z^{-1} \big\rrbracket_1$ satisfies the equation~\smash{$\mathcal{P} \series{\Psi}_\alpha = 0$}, then it is Borel summable at $\alpha$ along any ray in $\Omega_\alpha$, and its Borel sum is a scalar multiple of the analogous analytic solution $\Psi_\alpha$ from Theorem~{\rm\ref{re:thm:exist_uniq_ODE}}.
\end{Theorem}
\begin{proof}
By Proposition~\ref{prop:frac-der-int-borel}, the Borel transform \smash{$\series{\psi}_\alpha = \borel_\zeta \series{\Psi}_\alpha$} is a formal solution of the equation~\smash{$\hat{\mathcal{P}}_\alpha \series{\psi}_\alpha = 0$}, where $\hat{\mathcal{P}}_\alpha$ is the operator \eqref{eq:hat_P}. By Proposition~\ref{prop:gevrey_to_convergent}, the series $\series{\psi}_\alpha$ converges to a~holomorphic function $\hat{\psi}_\alpha$ on some open set $\Omega_\text{near}$ created by restricting $\Omega_\alpha$ to a small disk around of $\zeta = \alpha$. Since the integrals in $\hat{\mathcal{P}}_\alpha \hat{\psi}_\alpha$ converge absolutely, $\hat{\psi}_\alpha$ satisfies the same equation that its Taylor series does: $\hat{\mathcal{P}}_\alpha \hat{\psi}_\alpha = 0$.

Since the series $\series{\psi}_\alpha$ lies in $\zeta_\alpha^{\tau_\alpha-1} \C\{\zeta_\alpha\}$, the function $\hat{\psi}_\alpha$ lies in $\zeta_\alpha^{\tau_\alpha-1} \mathcal{O}(\Omega_\text{near})$. Since $\Omega_\text{near}$ is bounded, this implies that $\hat{\psi}_\alpha$ lies in $c \zeta_\alpha^{\tau_\alpha-1} + \singexpalg{\tau_\alpha}(\Omega_\text{near})$, where $c$ is the leading coefficient of \smash{$\series{\psi}_\alpha$}. We know from \cite[Theorem~4.6]{reg-sing-volterra} that the equation \smash{$\hat{\mathcal{P}}_\alpha \varphi = 0$} has a unique solution in the affine space \smash{$\zeta_\alpha^{\tau_\alpha-1} + \singexpalg{\tau_\alpha}(\Omega_\alpha)$}. The same theorem applies on the domain $\Omega_\text{near}$, where its uniqueness provision shows that~$\hat{\psi}_\alpha$ and $c\psi_\alpha$ must match. This means, in particular, that~$\hat{\psi}_\alpha$ can be analytically continued throughout $\Omega_\alpha$, and it is uniformly of exponential type on that sector. Thus, $\hat{\psi}_\alpha$ has a well-defined Laplace transform along any ray \smash{$\mathcal{J}^\theta_{\zeta, \alpha}$} in $\Omega_\alpha$, meaning that~$\series{\Psi}_\alpha$ is Borel summable along any such ray. Recalling how the solution $\Psi_\alpha$ was constructed in Theorem~\ref{re:thm:exist_uniq_ODE}, we see that the Borel sum of $\series{\Psi}_\alpha$ is $c\Psi_\alpha$.
\end{proof}

\subsection{Borel regularity for thimble integrals}\label{borel-reg-thimble}
Let us recall the setting from Section~\ref{borel-reg:explanatory-power}. Let $X$ be a $1$-dimensional complex manifold equipped with a volume form $\nu \in \Omega^{1,0}(X)$ and a holomorphic function $f\maps X\to\C$ with non--degenerate critical points. Let $I_a$ be the thimble integral
\begin{equation}\label{exp-int}
I_a := \int_{\mathcal{C}_a^\theta} {\rm e}^{-zf} \nu,
\end{equation}
where $\mathcal{C}_a^\theta$ is the thimble through the critical point $a$. Let $\alpha$ be the associated critical value, $f(a)$. In this subsection, we will prove that thimble integrals of the form \eqref{exp-int} are Borel regular. The first step is to rewrite $I_a$ as the Laplace transform of a function $\iota_a$, which is given explicitly by the well-known formula that we will call the \textit{thimble projection formula}.
\begin{Lemma}[adapted from {\cite[Section~3.3]{pham}}]\label{lem:thimble_proj_formula-proof}
A function $\iota_a$ with $I_a = \laplace_{\zeta, \alpha}^\theta \iota_a$ is given by the {\em thimble projection formula}
\begin{equation}\label{eqn:formula-proof}
\iota_a = \frac{\partial}{\partial \zeta} \bigg( \int_{\mathcal{C}_a^\theta(\zeta)}\nu \bigg),
\end{equation}
where $\mathcal{C}_a^\theta(\zeta)$ is the part of $\mathcal{C}_a^\theta$ that goes through $f^{-1}\bigl(\bigl[\alpha,\zeta {\rm e}^{{\rm i}\theta}\bigr]\bigr)$. Notice that $\mathcal{C}_a^\theta(\zeta)$ starts and ends in $f^{-1}(\zeta)$.
\end{Lemma}
\begin{Remark}
In \cite{pham}, Pham describes this formula in a slightly different form, integrating the variation of $\nu/{\rm d}f$, instead of differentiating the integral as we do in equation~\eqref{eqn:formula-proof}.
\end{Remark}
\begin{proof}[Proof of Lemma~\ref{lem:thimble_proj_formula-proof}]
We first recast the integral $I_a$ into the position domain. As $\zeta$ goes rightward from $\alpha$ with angle $\theta$, the start and end points of $\mathcal{C}_a^\theta(\zeta)$ sweep backward along $\mathcal{C}^-_a(\zeta)$ and forward along $\mathcal{C}^+_a(\zeta)$, respectively. Hence, we have
\begin{align*}
I_a & = \int_{\mathcal{C}_a^\theta} {\rm e}^{-zf} \nu= \int_{\mathcal{J}_{\zeta,\alpha}^\theta}{\rm e}^{-z\zeta} \left[\frac{\nu}{{\rm d}f}\right]_{\operatorname{start} \mathcal{C}_a^\theta(\zeta)}^{\operatorname{end} \mathcal{C}_a^\theta(\zeta)} \,{\rm d}\zeta.
\end{align*}
Noticing that the right-hand side is a Laplace transform, we learn that
\begin{equation}\label{thimble-difference}
{\iota}_a = \left[\frac{\nu}{ {\rm d}f}\right]_{\operatorname{start} \mathcal{C}_a^\theta(\zeta)}^{\operatorname{end} \mathcal{C}_a^\theta(\zeta)}.
\end{equation}
This is equivalent to the expression for $\iota_\alpha$ proposed in equation~\eqref{eqn:formula-proof}.
\end{proof}

We will now prove that $I_a$ is Borel regular.
\begin{Theorem}[restatement of Theorem~\ref{thm:maxim}]\label{thm:maxim-proof}
If the integral defining $I_a$ is absolutely convergent, then $I_a$ is Borel regular. More explicitly,
\begin{enumerate}\itemsep=0pt
\item[$(1)$] 
The function $I_a$ has an asymptotic expansion $\series{I}_a := \aexp^{-\theta} I_a$, which lies in the space ${\rm e}^{-z \alpha} z^{-\frac{1}{2}}\allowbreak\times \C\big\llbracket z^{-1}\big\rrbracket$. Here, $\theta$ is the direction of the ray $\mathcal{J}^\theta_{\zeta, \alpha}$ that defines the thimble and $\alpha=f(a)$.
\item[$(2)$] 
The Borel transform $\series{\iota}_a := \borel_\zeta \series{I}_a$ converges near $\zeta = \alpha$.
\item[$(3)$] 
The sum of $\series{\iota}_a$ can be analytically continued along the ray $\mathcal{J}_{\zeta, \alpha}^\theta$. Its Laplace transform along that ray is well defined, and equal to $I_a$.
\end{enumerate}
\end{Theorem}
\begin{proof}
Since $f$ has non--degenerate critical points, we can find a holomorphic chart $\tau$ around~$a$ with $\frac{1}{2} \tau^2 = f - \alpha$. Let $\mathcal{C}^-_a$ and $\mathcal{C}^+_a$ be the parts of $\mathcal{C}_a^\theta$ that go from the past $\bigl(-{\rm e}^{-{\rm i}\theta}\infty\bigr)$ to $a$ and from $a$ to the future $\bigl(+{\rm e}^{-{\rm i}\theta}\infty\bigr)$, respectively. By changing the sign of $\tau$, we can arrange for it to be valued in $(-\infty, 0]$ and $[0,\infty)$ on $\mathcal{C}^-_a$ and $\mathcal{C}^+_a$, respectively. Since $\nu$ is holomorphic, we can express it as a Taylor series
$\nu = \sum_{m \ge 0} b_m^a \tau^m{\rm d}\tau$
that converges in some disk $|\tau| < \varepsilon$.

We write $\approx$ to say that two functions are asymptotic to all orders. By the method of steepest descent,\footnote{The details can be found in \cite{miller2006applied}: follow the proof of Proposition~2.1 in through equation~(2.9).}
\[ {\rm e}^{z \alpha} I_a \approx \int_{\tau \in [-\varepsilon, \varepsilon]} {\rm e}^{-z\tau^2/2} \nu \]
as $z \to \infty$. Plugging in the Taylor series for $\nu$, we get
\begin{equation*}
{\rm e}^{z \alpha} I_a \approx \int_{-\varepsilon}^\varepsilon {\rm e}^{-z\tau^2/2} \sum_{m \ge 0} b_m^a \tau^m \,{\rm d}\tau = \int_{-\varepsilon}^\varepsilon {\rm e}^{-z\tau^2/2} \sum_{n \ge 0} b_{2n}^a \tau^{2n}\, {\rm d}\tau.
\end{equation*}
Applying the dominated convergence theorem on the right-hand side, we deduce that\footnote{Notice that the sum over $k$ is empty when $n = 0$. Following convention, we extend the double factorial to all odd integers by its recurrence relation, giving $(-1)!! = 1$.}
\begin{align*}
{\rm e}^{z \alpha} I_a & \approx \sum_{n \ge 0} b_{2n}^a \int_{-\varepsilon}^\varepsilon {\rm e}^{-z\tau^2/2} \tau^{2n} \,{\rm d}\tau \\
& = \sum_{n \ge 0} (2n-1)!! b_{2n}^a \Bigg[ \sqrt{2\pi} z^{-(n+1/2)} \operatorname{erf}\bigl(\varepsilon \sqrt{z/2}\bigr) - 2{\rm e}^{-z\varepsilon^2/2} \sum_{\substack{k \in \mathbb{N}_+ \\ k \le n}} \frac{\varepsilon^{2k-1}}{(2k-1)!!} z^{-n+k-1} \Bigg].
\end{align*}

The annoying ${\rm e}^{-z\varepsilon^2/2}$ correction terms are crucial for the convergence of the sum, but we can get rid of them and still have a formal series asymptotic to ${\rm e}^{-z \alpha} I_a$. In other words,
\[ {\rm e}^{z \alpha} I_a \sim \sqrt{2\pi} \sum_{n \ge 0} (2n-1)!! b_{2n}^a z^{-(n+1/2)} \operatorname{erf}\bigl(\varepsilon \sqrt{z/2}\bigr). \]
To see why, cut the sum off after $N$ terms, and observe that
\begin{align*}
 \Bigg | \sum_{n= 0}^N(2n-1)!! b_{2n}^a 2{\rm e}^{-z\varepsilon^2/2} \sum_{\substack{k \in \mathbb{N}_+ \\ k \le n}} \frac{\varepsilon^{2k-1}}{(2k-1)!!} z^{-n+k-1} \Bigg | &\leq 2 {\rm e}^{-\Re (z) \varepsilon^2/2} \sum_{n=0}^N(2n-1)!! b_{2n}^a n|z|^{-n}.
\end{align*}
The right-hand side goes to $0$ as $\Re(z)$ grows. The differences $1 - \operatorname{erf}\bigl(\varepsilon \sqrt{z/2}\bigr)$ shrink exponentially as $z$ grows \cite[inequality~(5)]{chiani-dardari-book}, allowing the simpler estimate
\[ {\rm e}^{z\alpha} I_a \sim \sqrt{2\pi} \sum_{n \ge 0} (2n-1)!! b_{2n}^a z^{-(n+1/2)}. \]
Hence, defining the formal series $\series{I}_a$
\[\series{I}_a := {\rm e}^{-z\alpha} z^{-1/2} \sqrt{2\pi} \sum_{n \ge 0} (2n-1)!! b_{2n}^a z^{-n}\]
we get $\series{I}_a \in {\rm e}^{-z\alpha}z^{-1/2}\C\big\llbracket z^{-1}\big\rrbracket$. This ends the proof of part~(1).

Using the properties of the Borel transform acting on trans-mononials (see Section~\ref{sec:action_transseries}) we get
\begin{equation*}
\series{\iota}_a = \sqrt{2\pi} \sum_{n \ge 0} (2n - 1)!! b_{2n}^a \frac{\zeta_\alpha^{n+\frac{1}{2}}}{\Gamma(n+\frac{1}{2})} = \sqrt{2\pi} \sum_{n \ge 0} \frac{2^n}{\sqrt{\pi}} b_{2n}^a \zeta_\alpha^{n+\frac{1}{2}}.
\end{equation*}
We know from the definition of $\varepsilon$ that $\left|b_n^a\right| \varepsilon^n \lesssim 1$, thus we deduce that $\series{\iota}$ has a finite radius of convergence. This ends the proof of part~(2).

The proof of part~(3) relies on the thimble projection formula \eqref{eqn:formula-proof}: we will show that the Taylor expansion of $\iota_a$ at $\zeta=\alpha$ agrees with $\series{\iota}_a$. We can rewrite the Taylor series for $\nu$ as
\begin{align*}
\nu & = \sum_{n \ge 0} b_n^a [2(f - \alpha)]^{(n - 1)/2} {\rm d}f,
\end{align*}
taking the positive branch of the square root on $\mathcal{C}^+_a$ and the negative branch on $\mathcal{C}^-_a$. Plugging this into our expression for $\iota_a$ (in Lemma~\ref{lem:thimble_proj_formula-proof}), we learn that
\begin{align*}
\iota_a & = \bigg[ \sum_{m \ge 0} b_m^a [2(f - \alpha)]^{(m - 1)/2} \bigg]_{\operatorname{start} \mathcal{C}_a(\zeta)}^{\operatorname{end} \mathcal{C}_a(\zeta)} = \sum_{m \ge 0} b_m^a \bigl( [2\zeta_\alpha]^{(m - 1)/2} - (-1)^{m-1}[2\zeta_\alpha]^{(m - 1)/2} \bigr) \\
& = \sum_{n \ge 0} 2 b_{2n}^a [2\zeta_\alpha]^{n - 1/2} = \sum_{n \ge 0} 2^{n+1/2} b_{2n}^a \zeta_\alpha^{n - 1/2}.
\end{align*}
In particular, $\series{\iota}_a$ can be analytically continued along $\mathcal{J}_{\zeta,\alpha}^\theta$ and its sum is given by~$\iota_a$.
\end{proof}

\begin{Remark}\label{rmk:1/2-deriv}
For some applications, it is more convenient to work with $\Phi_a := z^{-1/2} I_a$, whose asymptotic series has integer powers. The argument we use to prove of Theorem~\ref{thm:maxim-proof} can also be used to show that $\Phi_a$ is Borel regular, but the thimble projection formula is trickier in this case. Proposition~\ref{prop:L-int-op} tells us that $\Phi_a = \laplace^\theta_{\zeta, \alpha} \phi_a$ for
\smash{$\phi_a := \fracderiv{-1/2}{\zeta}{\alpha} \iota_a$}.
To get the thimble projection formula for $\phi_a$, take the $1/2$ fractional integral of both sides:
\begin{align*}
\fracderiv{-1/2}{\zeta}{\alpha} \phi_a & = \fracderiv{-1}{\zeta}{\alpha} \iota_a = \fracderiv{-1}{\zeta}{\alpha} \frac{\partial}{\partial \zeta} \left( \int_{\mathcal{C}_a^\theta(\zeta)}\nu \right) = \int_{\mathcal{C}_a^\theta(\zeta)}\nu - \int_{\mathcal{C}_a^\theta(\alpha)}\nu.
\end{align*}
The second term vanishes because $\mathcal{C}_a^\theta(\alpha)$ is a single point, leaving
\smash{$ \fracderiv{-1/2}{\zeta}{\alpha} \phi_a = \int_{\mathcal{C}_a^\theta(\zeta)}\nu$}.
The $1/2$ fractional derivative is a left inverse of the $1/2$ fractional integral. Applying it to both sides, we find that
\smash{$
\phi_a = \fracderiv{1/2}{\zeta}{\alpha} \bigl( \int_{\mathcal{C}_a^\theta(\zeta)}\nu \bigr)$}.
\end{Remark}
\begin{Remark}
We can use the thimble projection formula to show that each function $\iota_\alpha$ is resurgent, in the sense of \'{E}calle~\cite[Section~1]{EcalleI}. Let $\mathcal{J}_\alpha(\beta)$ be the straight path from $\zeta = \alpha$ to~${\zeta = \beta}$ in $\C$. As long as this path avoids the critical values of $f$, it lifts uniquely along $f$ to a~path $\mathcal{C}_a^\theta(\beta)$ in $X$. This lets us analytically continue \smash{$\int_{\mathcal{C}_a^\theta(\zeta)} \nu$} to a star-shaped domain $\Omega_\alpha \subset \C$. Intuitively, $\Omega_\alpha$ is constructed by drawing rays from $\zeta = \alpha$ through all of the other critical values, and then cutting $\C$ along the parts of those rays that go from the other critical values to $\zeta = \infty$, as shown in Figure~\ref{Fig:slit domain}.
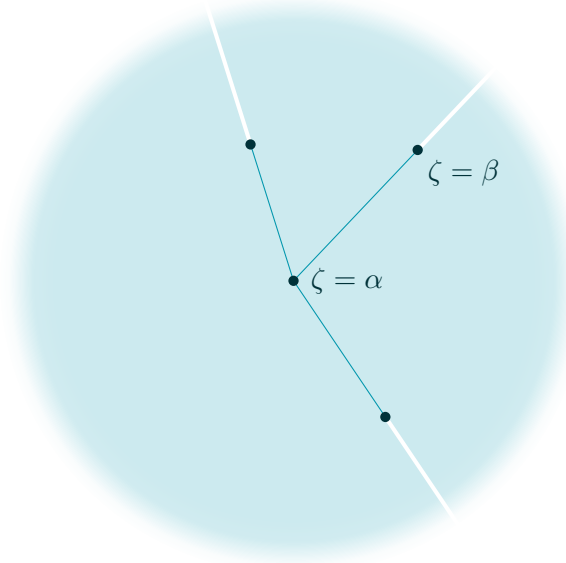
\begin{figure}[!ht]\centering\vspace{-2mm}

\begin{tikzpicture}[scale=0.97]
\newcommand{\spill}{4}
\newcommand{\slitsize}{0.4}
\fill[ietcoast!20, bezier bounding box, path fading=radial edge]
 (-\spill, -\spill) (\spill, \spill)
 (0, 0)--(180:\spill)
 arc (180:107.5+\slitsize:\spill)--(107.5+\spill/1.95*\slitsize:1.95)--(107.5-\spill/1.95*\slitsize:1.95)--(107.5-\slitsize:\spill)
 arc (107.5-\slitsize:46.5+\slitsize:\spill)--(46.5+\spill/2.46*\slitsize:2.46)--(46.5-\spill/2.46*\slitsize:2.46)--(46.5-\slitsize:\spill)
 arc (46.5-\slitsize:-56+\slitsize:\spill)--(-56+\spill/2.24*\slitsize:2.24)--(-56-\spill/2.24*\slitsize:2.24)--(-56-\slitsize:\spill)
 arc (-56-\slitsize:-180:\spill)--(0, 0);
\node[ietcoast!33!black, anchor=north west] at (46.5:2.46) {$\zeta = \beta$};
\foreach \p in {(107.5:1.95), (46.5:2.46), (-56:2.24)} {
 \draw[ietcoast] \p--(0,0);
 \fill[ietcoast!33!black] \p circle (0.7mm);
}
\fill[ietcoast!33!black] circle (0.7mm) node[anchor=180, outer sep=1mm] {$\zeta = \alpha$};
\end{tikzpicture}

\vspace{-1mm}

\caption{The domain $\Omega_\alpha$.}\label{Fig:slit domain}
\end{figure}

Initially, the function \smash{$\int_{\mathcal{C}_a(\zeta)} \nu$} is only defined for $\zeta$ on the ray $[\alpha, \infty)$, but it can be analytically continued to points off the ray by homotopy of the path $\mathcal{C}_a(\zeta)$. Since $\nu$ is holomorphic, the only way for \smash{$\int_{\mathcal{C}^\theta_a(\zeta)} \nu$} to become singular is for something to go wrong with the homotopy of $\mathcal{C}_a(\zeta)$, which can only happen at the critical values of $f$. Therefore, $\mathcal{C}_a(\zeta)$ is endlessly analytically continuable away from the critical values of $f$. Resurgent functions have resurgent derivatives,\footnote{This is because analytically continuing a function to a given open set automatically continues its derivative to the same open set.} so $\iota_\alpha$ is resurgent too.
\end{Remark}
\begin{Remark}\label{rmk:Pham formula}
When $f$ and $\nu$ are polynomials and $X$ is $N$-dimensional, a result of Pham~\cite[equation~(2.4), 2\`{e}me Partie]{pham} gives us an explicit expression for the singularity of the function $\iota_a$ from Lemma~\ref{lem:thimble_proj_formula-proof}. For some coefficients $c_k$ which depend on~$f$, $\nu$, and the thimble $\mathcal{C}_a^\theta$, we have
\[
\iota_a \in (-1)^{\frac{N(N-1)}{2}} \sum_{k=0}^\infty c_k \delta_{\zeta, \alpha}^{(-N/2 - k)} + \mathcal{O}_{\zeta = \alpha},
\]
where the singularities
\begin{equation*}
\delta_{\zeta,\alpha}^{(-\ell)}:=\begin{cases}
\displaystyle\frac{\zeta_\alpha^{\ell-1}}{-2 \pi {\rm i}(\ell-1)!}\log(\zeta_\alpha)+ \mathcal{O}_{\zeta = \alpha} & \text{ if } \ell \in \Z_{> 0},\vspace{1mm}\\
\displaystyle\frac{(-\zeta_\alpha)^{\ell-1}}{2 \pi {\rm i} \Gamma(-\ell)}+ \mathcal{O}_{\zeta = \alpha} & \text{ if } \ell \notin \Z_{> 0} .
\end{cases}
\end{equation*}
can be seen as fractional derivatives of the Dirac delta distribution. In particular, if $N=1$, the inverse Laplace transform of the function $I_a$ from equation~\eqref{exp-int} has a singularity of the form~\smash{$\zeta_\alpha^{1/2}$}, as we compute explicitly in the proof of Theorem~\ref{thm:maxim-proof}.
\end{Remark}
\section{Examples}\label{sec:examples}
\subsection{Airy}\label{sec:airy}
The Airy equation is
\[
\left[\left(\frac{\partial}{\partial y}\right)^2 - y\right] \Phi = 0.
\]
One solution is given by the Airy function,
\[ \Ai(y) = \frac{1}{2 \pi {\rm i}} \int_\gamma \exp\left(\frac{1}{3}t^3 - yt\right) {\rm d}t, \]
where $\gamma$ is a path that comes from $\infty$ at $-60^\circ$ and returns to $\infty$ at $60^\circ$.

\begin{figure}[!h]\centering
\includegraphics[width=6.5cm]{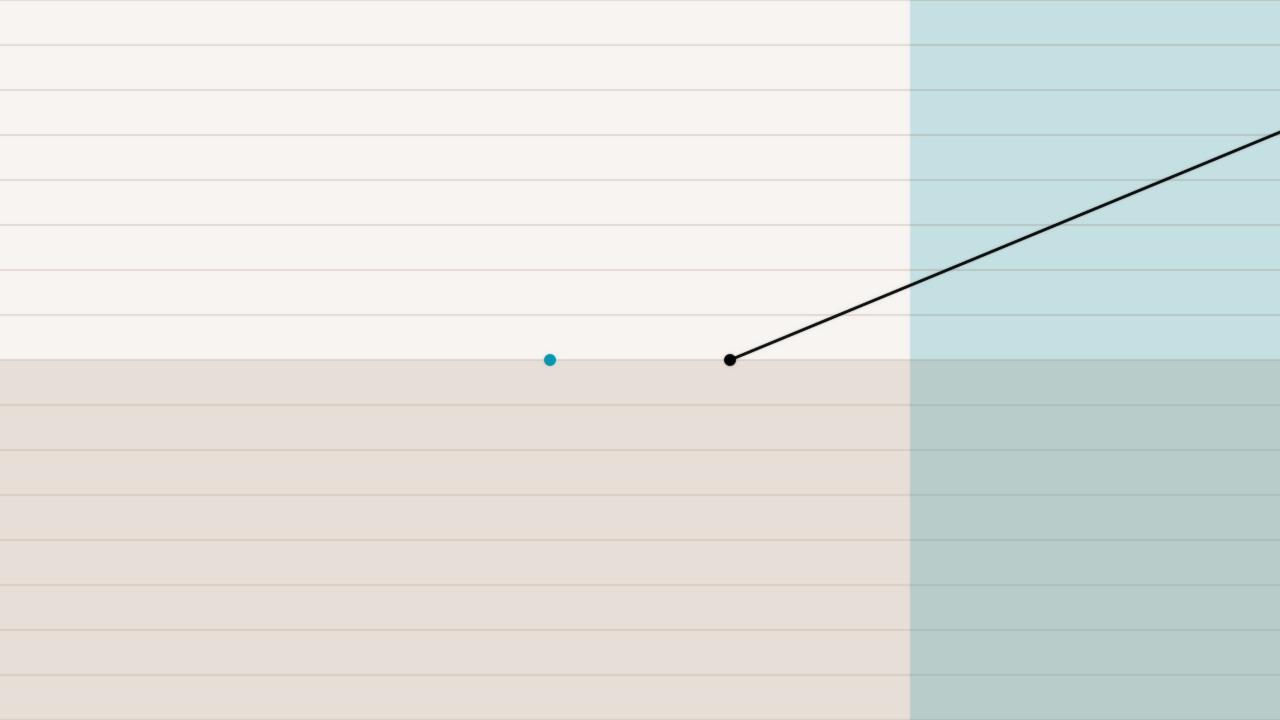}
\qquad
\includegraphics[width=6.5cm]{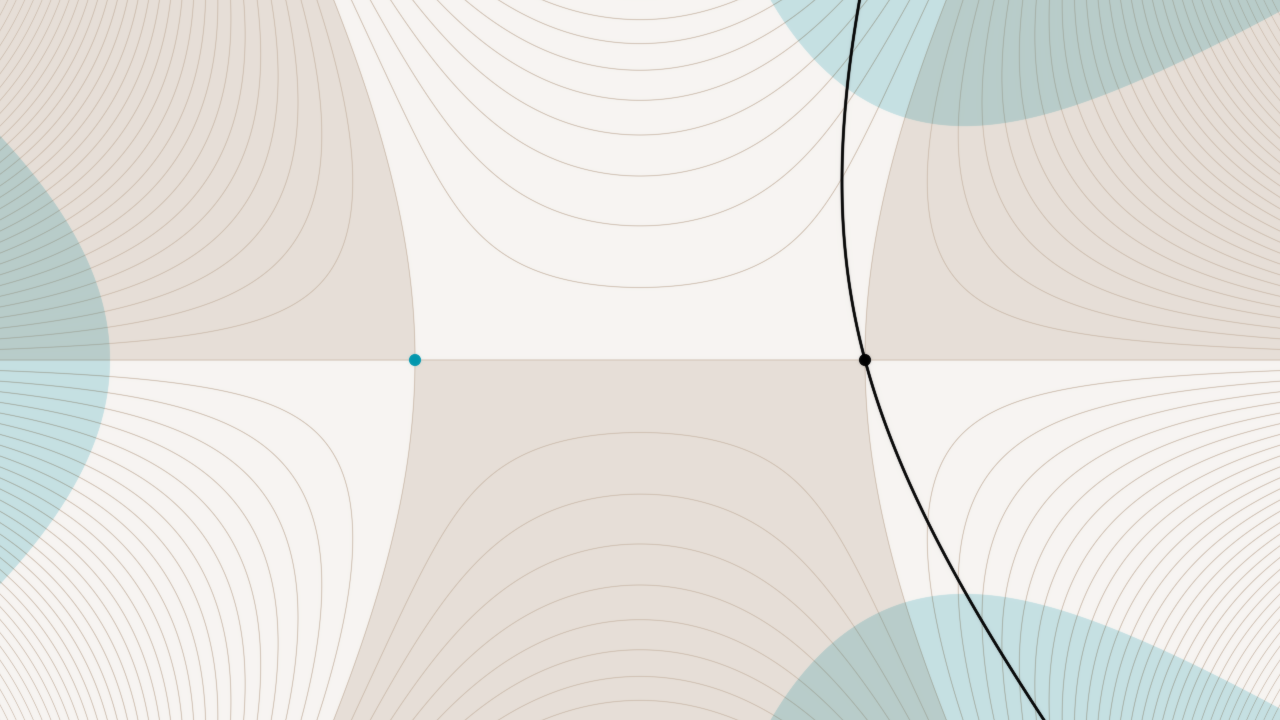}
\caption{The ray \smash{$\mathcal{J}^{\pi/8}_{\zeta, 1}$} in the $\zeta$ plane and its preimage \smash{$\mathcal{C}^{\pi/8}_1$} in the $u$ plane. The upper and lower halves of the $\zeta$ plane are coloured light and dark, so their preimages form a checker pattern on the $u$ plane. The region $\Re(\zeta) > 3$ and its preimage in the $u$ plane are shaded.}\label{fig:path_Airy}
\end{figure}

With the substitution $t = 2uy^{1/2}$, we can rewrite the Airy integral as
\begin{equation}\label{eqn:airy-distillation}
\Ai(y) = y^{1/2} \frac{1}{ \pi {\rm i}} \int_{\mathcal{C}^\theta_1} \exp\left[\frac{2}{3}y^{3/2} \bigl(4u^3 - 3u\bigr)\right] {\rm d}u,
\end{equation}
where $\mathcal{C}^\theta_1$ is the contour that passes through the point $u = \frac{1}{2}$ and projects to the ray \smash{$\mathcal{J}^\theta_{\zeta, 1}$} under the mapping $-\zeta = 4u^3 - 3u$. The contour and its projection are shown in Figure~\ref{fig:path_Airy}.

Note that $4u^3 - 3u$ is the third Chebyshev polynomial. By considering other Chebyshev polynomials, we can situate the Airy function within the family of {\em Airy--Lucas functions}. Treating these functions as a family adds more insight than complexity, so we will go straight to the general case. However, since the Airy function is a classic example in the study of Borel summation and resurgence, it may be worth seeing on its own. In Appendix~\ref{airy-appendix}, we give a detailed treatment of the Airy function, specializing our arguments.

\subsection{Airy--Lucas}\label{example_AL}
\begin{figure}[!ht]\centering
\includegraphics[width=10cm]{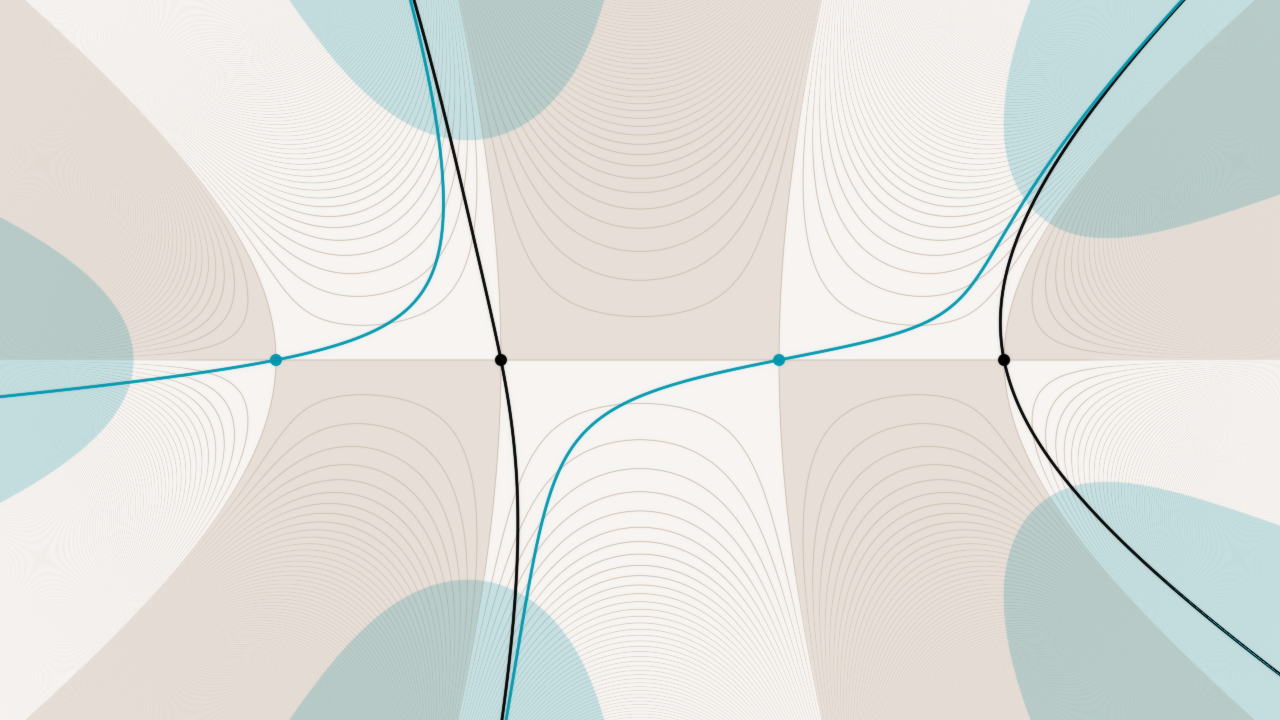}
\caption{The integration contours \smash{$\mathcal{C}^{\pi/8}_1, \dots, \mathcal{C}^{\pi/8}_4$} for the Airy--Lucas functions with $n=5$. The dark and light thimbles are the preimages of the rays \smash{$\mathcal{J}^{\pi/8}_{\zeta, 1}$} and \smash{$\mathcal{J}^{\pi/8}_{\zeta, -1}$}, respectively.}\label{fig:thimble_n_5}
\end{figure}

The Airy--Lucas equation is
\[
\left[\left(\frac{\partial}{\partial y}\right)^2 - (m-1) y^{-1} \frac{\partial}{\partial y} - y^{n-2}\right] \Phi = 0
\]
with $n \in \{3, 4, 5, \dots\}$ and $m \in \{1, 2, \dots, n-1\}$. A few solutions, indexed by $j \in \{1, \dots, n-1\}$, are given by the Airy--Lucas functions~\cite[equation~(3.6)]{charbonnier22}
\begin{equation}\label{integral:AL}
\widehat{\Ai}^{(j)}_{n, m-1}(y) = \left\{\begin{matrix}i & j \text{ odd} \\ 1 & j \text{ even}\end{matrix}\right\} \frac{y^{m/2}}{\pi} \int_{\mathcal{C}^\theta_j} \exp\left[\frac{2}{n} y^{n/2} T_n(u)\right] U_{m-1}(u) \,{\rm d}u,
\end{equation}
where $\mathcal{C}^\theta_j$ is the contour that passes through the point $u = \cos\bigl(\frac{j}{n}\pi\bigr)$ and projects to the ray~\smash{$\mathcal{J}^\theta_{\zeta, \pm 1}$} under the mapping $-\zeta = T_n(u)$. The ray starts at the critical value $1$ when $j$ is odd, and $-1$ when $j$ is even. To ensure that the integral converges, $\theta$ should be near $-\frac{n}{2} \arg(y)$. Figure~\ref{fig:thimble_n_5} shows the contours $\mathcal{C}^\theta_j$, at a suitable $\theta$, in the $n = 5$ case.

\subsubsection[Rewriting as a modified Bessel m/n equation]{Rewriting as a modified Bessel $\boldsymbol{\frac{m}{n}}$ equation}
We can distill the most interesting parts of the Airy--Lucas function by writing
\[ \widehat{\Ai}^{(j)}_{n, m-1}(y) = \frac{2 \sinh\bigl(\frac{m}{n} {\rm i}\pi\bigr)}{n\pi} \left\{\begin{matrix}{\rm i} & j \text{ odd} \\ 1 & j \text{ even}\end{matrix}\right\} y^{{m/2}} K^{(j)}_{m/n}\left(\frac{2}{n} y^{n/2}\right), \]
where
\begin{equation}\label{integral:mod-bessel-rational-AL}
K^{(j)}_{m/n}(z) = \frac{n}{2 \sinh\bigl(\frac{m}{n} {\rm i}\pi\bigr)}\int_{\mathcal{C}^\theta_j} \exp[z T_n(u)] U_{m-1}(u)\,{\rm d}u.
\end{equation}
Saying that \smash{$\widehat{\Ai}^{(j)}_{n, m-1}$} satisfies the Airy--Lucas equation is equivalent to saying that \smash{$K^{(j)}_{m/n}$} satisfies the modified Bessel equation with parameter $\frac{m}{n}$
\begin{equation}\label{eqn:mod-bessel-AL}
\left[z^2 \left(\frac{\partial}{\partial z}\right)^2 + z \frac{\partial}{\partial z} - \left[\left(\frac{m}{n} \right)^2 + z^2\right]\right] \Phi = 0.
\end{equation}
Let us put equation~\eqref{eqn:mod-bessel-AL} in the form $\mathcal{P}\Phi=0$ with $\mathcal{P}$ as in equation~\eqref{eqn:operator-P}
\begin{equation}\label{eqn:reg-mod-bessel-AL}
\left[ \left[ \left(\frac{\partial}{\partial z}\right)^2 - 1 \right] + z^{-1} \frac{\partial}{\partial z} - \left({\frac{m}{n}}\right)^2 z^{-2} \right] \Phi = 0.
\end{equation}
\subsubsection{Asymptotic analysis}\label{sec:asympt-AL}
From the work of Poincar\'e and the Ramis index theorem, as discussed at the beginning of Section~\ref{sec:new-summability-proof}, we know that equation~\eqref{eqn:reg-mod-bessel-AL} has a frame of formal $1$-Gevrey trans-monomial solutions
$\big\{ {\rm e}^{-\alpha z} z^{-\tau_\alpha} \series{W}_\alpha \mid \alpha^2 - 1 = 0 \big\}$,
where $\tau_\alpha = 1/2$ and $\series{W}_\alpha\in\C\big\llbracket z^{-1} \big\rrbracket_1$. From \cite[equations~(10.40.2) and (10.17.1)]{dlmf}, we learn that \smash{$K_{m/n} \sim \bigl(\frac{\pi}{2}\bigr)^{1/2} {\rm e}^{-z} z^{-1/2} \series{W}_1$}, with
\begin{gather*}
\series{W}_1 = 1 - \frac{\bigl(\frac{1}{2}\hspace{-0.5pt}-\hspace{-0.5pt}\frac{m}{n}\bigr)_1 \bigl(\frac{1}{2}\hspace{-0.5pt}+\hspace{-0.5pt}\frac{m}{n}\bigr)_1}{2^1 \cdot 1!} z^{-1} + \frac{\bigl(\frac{1}{2}\hspace{-0.5pt}-\hspace{-0.5pt}\frac{m}{n}\bigr)_2 \bigl(\frac{1}{2}\hspace{-0.5pt}+\hspace{-0.5pt}\frac{m}{n}\bigr)_2}{2^2 \cdot 2!} z^{-2} - \frac{\bigl(\frac{1}{2}\hspace{-0.5pt}+\hspace{-0.5pt}\frac{m}{n}\bigr)_3 \bigl(\frac{1}{2}\hspace{-0.5pt}+\hspace{-0.5pt}\frac{m}{n}\bigr)_3}{2^3 \cdot 3!} z^{-3} + \cdots\!.
\end{gather*}
The holomorphic analysis in Section~\ref{big-idea} will give us holomorphic solutions
$\big\{ {\rm e}^{-\alpha z} z^{-\tau_\alpha} W_\alpha \mid \alpha^2 - 1 = 0 \big\}$,
which seem analogous to the trans-monomials above. Borel summation makes the analogy precise. We will see in Section~\ref{bessel-regularity-AL} that each $z^{-\tau_\alpha} W_\alpha$ is a scalar multiple of the Borel sum of~$z^{-\tau_\alpha} \series{W}_\alpha$, as required by Theorem~\ref{thm:summability_ODE}.
\subsubsection{The big idea}\label{big-idea}
We are going to look for functions $v_\alpha$ whose Laplace transforms $\laplace_{\zeta, \alpha} v_\alpha$ satisfy equation~\eqref{eqn:reg-mod-bessel-AL}. We will succeed when $\alpha^2 - 1 = 0$, and we will see that $K_{m/n}$ is a scalar multiple of $\laplace_{\zeta, 1} v_1$.

By the properties of the Laplace transform laid out in Proposition~\ref{prop:L-int-op}, a function $\laplace_{\zeta, \alpha} v$ satisfies the differential equation~\eqref{eqn:reg-mod-bessel-AL} if and only if $v$ satisfies the integral equation
\begin{equation}\label{int-eq:spatial-mod-bessel-AL}
\left[ \bigl[ \zeta^2 - 1 \bigr] - \fracderiv{-1}{\zeta}{\alpha} \circ \zeta - \left(\frac{m}{n}\right)^2 \fracderiv{-2}{\zeta}{\alpha} \right] v = 0.
\end{equation}
It is tempting to differentiate both sides of this equation until we get
\begin{equation}\label{diff-eq:spatial-mod-bessel-AL}
\left[ \left(\frac{\partial}{\partial \zeta}\right)^2 \circ \bigl[ \zeta^2 - 1 \bigr] - \frac{\partial}{\partial \zeta} \circ \zeta - \left(\frac{m}{n}\right)^2 \right] v = 0,
\end{equation}
which is easier to solve. Unfortunately, an analytic solution of equation~\eqref{diff-eq:spatial-mod-bessel-AL} will not satisfy equation~\eqref{int-eq:spatial-mod-bessel-AL} in general. However, as we show in Appendix~\ref{shifting}, a solution of equation~\eqref{diff-eq:spatial-mod-bessel-AL} {\em will} satisfy equation~\eqref{int-eq:spatial-mod-bessel-AL} if it belongs to $(\zeta - \alpha)^\sigma \mathcal{O}_{\zeta = \alpha}$ for some $\sigma \in (-1, 0)$.

This is good news, because equation~\eqref{diff-eq:spatial-mod-bessel-AL} has a regular singularity at each root of $\zeta^2 - 1$, and the Frobenius method often gives a solution of the desired kind at each regular singular point. We can see the regular singularities by rewriting equation~\eqref{diff-eq:spatial-mod-bessel-AL} as follows
\[ \left[ \bigl(\zeta^2 - 1\bigr) \left(\frac{\partial}{\partial \zeta}\right)^2 + 3\zeta \frac{\partial}{\partial \zeta} + \left[ 1 - \left(\frac{m}{n}\right)^2 \right] \right] v = 0. \]

In Sections \ref{pos-root-AL}--\ref{neg-root-AL}, we will see this approach succeed. For each root $\alpha$, we will find a~solution~$v_\alpha$ of equation~\eqref{diff-eq:spatial-mod-bessel-AL} that belongs to $(\zeta - \alpha)^{\tau_\alpha-1} \mathcal{O}_{\zeta = \alpha}$ for some constant $\tau_\alpha \in (0, 1)$. We will express $v_\alpha$ explicitly enough to see that it extends to a function in $\singexpalg{\tau_\alpha-1}(\Omega_\alpha)$ on any sector $\Omega_\alpha$ that has its tip at $\zeta = \alpha$ and does not touch the other roots (see Figure~\ref{fig:sectorial_domain--with roots P}). We will also be able to see that $v_\alpha$ is uniformly of exponential type $\Lambda$ for any $\Lambda > 0$. It follows by Proposition~\ref{prop:laplace-cont} that for any ray \smash{$\mathcal{J}^\theta_{\zeta, \alpha}$} in $\Omega_\alpha$, the Laplace transform \smash{$\laplace^\theta_{\zeta, \alpha} v_\alpha$} is a well-defined element of \smash{$\dualsingexp{-\tau_\alpha}\bigl(\widehat{\Omega}_\alpha^\Lambda\bigr)$} that satisfies equation~\eqref{eqn:reg-mod-bessel-AL}. Using the technique from the proof of Lemma~\ref{lem:laplace-bridge}, we can see more precisely that $\laplace^\theta_{\zeta, \alpha} v_\alpha$ belongs to
\[ cz^{-\tau_\alpha} + \dualsingexp{-\tau_\alpha-1}\bigl(\widehat{\Omega}_\alpha^\Lambda\bigr) \]
for some non-zero constant $c$, confirming the existence part of Theorem~\ref{re:thm:exist_uniq_ODE}.

We know from Section~\ref{sec:change-translation} that
\smash{$\laplace^\theta_{\zeta, \alpha} v_\alpha = {\rm e}^{-\alpha z} V_\alpha$}
for \smash{$V_\alpha := \laplace^\theta_{\zeta_\alpha, 0} v_\alpha$} and $\zeta = \alpha + \zeta_\alpha$. Observing that multiplication by $z^{-\tau_\alpha}$ gives an isometry \smash{$\dualsingexp{0}\bigl(\widehat{\Omega}_\alpha^\Lambda\bigr) \to \dualsingexp{-\tau_\alpha}\bigl(\widehat{\Omega}_\alpha^\Lambda\bigr)$}, we can do the further decomposition
\smash{$ \laplace_{\zeta, \alpha} v_\alpha = {\rm e}^{-\alpha z} z^{-\tau_\alpha} W_\alpha$},
where $W_\alpha$ is a bounded holomorphic function on~\smash{$\widehat{\Omega}_\alpha^\Lambda$}.
\subsubsection[Focus on zeta = 1]{Focus on $\boldsymbol{\zeta = 1}$}\label{pos-root-AL}
Let us find a germ in $(\zeta-1)^{\tau_1-1} \mathcal{O}_{\zeta=1}$, for some $\tau_1 \in (0, 1)$, that satisfies equation~\eqref{diff-eq:spatial-mod-bessel-AL}. Define a new coordinate $\zeta_1$ on $\C$ so that $\zeta = 1 + \zeta_1$. In this coordinate, equation~\eqref{diff-eq:spatial-mod-bessel-AL} looks like
\[
\left[\zeta_1(2 + \zeta_1) \left(\frac{\partial}{\partial \zeta_1}\right)^2 + 3(1 + \zeta_1) \frac{\partial}{\partial \zeta_1} + \left[1 - \left(\frac{m}{n}\right)^2\right]\right] v = 0.
\]
With another change of coordinate, given by $\zeta_1 = -2\xi_1$, we can rewrite equation~\eqref{diff-eq:spatial-mod-bessel-AL} as the hypergeometric equation
\[
\left[\xi_1 (1 - \xi_1) \left(\frac{\partial}{\partial \xi_1}\right)^2 + 3(\frac{1}{2} - \xi_1) \frac{\partial}{\partial \xi_1} - \left[1 - \left(\frac{m}{n}\right)^2\right]\right] v = 0.
\]
Looking through the twenty-four expressions for Kummer's six solutions, we find one \cite[formula~(15.10.12)]{dlmf} that represents a germ in \smash{$\xi_1^{-1/2} \mathcal{O}_{\xi_1=0}$}
\begin{gather*}
v_1 = \xi_1^{-1/2} {}_2F_1\left(\frac{1}{2}-\frac{m}{n}, \frac{1}{2}+\frac{m}{n}; \frac{1}{2}; \xi_1\right) = -{\rm i}\sqrt{2} \zeta_1^{-1/2} {}_2F_1\left(\frac{1}{2}-\frac{m}{n}, \frac{1}{2}+\frac{m}{n}; \frac{1}{2}; -\frac{1}{2}\zeta_1\right).
\end{gather*}

From the expression above, we can see that $v_1$ has only two singularities in the complex plane: one at $\zeta_1 = 0$ from the factor of \smash{$\zeta_1^{-1/2}$}, and one at $\zeta_1 = -2$ from the hypergeometric function. Therefore, $v_1$ is holomorphic throughout the sector $\Omega_1 = \C \smallsetminus \mathcal{J}^\pi_{\zeta_1, 0}$. Since $v_1$ has a~power law singularity at $\infty$, it is uniformly of exponential type $\Lambda$ for any $\Lambda > 0$. Therefore, it belongs to the space $\singexp{-1/2}{\Lambda}(\Omega_1)$ for all $\Lambda > 0$. By thinking about its shifted Taylor expansion around~${\zeta_1 = 0}$, we can see more precisely that $v_1$ belongs to
\smash{$ -{\rm i}\sqrt{2} \zeta_1^{-1/2} + \singexp{1/2}{\Lambda}(\Omega_1) $}
for all~${\Lambda > 0}$. As discussed in Section~\ref{big-idea}, we conclude that $v_1$ has a well-defined Laplace transform along any ray from $\zeta = 1$ except the one going directly left. Its Laplace transform can be written as~${\laplace^\theta_{\zeta, 1} v_1 = {\rm e}^{-z} z^{-1/2} W_1}$,
where $W_1$ is a holomorphic function on the slit plane~${z \notin (-\infty, 0]}$ which is bounded outside any constant-radius neighborhood of the ray $z \in (-\infty, 0]$.
\subsubsection[Focus on zeta = -1]{Focus on $\boldsymbol{\zeta = -1}$}\label{neg-root-AL}
Let us find a germ in $(\zeta+1)^{\tau_{-1}-1} \mathcal{O}_{\zeta=-1}$, for some $\tau_{-1} \in (0, 1)$, that satisfies equation~\eqref{diff-eq:spatial-mod-bessel-AL}. In the rescaled coordinate from Section~\ref{pos-root-AL}, this is the point $\xi_1 = 1$. Looking again through Kummer's table of solutions, we find another expression \cite[formula~(15.10.14)]{dlmf} that represents a germ in $(1-\xi_1)^{-1/2} \mathcal{O}_{\xi_1=1}$
\begin{align*}
v_{-1}& = (1-\xi_1)^{-1/2} {}_2F_1\left(\frac{1}{2}-\frac{m}{n}, \frac{1}{2}+\frac{m}{n}; \frac{1}{2}; 1-\xi_1\right)\\
&= \sqrt{2} \zeta_{-1}^{-1/2} {}_2F_1\left(\frac{1}{2}-\frac{m}{n}, \frac{1}{2}+\frac{m}{n}; \frac{1}{2}; \frac{1}{2}\zeta_{-1}\right),
\end{align*}
where $\zeta_{-1}$ is the coordinate with $\zeta = -1 + \zeta_{-1}$. By the same reasoning as before, $v_{-1}$ is holomorphic throughout the sector $\Omega_{-1} = \C \smallsetminus \mathcal{J}^0_{\zeta_{-1}, 0}$, and it belongs to the space
\[ \sqrt{2} \zeta_{-1}^{-1/2} + \singexp{1/2}{\Lambda}(\Omega_1) \]
for all $\Lambda > 0$. It has a well-defined Laplace transform along any ray from $\zeta = -1$ except the one going directly right, and its Laplace transform can be written as
$\laplace^\theta_{\zeta, -1} v_1 ={
\rm e}^z z^{-1/2} W_{-1}$,
where $W_{-1}$ is a holomorphic function on the slit plane $z \notin [0, \infty)$ which is bounded outside any constant-radius neighborhood of the ray $z \in [0, \infty)$.

In this example, $v_1$ and $v_{-1}$ happen to be related by a symmetry: the M\"{o}bius transformation that pulls $\zeta$ back to $-\zeta$. Kummer's solutions typically come from six different hypergeometric equations, which are related by the M\"{o}bius transformations that permute their singularities. In our case, though, exchanging $1$ with $-1$ keeps equation~\eqref{diff-eq:spatial-mod-bessel-AL} the same.
\subsubsection{Abstract argument for Borel regularity}\label{bessel-regularity-AL}
The analysis in Sections~\ref{big-idea}--\ref{neg-root-AL} picks out a frame in the space of analytic solutions of~\eqref{eqn:reg-mod-bessel-AL}. The frame is generated by solutions of the form $\laplace_{\zeta, 1} v_1$ and $\laplace_{\zeta, -1} v_{-1}$, with $v_\alpha \in \singexpalg{-1/2}(\Omega_\alpha)$.

The Poincar\'e algorithm, as discussed in Section~\ref{sec:asympt-AL} picks out a frame in a space of formal solutions of equation~\eqref{eqn:reg-mod-bessel-AL}. The frame is generated by $1$-Gevrey trans-monomial solutions of the form ${\rm e}^{-z} z^{-1/2} \tilde{W}_1$ and $e^z z^{-1/2} \tilde{W}_{-1}$, with $\tilde{W}_j \in \C\big\llbracket z^{-1} \big\rrbracket_1$. Reprising the proof of Theorem~\ref{thm:summability_ODE}, we will show that these solutions are Borel summable, and their Borel sums generate the same frame as $\laplace_{\zeta, 1} v_1$ and $\laplace_{\zeta, -1} v_{-1}$.

The Borel transform $\borel_\zeta$ maps \smash{${\rm e}^{-\alpha z} z^{-1/2} \C\big\llbracket z^{-1} \big\rrbracket_1$} into \smash{$\zeta_\alpha^{-1/2} \C\{\zeta_\alpha\}$}, as discussed in Sections~\ref{sec:action_transseries} and \ref{sec:new-summability-proof}. In particular, it sends each Poincar\'e solution \smash{${\rm e}^{-\alpha z} z^{-1/2} \series{W}_\alpha$} to a convergent series \smash{$\series{v}_\alpha \in \zeta_\alpha^{-1/2} \C\{\zeta_\alpha\}$}. Make a bounded sector $\Omega'_\alpha$ by intersecting $\Omega_\alpha$ with a disk around $\zeta_\alpha = 0$ whose radius is smaller than the radius of convergence of $\series{v}_\alpha$. On $\Omega'_\alpha$, the sum $\hat{v}_\alpha$ of $\series{v}_\alpha$ can be bounded by a shifted geometric series, showing that \smash{$\zeta_\alpha^{1/2} \hat{v}_\alpha$} is bounded by a constant. Therefore, $\hat{v}_\alpha$ belongs to $\singexpalg{-1/2}(\Omega'_\alpha)$.

The Borel transform
\[ \borel_\zeta \maps\ {\rm e}^{-\alpha z} z^{-1/2} \C\big\llbracket z^{-1} \big\rrbracket_1 \to \zeta_\alpha^{-1/2} \C\{\zeta_\alpha\} \]
also turns formal trans-monomial solutions of equation~\eqref{eqn:reg-mod-bessel-AL} into formal power series solutions of equation~\eqref{int-eq:spatial-mod-bessel-AL}. Thus, by the dominated convergence theorem, $\hat{v}_\alpha$ is an analytic solution of equation~\eqref{int-eq:spatial-mod-bessel-AL} on $\Omega'_\alpha$. The uniqueness part of \cite[Theorem~4.6]{reg-sing-volterra} then shows that $\hat{v}_\alpha$ must be a scalar multiple of the solution $v_\alpha$ we found in Sections \ref{big-idea}--\ref{neg-root-AL}. Since $v_\alpha$ belongs to~\smash{$\singexpalg{-1/2}(\Omega_\alpha)$} for the full, unbounded sector $\Omega_\alpha$, it has a well-defined Laplace transform along any ray \smash{$\mathcal{J}^\theta_{\zeta, \alpha}$} in $\Omega_\alpha$. Therefore, the Poincar\'e solution \smash{${\rm e}^{-\alpha z} z^{-1/2} \series{W}_\alpha$} has a well-defined Borel sum along \smash{$\mathcal{J}^\theta_{\zeta, \alpha}$}, which is a scalar multiple of \smash{${\rm e}^{-\alpha z} V_\alpha$}. It follows, by Lemma~\ref{lem:laplace-bridge}, that~${\rm e}^{-\alpha z} V_\alpha$ is Borel regular.
\subsubsection{Confirmation of Borel regularity}\label{confirmation-borel-regularity}
We can verify the conclusions of Section~\ref{bessel-regularity-AL} using our explicit expressions for the formal power series $\tilde{W}_\alpha$ and the functions $v_\alpha$. We found in Section~\ref{sec:asympt-AL} that
\begin{gather*}
\tilde{W}_1 = \sum_{k = 0}^{\infty} \frac{\bigl(\frac{1}{2}-\frac{m}{n}\bigr)_k \bigl(\frac{1}{2}+\frac{m}{n}\bigr)_k}{k!} \left(-\frac{1}{2}\right)^k z^{-k} ,\qquad
\tilde{W}_{-1} = \sum_{k = 0}^{\infty} \frac{\bigl(\frac{1}{2}-\frac{m}{n}\bigr)_k \bigl(\frac{1}{2}+\frac{m}{n}\bigr)_k}{k!} \left(\frac{1}{2}\right)^k z^{-k}.
\end{gather*}
Computing
\begin{align*}
\borel_\zeta \bigl[ {\rm e}^{-z} z^{-1/2} \tilde{W}_1 \bigr] & = \borel_\zeta \left[ {\rm e}^{-z} \sum_{k = 0}^{\infty} \frac{\bigl(\frac{1}{2}-\frac{m}{n}\bigr)_k \bigl(\frac{1}{2}+\frac{m}{n}\bigr)_k}{k!} \left(-\frac{1}{2}\right)^k z^{-k-\frac{1}{2}} \right] \\
& = \sum_{k = 0}^{\infty} \frac{\bigl(\frac{1}{2}-\frac{m}{n}\bigr)_k \bigl(\frac{1}{2}+\frac{m}{n}\bigr)_k}{k!} \left(-\frac{1}{2}\right)^k \frac{\zeta_1^{k-\frac{1}{2}}}{\Gamma\bigl(k+\frac{1}{2}\bigr)} \\
& = \sum_{k = 0}^{\infty} \frac{\bigl(\frac{1}{2}-\frac{m}{n}\bigr)_k \bigl(\frac{1}{2}+\frac{m}{n}\bigr)_k}{k!} \left(-\frac{1}{2}\right)^k \frac{\zeta_1^{k-\frac{1}{2}}}{\Gamma\bigl(\frac{1}{2}\bigr) \bigl(\frac{1}{2}\bigr)_k} \\
& = \frac{\zeta_1^{-\frac{1}{2}}}{\Gamma\bigl(\frac{1}{2}\bigr)} \sum_{k = 0}^{\infty} \frac{\bigl(\frac{1}{2}-\frac{m}{n}\bigr)_k \bigl(\frac{1}{2}+\frac{m}{n}\bigr)_k}{\bigl(\frac{1}{2}\bigr)_k} \left(-\frac{1}{2}\right)^k \frac{\zeta_1^k}{k!},
\end{align*}
we see that $\borel\bigl[ {\rm e}^{-z} z^{-1/2} \tilde{W}_1 \bigr]$ sums to
\[ \frac{1}{\Gamma(1/2)} \zeta_1^{-1/2} {}_2F_1\left(\frac{1}{2}-\frac{m}{n}, \frac{1}{2}+\frac{m}{n}; \frac{1}{2}; -\frac{1}{2}\zeta_1\right). \]
Looking back at Section~\ref{pos-root-AL}, we recognize this as a scalar multiple of $v_1$.

Through a similar calculation, we see that $\borel\bigl[{\rm e}^z z^{-1/2} \tilde{W}_{-1} \bigr]$ sums to
\[ \frac{1}{\Gamma(1/2)} \zeta_{-1}^{-1/2} {}_2F_1\left(\frac{1}{2}-\frac{m}{n}, \frac{1}{2}+\frac{m}{n}; \frac{1}{2}; \frac{1}{2}\zeta_{-1}\right). \]
Looking back at Section~\ref{neg-root-AL}, we recognize this as a scalar multiple of $v_{-1}$.
\subsubsection{Thimble projection reasoning}\label{contour-argument-AL}
We can also study the Airy--Lucas functions by applying the thimble projection formula to integral~\eqref{integral:mod-bessel-rational-AL}, specializing the reasoning behind Lemma~\ref{lem:thimble_proj_formula-proof}. We recast integral~\eqref{integral:mod-bessel-rational-AL} into the~$\zeta$ plane by setting $-\zeta = T_n(u)$, which implies that $-{\rm d}\zeta = n U_{n-1}(u) {\rm d}u$. Recall from below equation~\eqref{integral:AL} that the integration contour $\mathcal{C}^\theta_j$ runs through the critical point $u = \cos\bigl(\frac{j}{n}\pi\bigr)$, with its direction determined by the parameter $\theta$. The critical point splits \smash{$\mathcal{C}^\theta_j$} into two pieces: the {\em incoming} branch, where the orientation of the thimble runs toward the critical point, and the {\em outgoing} branch, where the orientation points away. Recalling that $\mathcal{C}^\theta_j$ is a preimage of $\mathcal{J}^\theta_{\zeta, \mp 1}$, we get
\begin{align*}
K^{(j)}_{m/n}(z) & = \frac{n}{2 \sinh\bigl(\frac{m}{n} {\rm i}\pi\bigr)} \int_{\mathcal{C}^\theta_j} \exp\left[z T_n(u)\right] U_{m-1}(u) \,{\rm d}u \\
& = -\frac{1}{2\sinh\bigl(\frac{m}{n} {\rm i}\pi\bigr)} \left[ \int_{\mathcal{J}^\theta_{\zeta, \mp 1}} {\rm e}^{-z\zeta} \frac{U_{m-1}(u_+)}{U_{n-1}(u_+)} \,{\rm d}\zeta - \int_{\mathcal{J}^\theta_{\zeta, \mp 1}} {\rm e}^{-z\zeta} \frac{U_{m-1}(u_-)}{U_{n-1}(u_-)} \,{\rm d}\zeta \right],
\end{align*}
where $u_-$ and $u_+$ are the lifts to the incoming and outgoing branches of $\mathcal{C}^\theta_j$, respectively.

Since the Chebyshev polynomial $T_n$ is defined by the identity $T_n(\cos(\phi)) = \cos(n\phi)$, we introduce a new variable $\phi$ with $u = \cos(\phi)$ and $-\zeta = \cos(n\phi)$. On the $\phi$ plane, which is an infinite branched cover of the $u$ plane, we can lift $\mathcal{C}^\theta_j$ to the path $\frac{j}{n}\pi + {\rm i}\R$. The incoming and outgoing branches lift to \smash{$\frac{j}{n}\pi - {\rm i}[0, \infty)$} and \smash{$\frac{j}{n}\pi + {\rm i}[0, \infty)$}, respectively.

We can now use \cite[identity~(15.4.16)]{dlmf} to write the integrand explicitly in terms of $\zeta$
\begin{align*}
\frac{U_{m-1}(\cos(\phi))}{U_{n-1}(\cos(\phi))} &= \frac{\sin(m\phi)}{\sin(\phi)}\frac{\sin(\phi)}{\sin(n \phi)}= \frac{\sin(m\phi)}{\sin(n \phi)}\\
& = \frac{m}{n} {}_2F_1\left(\frac{1}{2} - \frac{m}{2n}, \frac{1}{2} + \frac{m}{2n}; \frac{3}{2}; \sin(n \phi)^2\right) \\
& = \frac{m}{n} {}_2F_1\left(\frac{1}{2} - \frac{m}{2n}, \frac{1}{2} + \frac{m}{2n}; \frac{3}{2}; 1 - \zeta^2\right).
\end{align*}
Initially, we only know that this equation holds in the disk $|\phi| < \frac{\pi}{2n}$. However, observing that the left-hand side is meromorphic throughout the $\phi$ plane, we can deduce by analytic continuation that the equation holds throughout the $\phi$ plane.

Next, we simplify the integrand using identities (15.8.4) and (15.8.27)--(15.8.28) from \cite{dlmf}, which tell us that
\begin{gather*}
{}_2F_1\left(\frac{1}{2} - \frac{m}{2n}, \frac{1}{2} + \frac{m}{2n}; \frac{3}{2}; 1 - \zeta^2\right) \\
\qquad = \hphantom{+} \frac{\pi}{\Gamma\bigl(1-\frac{m}{2n}\bigr) \Gamma\bigl(1 + \frac{m}{2n}\bigr)} {}_2F_1\left(\frac{1}{2} - \frac{m}{2n}, \frac{1}{2} + \frac{m}{2n}; \frac{1}{2}; \zeta^2\right) \\
\qquad\hphantom{=} - \frac{\pi \zeta}{\Gamma\bigl(\frac{1}{2} - \frac{m}{2n}\bigr) \Gamma\bigl(\frac{1}{2} + \frac{m}{2n}\bigr)} {}_2F_1\left(1 - \frac{m}{2n}, 1 + \frac{m}{2n}; \frac{3}{2}; \zeta^2\right) \\
\qquad = {}_2F_1\left(1 - \frac{m}{n}, 1 + \frac{m}{n}; \frac{3}{2}; \frac{1}{2} - \frac{\zeta}{2}\right) + {}_2F_1\left(1-\frac{m}{n}, 1 + \frac{m}{n}; \frac{3}{2}; \frac{1}{2} + \frac{\zeta}{2}\right) \\
\qquad \hphantom{=} + \frac{1}{2} {}_2F_1\left(1-\frac{m}{n}, 1+\frac{m}{n}; \frac{3}{2}; \frac{1}{2} - \frac{\zeta}{2}\right) - \frac{1}{2} {}_2F_1\left(1-\frac{m}{n}, 1 + \frac{m}{n}; \frac{3}{2}; \frac{1}{2} + \frac{\zeta}{2}\right) \\
\qquad = \frac{3}{2} {}_2F_1\left(1 - \frac{m}{n}, 1 + \frac{m}{n}; \frac{3}{2};\frac{1}{2} - \frac{\zeta}{2}\right) + \frac{1}{2} {}_2F_1\left(1 - \frac{m}{n}, 1 + \frac{m}{n}; \frac{3}{2}; \frac{1}{2} + \frac{\zeta}{2}\right)
\end{gather*}
away from the line $\Re(\zeta) = 0$ and the rays $\mathcal{J}^0_{\zeta,1}$ and $\mathcal{J}^\pi_{\zeta,-1}$.

Returning to the projected thimble integral
\[ K^{(j)}_{m/n}(z) = -\frac{1}{2\sinh\bigl(\frac{m}{n} {\rm i}\pi\bigr)} \int_{\mathcal{J}^\theta_{\zeta, \mp 1}} {\rm e}^{-z\zeta} \left[ \frac{U_{m-1}(u_+)}{U_{n-1}(u_+)} - \frac{U_{m-1}(u_-)}{U_{n-1}(u_-)} \right] {\rm d}\zeta, \]
we can see the integrand as the variation of the function
\[ \frac{U_{m-1}(\cos(\phi))}{U_{n-1}(\cos(\phi))} = \frac{3m}{2n} {}_2F_1\left(1 - \frac{m}{n}, 1 + \frac{m}{n}; \frac{3}{2};\frac{1}{2} - \frac{\zeta}{2}\right) + \frac{m}{2n} {}_2F_1\left(1 - \frac{m}{n}, 1 + \frac{m}{n}; \frac{3}{2}; \frac{1}{2} + \frac{\zeta}{2}\right) \]
across the branch cut $\mathcal{J}^\theta_{\zeta, \pm 1}$. When $j$ is odd, only the second term will contribute to the jump, because the first term is regular at $\zeta = 1$. Similarly, when $j$ is even, only the first term will contribute to the jump. We can write the jump explicitly using identity~(15.2.3) from \cite{dlmf}. For odd $j$, we have
\[ K^{(j)}_{m/n}(z) = -\frac{1}{2} \int_{\mathcal{J}^\theta_{\zeta, 1}} {\rm e}^{-z\zeta} \left(-\frac{1}{2}+\frac{\zeta}{2}\right)^{-1/2} {}_2F_1\left(\frac{1}{2} - \frac{m}{n}, \frac{1}{2} + \frac{m}{n}; \frac{1}{2}; \frac{1}{2} - \frac{\zeta}{2}\right) {\rm d}\zeta, \]
and for even $j$, we have
\[ K^{(j)}_{m/n}(z) = \frac{3}{2} \int_{\mathcal{J}^\theta_{\zeta, -1}} {\rm e}^{-z\zeta} \left(-\frac{1}{2}-\frac{\zeta}{2}\right)^{-1/2} {}_2F_1\left(\frac{1}{2} - \frac{m}{n}, \frac{1}{2} + \frac{m}{n}; \frac{1}{2}; \frac{1}{2} + \frac{\zeta}{2}\right) {\rm d}\zeta. \]

Comparing the expressions above with the expressions for $v_{\pm 1}$ computed in Sections~\ref{pos-root-AL}--\ref{neg-root-AL}, we notice that
\[ K^{(j)}_{m/n}(z) = \begin{cases}
\hphantom{-}\dfrac{\rm i}{2}\laplace_{\zeta, 1} v_1 & j \text{ odd}, \vspace{1mm}\\
-\dfrac{3{\rm i}}{2}\laplace_{\zeta, -1} v_{-1} & j \text{ even}.
\end{cases} \]

\subsubsection{A flavor of resurgence: the Stokes phenomenon in the position domain}\label{resurgence-AL}

So far, we have treated the Airy--Lucas functions as separate elements of a frame of solutions. However, \'{E}calle's theory of {\em resurgence} reveals that these solutions are deeply intertwined. In fact, you can recover the whole frame from any one of its elements by studying the analytic continuation of the corresponding function on the position domain.

As an example, consider the Borel regular solution ${\rm e}^{-z} V_1$ of the modified Bessel equation~\eqref{eqn:mod-bessel-AL}, which we found by solving an integral equation in Section~\ref{big-idea} and by evaluating a thimble integral in Section~\ref{contour-argument-AL}. This solution arises as the Laplace transform $\laplace^\theta_{\zeta, 1} v_1$ of the function
\[ v_1 = \left(-\frac{\zeta_1}{2}\right)^{-1/2} {}_2F_1\left(\frac{1}{2}-\frac{m}{n},\frac{1}{2}+\frac{m}{n};\frac{1}{2};-\frac{\zeta_1}{2}\right), \]
on the position domain, with $\zeta=1+\zeta_1$ as before. The function $v_1$ is singular at $\zeta = 1$, but the hypergeometric factor ${}_2F_1\bigl(\frac{1}{2}-\frac{m}{n},\frac{1}{2}+\frac{m}{n};\frac{1}{2};-\frac{\zeta_1}{2}\bigr)$ is analytic at that point, so $v_1$ has the typical form of a {\em singularity} in the formalism of resurgent functions~\cite[Section~2]{sauzin2021variations}. We can see~${{}_2F_1\bigl(\frac{1}{2}-\frac{m}{n},\frac{1}{2}+\frac{m}{n};\frac{1}{2};-\frac{\zeta_1}{2}\bigr)}$ as a new holomorphic function, analytically continued from a~germ at $\zeta = 1$. This new function has a branch cut singularity at $\zeta=-1$, and we can find its jump across the branch cut $\zeta \in (-\infty, -1]$ using \cite[equation~(15.2.3)]{dlmf}. Along the branch cut, we have
\begin{gather*}
{}_2F_1\left(\frac{1}{2}-\frac{m}{n},\frac{1}{2}+\frac{m}{n};\frac{1}{2};-\frac{\zeta_1}{2}+{\rm i}\varepsilon\right)-{}_2F_1\left(\frac{1}{2}-\frac{m}{n},\frac{1}{2}+\frac{m}{n};\frac{1}{2};-\frac{\zeta_1}{2}-{\rm i}\varepsilon\right) \\
\qquad=\frac{2 \pi {\rm i}}{\Gamma\bigl(\frac{1}{2}-\frac{m}{n}\bigr)\Gamma\bigl(\frac{1}{2}+\frac{m}{n}\bigr)} \left(-\frac{\zeta_{-1}}{2}\right)^{-1/2} {}_2F_1\left(\frac{m}{n},-\frac{m}{n};\frac{1}{2};\frac{\zeta_{-1}}{2}\right) \\
\qquad=2{\rm i}\cos\left(\frac{m}{n}\pi\right) \left(-\frac{\zeta_{-1}}{2}\right)^{-1/2} {}_2F_1\left(\frac{m}{n},-\frac{m}{n};\frac{1}{2};\frac{\zeta_{-1}}{2}\right) \\
\qquad=2\cos\left(\frac{m}{n}\pi\right) \left(\frac{\zeta_{-1}}{2}\right)^{-1/2} {}_2F_1\left(\frac{m}{n},-\frac{m}{n};\frac{1}{2};\frac{\zeta_{-1}}{2}\right) \\
\qquad=2\cos\left(\frac{m}{n}\pi\right) \left(\frac{\zeta_{-1}}{2}\right)^{-1/2} \left(-\frac{\zeta_{1}}{2}\right)^{1/2} {}_2F_1\left(\frac{1}{2}-\frac{m}{n},\frac{1}{2}+\frac{m}{n};\frac{1}{2};\frac{\zeta_{-1}}{2}\right) \\
\qquad=\left(-\frac{\zeta_{1}}{2}\right)^{1/2} 2\cos\left(\frac{m}{n}\pi\right) v_{-1}.
\end{gather*}
Multiplying both sides by \smash{$(-\zeta_1/2)^{-1/2}$} reveals a formula for the jump of $v_1$ across the branch cut. For any point $p$ on the branch cut,
\begin{equation}\label{eqn:resurg-rel-1}
v_1(p_-) - v_1(p_+) = 2\cos\bigl(\frac{m}{n}\pi\bigr) v_{-1}(p),
\end{equation}
where $v_1(p_-)$ and $v_1(p_+)$ respectively denote the evaluation of $v_1(p)$ by analytic continuation from the regions $\Im(\zeta) < 0$ and $\Im(\zeta) > 0$. The appearance of the function $v_{-1}$ corresponding to the other Borel regular solution ${\rm e}^z V_{-1}$ is an example of the {\em resurgence} phenomenon~\cite{EcalleI,EcalleII,EcalleIII}.

Reviews like \cite{aniceto2019primer,Dorigoni,diverg-resurg-i} discuss many details and applications of resurgence. Berry and Howls also discuss the optimal truncation of thimble integrals from a resurgence point of view~\cite{Berry_Howls}. We would like to emphasize one application: through resurgence, we can use calculations in the position domain to understand the Stokes phenomenon in the frequency domain. This lets us skip the last and often trickiest step of Borel summation: using the Laplace transform to map the functions under consideration back into the frequency domain. The {\em formalism of singularities} and the {\em alien calculus} are two of the tools that make this possible.

To compute the Stokes constants for the Airy--Lucas functions, we first repeat the calculation above for
\[ v_{-1} = \left(\frac{\zeta_{-1}}{2}\right)^{-1/2} {}_2F_1\left(\frac{1}{2}-\frac{m}{n},\frac{1}{2}+\frac{m}{n};\frac{1}{2};\frac{\zeta_{-1}}{2}\right). \]
We find that
\begin{gather*}
{}_2F_1\left(\frac{1}{2}-\frac{m}{n},\frac{1}{2}+\frac{m}{n};\frac{1}{2};\frac{\zeta_{-1}}{2}+i\varepsilon\right)-{}_2F_1\left(\frac{1}{2}-\frac{m}{n},\frac{1}{2}+\frac{m}{n};\frac{1}{2};\frac{\zeta_{-1}}{2}-i\varepsilon\right) \\
\qquad=\frac{2 \pi {\rm i}}{\Gamma\bigl(\frac{1}{2}-\frac{m}{n}\bigr)\Gamma\bigl(\frac{1}{2}+\frac{m}{n}\bigr)} \left(\frac{\zeta_1}{2}\right)^{-1/2} {}_2F_1\left(\frac{m}{n},-\frac{m}{n};\frac{1}{2};-\frac{\zeta_1}{2}\right) \\
\qquad=2{\rm i}\cos\left(\frac{m}{n}\pi\right) \left(\frac{\zeta_1}{2}\right)^{-1/2} {}_2F_1\left(\frac{m}{n},-\frac{m}{n};\frac{1}{2};-\frac{\zeta_1}{2}\right) \\
\qquad=-2\cos\left(\frac{m}{n}\pi\right) \left(-\frac{\zeta_1}{2}\right)^{-1/2} {}_2F_1\left(\frac{m}{n},-\frac{m}{n};\frac{1}{2};-\frac{\zeta_1}{2}\right) \\
\qquad=-2\cos\left(\frac{m}{n}\pi\right) \left(-\frac{\zeta_{1}}{2}\right)^{-1/2} \left(\frac{\zeta_{-1}}{2}\right)^{1/2} {}_2F_1\left(\frac{1}{2}-\frac{m}{n},\frac{1}{2}+\frac{m}{n};\frac{1}{2};-\frac{\zeta_1}{2}\right) \\
\qquad=\left(\frac{\zeta_{-1}}{2}\right)^{1/2} (-2)\cos\left(\frac{m}{n}\pi\right) v_1.
\end{gather*}
Multiplying both sides by \smash{$(\zeta_{-1}/2)^{-1/2}$} gives the jump formula
\begin{equation}\label{eqn:resurg-rel-2}
v_{-1}(p_+) - v_{-1}(p_-) = -2\cos\left(\frac{m}{n}\pi\right) v_1(p)
\end{equation}
for any point $p$ on the branch cut $\zeta \in [1, \infty)$, where again $v_{-1}(p_+)$ and $v_{-1}(p_-)$ denote evaluation by analytic continuation from $\Im(\zeta) > 0$ and $\Im(\zeta) < 0$.

Equations~\eqref{eqn:resurg-rel-1} and~\eqref{eqn:resurg-rel-2} are a manifestation of the Stokes phenomenon in the position domain. When we vary the direction of integration for the Laplace transform, we expect discontinuities at $\laplace^\pi_{\zeta,1}$ and $\laplace^0_{\zeta,-1}$, because each of the contours $\mathcal{J}^\pi_{\zeta,1}$ and $\mathcal{J}^0_{\zeta,-1}$ hits a singularity other than its starting point. The Stokes constant for each critical direction relates the results of the Laplace transforms on either side. In our example,
\begin{gather*}
\laplace_{\zeta,1}^{\pi+\varepsilon} v_1 - \hphantom{{}_{-}}\laplace_{\zeta,1}^{\pi-\varepsilon} v_1 = \hphantom{-}2\cos\left(\frac{m}{n}\pi\right) {\rm e}^{-z} \laplace_{\zeta,-1}^{\pi} v_{-1}, \\
\laplace_{\zeta,-1}^{\varepsilon} v_{-1} - \laplace_{\zeta,-1}^{-\varepsilon} v_{-1} = -2\cos\left(\frac{m}{n}\pi\right) {\rm e}^z \laplace_{\zeta,1} v_{1},
\end{gather*}
so the Stokes constants for the directions $\pi$ and $0$ are $2\cos\bigl(\frac{m}{n}\pi\bigr)$ and $-2\cos\bigl(\frac{m}{n}\pi\bigr)$, respectively.

Because the Airy--Lucas functions can always be expressed as thimble integrals, it may seem surprising that the Stokes constants are not always integers. For thimble integrals given by Morse functions, the Stokes constants are the intersection numbers of \textit{dual pairs of thimbles}, according to Picard--Lefschetz formula (see \cite[Section~5]{pham} and \cite[Chapter~1]{Arnold}). However, the functions~$T_n(u)$ defining the Airy--Lucas thimble integrals~\eqref{integral:mod-bessel-rational-AL} are only Morse in the Airy case $n = 3$ and the degenerate limit $n \to \infty$ that we will discuss in Section~\ref{sec:bessel-0}. In the other cases, we speculate that the Stokes constants are non-integer because there are multiple thimbles over each critical value.
\subsection{Modified Bessel}\label{sec:mod-bessel-lift}
Our analysis of the Airy--Lucas functions boiled down to an analysis of equation~\eqref{eqn:mod-bessel-AL}: the modified Bessel equation with a rational parameter $\mu = \frac{m}{n}$. We now do a general analysis of the modified Bessel equation, allowing the parameter $\mu$ to be any complex number
\begin{equation}\label{eqn:mod-bessel}
\left[z^2 \left(\frac{\partial}{\partial z}\right)^2 + z \frac{\partial}{\partial z} - \bigl[\mu^2 + z^2\bigr]\right] \Phi = 0.
\end{equation}
The ODE reasoning of Sections \ref{big-idea}--\ref{confirmation-borel-regularity} works basically the same in this more general setting, so we will only briefly state the analogous results.

\textbf{Asymptotic analysis.}
Equation~\eqref{eqn:mod-bessel} has a basis of formal solutions ${\rm e}^{-z} z^{-1/2} \series{W}_1$ and $ {\rm e}^z z^{-1/2} \series{W}_{-1}$,
given by
\begin{gather*}
\series{W}_{\mu,1} = 1 - \frac{\bigl(\frac{1}{2}-\mu\bigr)\bigl(\frac{1}{2}+\mu\bigr)}{2 \cdot 1!} z^{-1} + \frac{\bigl(\frac{1}{2}-\mu\bigr)_2\bigl(\frac{1}{2}+\mu\bigr)_2}{2^2 \cdot 2!} z^{-2}
 - \frac{\bigl(\frac{1}{2}-\mu\bigr)_3\bigl(\frac{1}{2}+\mu\bigr)_3}{2^3 \cdot 3!} z^{-3} + \cdots, \\
\series{W}_{\mu,-1} = 1 + \frac{\bigl(\frac{1}{2}-\mu\bigr)\bigl(\frac{1}{2}+\mu\bigr)}{2 \cdot 1!} z^{-1} + \frac{\bigl(\frac{1}{2}-\mu\bigr)_2\bigl(\frac{1}{2}+\mu\bigr)_2}{2^2 \cdot 2!} z^{-2}
+ \frac{\bigl(\frac{1}{2}-\mu\bigr)_3\bigl(\frac{1}{2}+\mu\bigr)_3}{2^3 \cdot 3!} z^{-3} + \cdots.
\end{gather*}

\textbf{Frame of analytic solutions.}
Equation~\eqref{eqn:mod-bessel} also has a basis of analytic solutions~$\laplace^\theta_{\zeta, 1} v_1$ and~$\laplace^\theta_{\zeta, -1} v_{-1}$, which are the Laplace transforms of the functions
\begin{gather*}
v_1 = -{\rm i}\sqrt{2} \zeta_{1}^{-1/2} {}_2F_1\left(\frac{1}{2}-\mu, \frac{1}{2}+\mu; \frac{1}{2}; -\frac{1}{2}\zeta_{1}\right),\\
v_{-1} =\sqrt{2} \zeta_{-1}^{-1/2} {}_2F_1\left(\frac{1}{2}-\mu, \frac{1}{2}+\mu; \frac{1}{2}; \frac{1}{2}\zeta_{-1}\right)
\end{gather*}
along the directions $\theta \neq \pi$ and $\theta \neq 0$, respectively.

\textbf{Borel regularity.}
The Borel transforms $\mathcal{B}_\zeta\bigl[z^{-1/2} {\rm e}^{-\alpha z} \series{W}_{\pm 1}\bigr]$ sum to scalar multiples of~$v_{\pm 1}$, implying through Lemma~\ref{lem:laplace-bridge} that the solutions $\laplace^\theta_{\zeta, \pm 1} v_{\pm 1}$ are Borel regular.

On the other hand, the thimble projection reasoning of Section~\ref{contour-argument-AL} has to be generalized, as we will discuss in Section~\ref{countable-cover} below.

\subsubsection{Lifting to a countable cover}\label{countable-cover}
\begin{figure}\centering
\includegraphics[width=10cm]{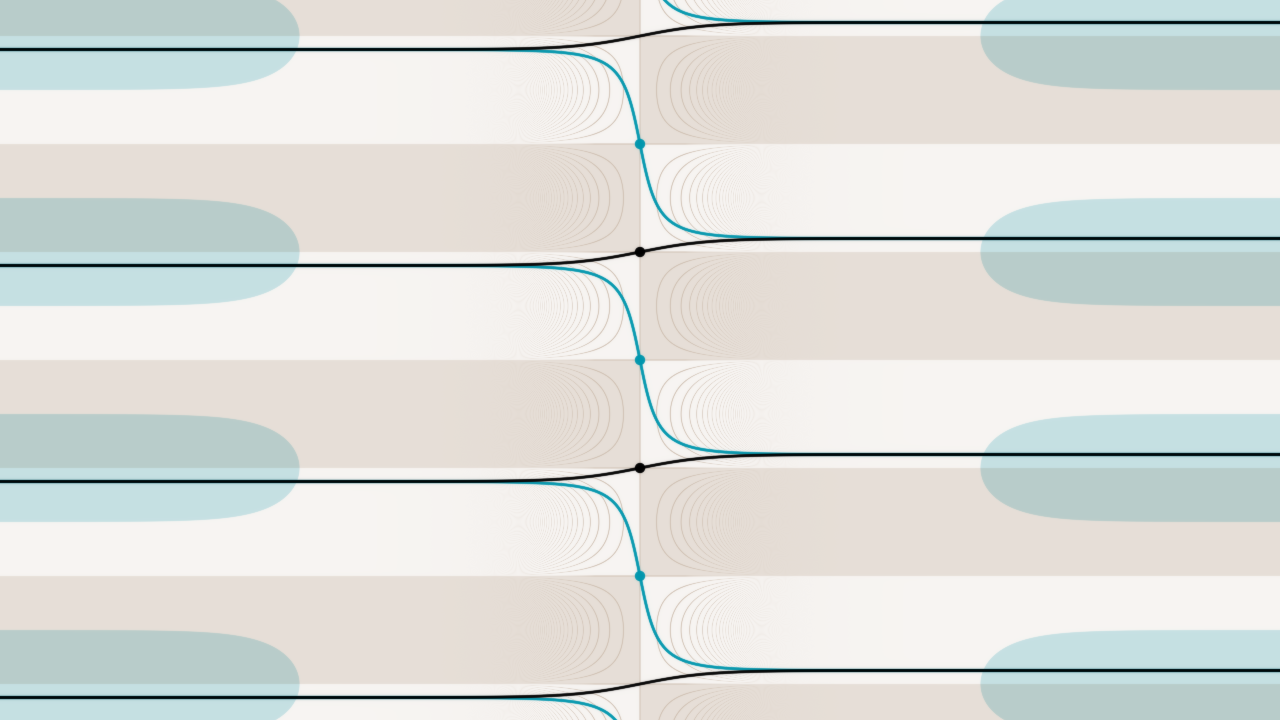}
\caption{The paths $\mathcal{C}^{\pi/8}_j$ in the $t$ plane.}\label{fig:bessel_unrolled}
\end{figure}

Formula~\eqref{integral:mod-bessel-rational-AL} expresses the modified Bessel function $K_{m/n}$ as an exponential integral on a finite cover of $\C$. Lifting to a countable cover reveals this formula as a special case of a general integral formula for modified Bessel functions.
Setting $u = \cosh(t/n)$, and recalling that
$
\cosh(n\tau) = T_n(\cosh(\tau)) $ and $
\sinh(m\tau) = U_{m-1}(\cosh(\tau))  \sinh(\tau)$,
we can rewrite formula~\eqref{integral:mod-bessel-rational-AL} as
\begin{align}
 K^{(j)}_{m/n}(z) & = \frac{n}{2 \sinh\bigl(\frac{m}{n} {\rm i}\pi\bigr)} \int_\gamma \exp[z T_n(u)] U_{m-1}(u) \,{\rm d}u \nonumber\\
 & = \frac{n}{2 \sinh\bigl(\frac{m}{n} {\rm i}\pi\bigr)} \int_{\mathcal{C}^\theta_j} \exp[z \cosh(t)] U_{m-1}(\cosh(t/n)) \sinh(t/n)\, {\rm d}(t/n)\nonumber \\
\label{integral:mod-bessel-lifted} & = \frac{1}{2 \sinh\bigl(\frac{m}{n} {\rm i}\pi\bigr)} \int_{\mathcal{C}^\theta_j} \exp[z \cosh(t)] \sinh\left(\frac{m}{n} t\right) {\rm d}t,
\end{align}
using the path $\mathcal{C}^\theta_j$ described below.
\begin{center}
\begin{tabular}{l|l|l}
When $j$ is \dots & $\mathcal{C}^\theta_j$ comes from $-\infty$ along \dots & and goes to $+\infty$ along \dots \\[1mm] \hline
\vphantom{\rule{0mm}{5mm}} Odd & $(-\infty,0) + (j\pi-\theta){\rm i}$ & $(0, \infty) + (j\pi+\theta){\rm i}$ \\[1mm] \hline
\vphantom{\rule{0mm}{5mm}} Even & $(-\infty,0) + \bigl((j+1)\pi-\theta\bigr){\rm i}$ & $(0, \infty) + \bigl((j-1)\pi+\theta\bigr){\rm i}$ \\[1mm]
\end{tabular}
\end{center}
Because $\cosh(t)$ is periodic, $j$ only affects the integral through its parity.
The path $\mathcal{C}^\theta_j$ is shown, for various $j$, in Figure~\ref{fig:bessel_unrolled}.

For any $\mu \in \C \smallsetminus \Z$, the classical modified Bessel function $K_\mu(z)$ can be expressed as the thimble integral
\begin{equation}\label{int:mod-bessel-gen}
 K_\mu(z) = \frac{1}{2 \sinh(\mu {\rm i}\pi)} \int_{\mathcal{C}^\theta_1} \exp[z \cosh(t)] \sinh(\mu t) \,{\rm d}t
\end{equation}
for $|\theta| < \frac{\pi}{2}$. This formula can be deduced, with a lot of finagling, from \cite[formulas~(10.32.12) and~(10.27.4)]{dlmf}. The integral converges when $z$ is in the right half-plane. Choosing a rational parameter $\mu = \frac{m}{n}$ gives formula~\eqref{integral:mod-bessel-lifted}, showing that \smash{$K^{(j)}_{m/n}(z)$} is the classical modified Bessel function $K_{m/n}(z)$ when $j$ is odd.

We can now apply the thimble projection formula, using the same reasoning as in Section~\ref{contour-argument-AL}. We first recast the integral into the $\zeta$ plane by setting $-\zeta = \cosh(t)$. This gives us the expression
\[ K^{(j)}_\mu(z) = -\frac{1}{2\sinh(\mu {\rm i}\pi)} \int_{\mathcal{J}^\theta_{\zeta, \mp 1}} {\rm e}^{-z\zeta} \left[ \frac{\sinh(\mu t_+)}{\sinh(t_+)} - \frac{\sinh(\mu t_-)}{\sinh(t_-)} \right] {\rm d}\zeta, \]
where $t_-$ and $t_+$ are the lifts to the incoming and outgoing branches of $\mathcal{C}^\theta_j$. We then use identity~(15.4.16) from \cite{dlmf} to write the integral explicitly in terms of $\zeta$
\begin{align*}
\frac{\sinh(\mu t)}{\sinh(t)} & = \mu {}_2F_1\left(\frac{1}{2} - \frac{\mu}{2}, \frac{1}{2} + \frac{\mu}{2}; \frac{3}{2}; -\sinh(t)^2\right) = \mu {}_2F_1\left(\frac{1}{2} - \frac{\mu}{2}, \frac{1}{2} + \frac{\mu}{2}; \frac{3}{2}; 1 - \zeta^2\right).
\end{align*}
Like before, we can deduce by analytic continuation that this identity holds throughout the $t$ plane.

Using \cite[formula~(15.8.4)]{dlmf}, followed by formulas (15.8.27) and (15.8.28) from the same source, we can repeat the arguments of Section~\ref{contour-argument-AL}, eventually rewriting the integrand as the variation of the function
\[ \frac{\sinh(\mu t)}{\sinh(t)} = \frac{3\mu}{2} {}_2F_1\left(1 - \mu, 1 + \mu; \frac{3}{2};\frac{1}{2} - \frac{\zeta}{2}\right) + \frac{\mu}{2} {}_2F_1\left(1 - \mu, 1 + \mu; \frac{3}{2}; \frac{1}{2} + \frac{\zeta}{2}\right) \]
across the branch cut $\mathcal{J}^\theta_{\zeta, \mp 1}$. As before, the term that contributes to the jump is the one which is singular at the critical value where the branch cut starts. We use identity~(15.2.3) from \cite{dlmf} to write the jump explicitly. For odd $j$, we get
\begin{equation}\label{eq: K_j_mu at 1}
 K^{(j)}_\mu(z) = -\frac{1}{2} \int_{\mathcal{J}^\theta_{\zeta, 1}} {\rm e}^{-z\zeta} \left(-\frac{1}{2}+\frac{\zeta}{2}\right)^{-1/2} {}_2F_1\left(\frac{1}{2} - \mu, \frac{1}{2} + \mu; \frac{1}{2}; \frac{1}{2} - \frac{\zeta}{2}\right) {\rm d}\zeta,
\end{equation}
and for even $j$, we get
\[
 K^{(j)}_\mu(z) = \frac{3}{2} \int_{\mathcal{J}^\theta_{\zeta, -1}} {\rm e}^{-z\zeta} \left(-\frac{1}{2}-\frac{\zeta}{2}\right)^{-1/2} {}_2F_1\left(\frac{1}{2} - \mu, \frac{1}{2} + \mu; \frac{1}{2}; \frac{1}{2} + \frac{\zeta}{2}\right) {\rm d}\zeta.
\]
\subsubsection[The modified Bessel function with parameter 0]{The modified Bessel function with parameter 0} \label{sec:bessel-0}

When $\mu$ goes to $0$, formula~\eqref{int:mod-bessel-gen} becomes
\[ K_0^{(j)}(z) = \frac{1}{2 \pi {\rm i}} \int_{\mathcal{C}_j^\theta} \exp[z \cosh(t)] t \,{\rm d}t. \]
Let us compute \smash{$K^{(1)}_{0}(z)$} using the contour $\mathcal{C}_1^0$, which runs rightward along the line $\Im(t) = \pi$. In the translated coordinate $w$ defined by $t = w + {\rm i}\pi$, the integral becomes
\begin{align*}
K_0^{(1)}(z) & = \frac{1}{2 \pi {\rm i}} \int_{-\infty}^\infty \exp\left[z \cosh(w + {\rm i}\pi)\right] (w + {\rm i}\pi)\,{\rm d}w \\
& = \frac{1}{2 \pi {\rm i}} \int_{-\infty}^\infty \exp[-z \cosh(w)] w {\rm d}w + \frac{1}{2 \pi {\rm i}} \int_{-\infty}^\infty \exp[-z \cosh(w)] {\rm i}\pi\, {\rm d}w.
\end{align*}
The first integrand is odd with respect to $w = 0$, so it vanishes, leaving
\begin{align}
K_0^{(1)}(z) = \frac{1}{2} \int_{-\infty}^\infty \exp [-z \cosh(w) ] \,{\rm d}w = \int_0^\infty \exp[-z \cosh(w)] \,{\rm d}w. \label{int:k0_cosh}
\end{align}
This is a special case of \cite[formula~(10.32.9)]{dlmf}. Rolling the $t$ plane up into a cylinder, parameterized by the coordinate $s = {\rm e}^t$, we can express \smash{$K_0^{(1)}(z)$} as an exponential period:
\[ K_0^{(1)}(z) = \int_{\mathcal{J}^0_{s, 1}} \exp\left[-z \frac{1}{2}\left(s + \frac{1}{s}\right)\right] \frac{{\rm d}s}{s}. \]

To confirm that the calculation above is consistent with Section~\ref{countable-cover}, recall formula~\eqref{eq: K_j_mu at 1}:
\[ K_0^{(1)}(z) = -\frac{1}{2} \int_{\mathcal{J}^\theta_{\zeta, 1}} {\rm e}^{-z\zeta} \left(-\frac{1}{2}+\frac{\zeta}{2}\right)^{-1/2} {}_2F_1\left(\frac{1}{2}, \frac{1}{2}; \frac{1}{2}; \frac{1}{2} - \frac{\zeta}{2}\right) {\rm d}\zeta. \]
The hypergeometric function in the integrand can be expressed algebraically using identi\-ty~(15.4.13) from \cite{dlmf}. This leads to an expression for the integrand that makes its \smash{$\zeta_{\pm 1}^{1/2}$} singularities even more apparent
\[
K_0^{(1)}(z) = \int_{\mathcal{J}^\theta_{\zeta, 1}} {\rm e}^{-z\zeta} \bigl(\zeta^2 - 1\bigr)^{-1/2}\, {\rm d}\zeta.
\]
We could get the same result from formula~\eqref{int:k0_cosh} by trigonometric substitution. This is a special case of \cite[formula~(10.32.8)]{dlmf}.
\subsection{Generalized Airy}\label{sec:gen-airy}
\begin{figure}[!ht]
\centering
\newcommand{\apathcolor}{ietcoast}
\newcommand{\bpathcolor}{black}
\begin{tikzpicture}
\renewcommand{\dotsize}{0.08}
\fill[pwbeige!10] circle(3.95);
\begin{scope}[\bpathcolor, very thick, -stealth]
 \draw (80:0.8) arc (80:400:0.8);
 \node at (60:0.8) {$\mathscr{b}_0$};
 \foreach \ang/\name in {0/1, 120/2, 240/3} {
 \draw[very thick,-stealth] (0, 0)--++(\ang:4) node[anchor=\ang-180] {$\mathscr{b}_\name$};
 };
 \fill (1.25, 0) circle (\dotsize) node[anchor=north west] {$\frac{1}{2}$};
 \draw[fill=white, thin] circle (\dotsize);
\end{scope}
\begin{scope}[\apathcolor, very thick, -stealth]
 \foreach \ang/\name in {0/1, 120/2, 240/3} {
 \draw[rotate=\ang] (-0.5, 0) +(-120:3.75) .. controls +(62:1) and +(0, -1.2) .. +(-0.75, 0) .. controls +(0, 1.2) and +(-62:1) .. +(120:3.75) node[midway, anchor=\ang+30] {$\mathscr{a}_\name$};
 };
 \fill (-1.25, 0) circle (\dotsize) node[anchor=north east] {$-\frac{1}{2}$};
\end{scope}
\end{tikzpicture}
\caption{The integration contours in the $t$ plane that define $\mathrm{A}_k$, $\mathrm{B}_0$, $\mathrm{B}_k$.}\label{fig:path-generalized-Airy}
\end{figure}
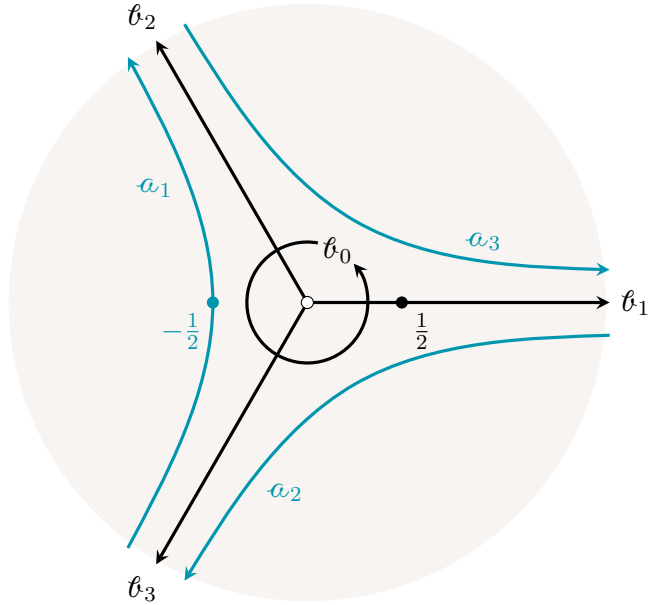

In \cite{Reid} and the appendix of \cite{drazin-reid}, Drazin and Reid construct approximate solutions of the Orr--Sommerfield fluid equation using the generalized Airy functions
\begin{gather*}
\mathrm{A}_k(y,p) = \frac{1}{2 \pi {\rm i}}\int_{\mathscr{a}_k} \exp\left[yt-\frac{t^3}{3}\right] \frac{{\rm d}t}{t^p} ,\qquad p \in \C ,\quad k \in \{1,2,3\}, \\
\mathrm{B}_0(y,p) = \frac{1}{2 \pi {\rm i}}\int_{\mathscr{b}_0} \exp\left[yt-\frac{t^3}{3}\right] \frac{{\rm d}t}{t^p} ,\qquad p \in \Z ,\\
\mathrm{B}_k(y,p) = \int_{\mathscr{b}_k} \exp\left[yt-\frac{t^3}{3}\right] \frac{{\rm d}t}{t^p} ,\qquad p \in \Z_{\le 0},\quad k \in \{1,2,3\},
\end{gather*}
which are defined by integrals along the contours $\mathscr{a}_k$, $\mathscr{b}_0$, $\mathscr{b}_k$ shown in Figure~\ref{fig:path-generalized-Airy} (see \cite{drazin-reid} and \cite[Section~9.13\,(ii)]{dlmf}).

\noindent
These functions satisfy the generalized Airy equation
\[ \left[\left(\frac{\partial}{\partial y}\right)^3 - y \frac{\partial}{\partial y} + (p-1)\right] \Phi = 0. \]
When $p = 0$, this equation reduces to
\[ \frac{\partial}{\partial y} \circ \left[\left(\frac{\partial}{\partial y}\right)^2 - y\right] \Phi = 0, \]
which is equivalent to an inhomogeneous version of the classical Airy equation.

As we did with the Airy--Lucas functions, we will use the substitution $t = 2uy^{1/2}$ to rewrite the generalized Airy function $\mathrm{A}_1(y,p)$ in terms of a thimble integral. We get
\[ \mathrm{A}_1(y,p) = (12 z)^{(1-p)/3} I_+\left(\frac{2}{3}y^{3/2},p\right) \]
with
\[ I_+(z,p) = \frac{1}{2 \pi {\rm i}} \int_{\mathcal{C}^\theta_1} \exp\left[-z\bigl(4u^3 - 3u\bigr)\right] \frac{{\rm d}u}{u^p}. \]
To find the right contour $\mathcal{C}_1^\theta$, we first note that the mapping $\zeta = 4u^3 - 3u$ is familiar, up to a sign, from Section~\ref{sec:airy}. In particular, we already know its critical points $u = \pm\frac{1}{2}$ and the corresponding critical values $\zeta = \pm 1$. To get the desired integral, $\mathcal{C}_1^\theta$ should be the thimble through $u = -\frac{1}{2}$ over the ray $\mathcal{J}^\theta_1$. However, this thimble does not live on the complex plane, like it did for the classical Airy function. When $p$ is a positive integer, the volume form ${\rm d}u/u^p$ has a~pole at $u = 0$, so we put $\mathcal{C}_1^\theta$ as a contour on $\C^\times$. More generally, for any complex value of $p$, we can put $\mathcal{C}_1^\theta$ on the universal cover \smash{$\widetilde{\C^\times}$}. We will now show, in the case $p = 1$, that equation~\eqref{eqn:I} has a frame of Borel regular solutions---as long as the constant solution $\mathrm{B}_0(z, 1)$ counts as Borel regular. For $\mathrm{A}_1(y, p)$ to satisfy the generalized Airy equation, $I_+(z,1)$ must satisfy the equation
\begin{gather}\label{eqn:I}
\left[\left[\left(\frac{\partial}{\partial z}\right)^2 - 1\right] + z^{-1} \left(\frac{\partial}{\partial z}\right)^2 - \frac{1}{9} z^{-2} \right] \frac{\partial}{\partial z} \Phi = 0.
\end{gather}
Equivalently, $\frac{\partial}{\partial z} I_+(z, 1)$ satisfies the modified Bessel equation \eqref{eqn:reg-mod-bessel-AL} with parameter $1/3$. We can now follow the reasoning of Section~\ref{big-idea}, viewing that frequency domain differential equation~\eqref{eqn:I} as the image under $\laplace_{\zeta, \alpha}$ of the position domain integral equation
\begin{gather}\label{int-eq:position-I}
\left[ \left[ \zeta^2 - 1 \right] - \fracderiv{-1}{\zeta}{\alpha} \circ \zeta - \left(\frac{1}{3}\right)^2 \fracderiv{-2}{\zeta}{\alpha} \right] (-\zeta) v^{(1)} = 0.
\end{gather}
In Sections \ref{pos-root-AL} and \ref{neg-root-AL}, we found solutions $v_{\pm 1}$ of the closely related equation~\eqref{int-eq:spatial-mod-bessel-AL}, which we can turn into solutions \smash{$v^{(1)}_{\pm 1} = \zeta^{-1} v_{\pm 1}$} of equation~\eqref{int-eq:position-I}.

At the critical values $\zeta = \pm 1$, we can build solutions \smash{$v^{(1)}_{\pm 1}$} of this equation from the solutions~$v_{\pm 1}$ of equation~\eqref{int-eq:spatial-mod-bessel-AL} that we found in Sections \ref{pos-root-AL} and \ref{neg-root-AL}. Explicitly, in terms of the coordinates defined by $\zeta = 1 + \zeta_1$ and $\zeta = -1 + \zeta_{-1}$, these solutions are
\begin{gather*}
v^{(1)}_1 = -{\rm i}\sqrt{2} (1 + \zeta_1)^{-1} \zeta_1^{-1/2} {}_2F_1\left(\frac{1}{6},\frac{5}{6};\frac{1}{2};-\frac{1}{2}\zeta_{1}\right), \\
v^{(1)}_{-1} = \sqrt{2} (-1 + \zeta_{-1})^{-1} \zeta_{-1}^{-1/2} {}_2F_1\left(\frac{1}{6},\frac{5}{6};\frac{1}{2};\frac{1}{2}\zeta_{-1}\right).
\end{gather*}
Their Laplace transforms \smash{$V^{(1)}_1 = \laplace_{\zeta, 1} v^{(1)}_1$} and \smash{$V^{(1)}_{-1} = \laplace_{\zeta, -1} v^{(1)}_{-1}$} satisfy equation~\eqref{eqn:I}. Notice~that \smash{$v^{(1)}_1$} belongs to \smash{$\singexpalg{-1/2}(\Omega_1)$}, has a~simple pole at $\zeta_1 = -1$, and has a~logarithmic singularity~at $\zeta_1 = -2$. Similarly, \smash{$v^{(1)}_{-1}$} belongs to \smash{$\singexpalg{1/2}(\Omega_{-1})$}, has a~simple pole at $\zeta_{-1} = 1$, and has a~logarithmic singularity at $\zeta_{-1} = 2$.

With an argument analogous to the one in Section~\ref{bessel-regularity-AL}, one should be able to show that~\smash{$V^{(1)}_1$} and~\smash{$V^{(1)}_{-1}$} are Borel regular. Since equation~\eqref{eqn:I} is third-order, we need one more Borel regular solution to make a frame. We choose the constant solution $1$, which can be seen as the Laplace transform of the formal convolution unit $\delta$ introduced in Section~\ref{sec:borel-laplace-homom}. Interestingly, $1$ is actually the generalized Airy function $\mathrm{B}_0(z,1)$.
\subsection{Third-degree thimble integrals}\label{sec:deg3}
In this subsection, we consider the thimble integral
$I(z) = \int_{\mathcal{C}} {\rm e}^{-zg}\, {\rm d}u$
defined by a general third-degree polynomial map $g \maps \C \to \C$. By definition, the Lefschetz thimble $\mathcal{C}$ is a connected component of the preimage of a critical value of $g$. After taking various symmetries into account, we will find that there are essentially only two cases, distinguished by whether $g$ has two distinct critical points or a single degenerate one.
\subsubsection{Symmetries of the integral}\label{sec:int-symm}
First, we step back and consider the completely general thimble integral
\[
J(z) = \int_{\mathcal{C}} {\rm e}^{-zf} \nu
\]
defined by an arbitrary holomorphic map $f \maps X \to \C$ and an arbitrary $1$-form $\nu$ on a complex manifold~$X$. The affine group of $\C$ acts on $f$ by post-composition, and this action affects $J$ in a simple way. The affine group is generated by two kinds of transformations, which we consider separately.

\textbf{Translations.}
For any $c \in \C$, the translation action $f \mapsto f + c$ on maps corresponds to the action
$J(z) \mapsto {\rm e}^{-cz} J(z)$
on integrals. This is because
\begin{align*}
\int_\mathcal{C} {\rm e}^{-z(f + c)} \nu &= {\rm e}^{-cz} \int_\mathcal{C} {\rm e}^{-zf} \nu = {\rm e}^{-cz} J(z).
\end{align*}

\textbf{Scaling-rotations.}
For any $r \in \C^\times$, the scaling-rotation action $f \mapsto rf$ corresponds to the action
$J(z) \mapsto J(rz)$
on integrals. This is because
\begin{align*}
\int_\mathcal{C} {\rm e}^{-z(rf)} \nu & = \int_\mathcal{C} {\rm e}^{-(rz)f} \nu = J(rz).
\end{align*}

Next, we specialize to the case where $X$ is the complex plane and $\nu = {\rm d}u$, where $u$ is the standard coordinate. The affine group of $\C$ now acts on $f$ by pre-composition as well. To make sure that $\mathcal{C}$ remains a Lefschetz thimble, we act on it simultaneously by the inverse of whatever we precompose with $f$.

 \textbf{Translations.}
For any $c \in \C$, pulling $f$ back along the translation $\mathsf{T}_c^* u = u - c$ has no effect on the integral. This is because ${\rm d}u$ is translation-invariant.

\textbf{Scaling-rotations.}
For any $s \in \C^\times$, the scaling-rotation action $\mathsf{M}_s^* u = s^{-1} u$ corresponds to the action
$J(z) \mapsto s J(z)$
on integrals. This is because
\begin{align*}
\int_{\mathsf{M}_{1/s} \mathcal{C}} {\rm e}^{-z \mathsf{M}_s^* f}\, {\rm d}u = \int_\mathcal{C} {\rm e}^{-zf} \mathsf{M}_{1/s}^* \, {\rm d}u
 = \int_\mathcal{C} {\rm e}^{-zf} s \, {\rm d}u
 = s J(z).
\end{align*}

\subsubsection{Reduction of the integral}
We now return to our original thimble integral $I(z)$, defined by the third-degree polynomial map~${ g = a_3 u^3 + a_2 u^2 + a_1 u + a_0}$,
and use the group actions from Section~\ref{sec:int-symm} to reduce it to its two essential cases. We first put $g$ in depressed form by scaling the leading coefficient to $1$ and translating the mean of the critical points to $0$.\footnote{The coefficients of the depressed form are related to the original coefficients as follows
\smash{$
p = \frac{a_1}{a_3} - \frac{1}{3} \big(\frac{a_2}{a_3}\big)^2$}, \smash{$
q = \frac{a_0}{a_3} - \frac{1}{3} \big(\frac{a_1}{a_3}\big)\big(\frac{a_2}{a_3}\big) + \frac{2}{27} \big(\frac{a_2}{a_3}\big)^3$}.
} We then also translate the mean of the critical values to $0$.
\begin{center}
\begin{tabular}{c|c|c}
\textbf{Symmetry} & \textbf{Resulting polynomial} & \textbf{Resulting integral} \\[2mm]
$\displaystyle f \mapsto \frac{f}{a_3}$ & $\displaystyle u^3 + \frac{a_2}{a_3} u^2 + \frac{a_1}{a_3} u + \frac{a_0}{a_3}$ & $\displaystyle \hphantom{{\rm e}^{qz}} I\left(\frac{z}{a_3}\right)$ \\[5mm]
$\displaystyle \mathsf{T}_{(1/3) a_2/a_3}^*$ & $\displaystyle u^3 + pu + q$ & $\displaystyle \hphantom{{\rm e}^{qz}} I\left(\frac{z}{a_3}\right)$ \\[5mm]
$\displaystyle f \mapsto f - q$ & $\displaystyle u^3 + pu$ & $\displaystyle {\rm e}^{qz} I\left(\frac{z}{a_3}\right)$ \\[5mm]
\end{tabular}
\end{center}
Our next move depends on how many critical points $g$ has.

 \textbf{Two distinct critical points.}
When $p$ is non-zero, $g$ has two distinct critical points. We can use our symmetries turn it into any other third-degree polynomial with distinct critical points. To illustrate, we turn~$g$ into the Chebyshev polynomial $T_3$ by scaling and rotating its critical points to $\pm\frac{1}{2}$ and its critical values to $\pm 1$. Setting
\smash{$ r := \frac{{\rm i}}{2} \big(\frac{3}{p}\big)^{3/2}$},
we get
\begin{center}
\begin{tabular}{c|c|c}
\textbf{Symmetry} & \textbf{Resulting polynomial} & \textbf{Resulting integral} \\[2mm]
$\displaystyle \mathsf{M}_{rp/3}^*$ & $\displaystyle r^{-1} \left(-4u^3 + 3u\right)$ & $\displaystyle \frac{rp}{3} {\rm e}^{qz} I\left(\frac{z}{a_3}\right)$ \\[5mm]
$\displaystyle f \mapsto rf$ & $\displaystyle -4u^3 + 3u$ & $\displaystyle \frac{rp}{3} {\rm e}^{qrz} I\left(\frac{rz}{a_3}\right)$ \\[5mm]
\end{tabular}
\end{center}

We can now relate the general third-degree thimble integral $I$ to the modified Bessel function $K_{1/3}$, defined in equation~\eqref{integral:mod-bessel-rational-AL}
\[
K_{1/3}(z) = \frac{{\rm i}rp}{\sqrt{3}} {\rm e}^{qrz} I\left(\frac{rz}{a_3}\right).
\]
From this identity, we can see that $I(z)$ satisfies a transformed version of the modified Bessel equation
\[ \left[ \left[\left(\frac{\partial}{\partial z}\right)^2 + 2a_3 q \frac{\partial}{\partial z} + a_3^2 \big(q^2 - r^{-2}\big)\right] + z^{-1} \left(a_3 q + \frac{\partial}{\partial z}\right) - \frac{1}{9} z^{-2} \right] \Phi = 0. \]

\textbf{One degenerate critical point.}
When $p$ is zero, $g$ has a single degenerate critical point. Like we did in Appendix~\ref{airy-appendix} for the non-degenerate case, we will first study $I(z)$ by viewing it as a solution of an ODE, and we will then confirm our results using the thimble projection formula.

By differentiating under the integral, we can see that $I(z)$ satisfies the equation
\[
\left[\frac{\partial}{\partial z} + \frac{1}{3z} +q \right] \Phi = 0.
\]
A function $\laplace_{\zeta,\alpha}^{\theta} v$ satisfies this equation if and only if $v$ satisfies the integral equation
\[ \left[ -\zeta + \frac{1}{3} \partial^{-1}_{\zeta, \alpha} +q \right] v = 0. \]
When $\alpha = q$, this equation has the solution
$ v_q = \zeta_q^{-2/3}$,
whose Laplace transform
\[ V_q = {\rm e}^{-zq} \Gamma\left(\frac{1}{3}\right) z^{-1/3} \]
must be proportional to $I(z)$.

We can confirm this result, and find the correct normalization, using the thimble projection formula
\begin{align*}
\int_{\mathcal{C}} {\rm e}^{-z(u^3+q)} \, {\rm d}u = \bigl(1 - {\rm e}^{{\rm i}\pi/3}\bigr) \int_{\mathcal{J}_{\zeta_q, 0}^\theta} {\rm e}^{-z\zeta} \frac{1}{3}\zeta_q^{-2/3} \, {\rm d}\zeta_q = {\rm e}^{-zq}\frac{1}{{\rm i}\sqrt{3}} {\rm e}^{{\rm i}\pi/3} \Gamma\left(\frac{1}{3}\right) z^{-1/3}.
\end{align*}
The factor of ${\rm e}^{{\rm i}\pi/3} - 1$ comes from the monodromy of $\zeta_q^{-2/3}$ around $\zeta_q = 0$.

\subsubsection{Coalescence of critical values}
We can recover the degenerate case by taking the limit of the non-degenerate case as $p$ goes to zero. Rearranging the relationship above between $I(z)$ and the modified Bessel function $K_{1/3}(z)$, we get the formula
\[ I(z) = \frac{\sqrt{3}}{{\rm i}rp} {\rm e}^{-qa_3z} K_{1/3}\left(\frac{a_3 z}{r}\right). \]
Setting $q=0$, and using the Taylor expansion of $K_{1/3}(z)$ around $z = 0$, given by \cite[equations~(10.25.2) and~(10.27.4)]{dlmf}, we can see that $I(z)$ converges pointwise to \smash{$\frac{1}{{\rm i}\sqrt{3}} {\rm e}^{{\rm i}\pi/3} \Gamma\bigl(\frac{1}{3}\bigr) z^{-1/3}$} as~$p$ goes to zero with~$a_3$ held constant.
\subsection{The triangular cantilever}\label{sec:catilever}
\subsubsection{Setting}
A {\em triangular cantilever} is a flexible strip with constant thickness and linearly tapered width, clamped at the broad end so it sticks out horizontally like a diving board.
\begin{center}
\begin{tikzpicture}
\newcommand{\clen}{6}
\newcommand{\ctaper}{1}
\newcommand{\cthic}{0.1}
\fill[pwbeige!10] (0, \cthic, 0)--++(6, 0, -\ctaper)--++(0, -2*\cthic, 0)--++(0, 0, 2*\ctaper)--(0, -\cthic, 0);
\draw (0, -\cthic, 0)--++(\clen, 0, \ctaper)--++(0, 0, -2*\ctaper);
\draw[pwbeige] (0, 0, 0)--++(\clen, 0, 0);
\draw[->] (\clen, 0, 0)--++(1, 0, 0) node[anchor=west] {$z$};
\draw (0, \cthic, 0)--++(\clen, 0, \ctaper)--++(0, 0, -2*\ctaper)--cycle;
\foreach \x / \y in {0 / 0, 6 / 1, 6 / -1} {
 \draw (\x, 0.1, \y)--(\x, -0.1, \y);
}
\end{tikzpicture}
\end{center}
If you strike the strip from above, it vibrates up and down. Let us suppose the vibrations are small, and uniform across the width of the strip, so we can describe them using the Euler--Bernoulli beam model~\cite[Section 12.4]{genta2009vibration}. The vibration modes of frequency $\omega$ are described by the vertical displacement profiles $\Phi$ that satisfy the equation
\begin{equation}\label{eqn:triangular_cantilever}
 \left[\left[\left(\frac{\partial}{\partial z}\right)^4 - \omega^2\right] + \frac{2}{z}\left(\frac{\partial}{\partial z}\right)^3\right] \Phi = 0,
\end{equation}
where $z$ is the distance from the tip along the strip's axis.\footnote{To keep the equation simple, we have adjusted the units of time so that the strip's elasticity, density, and taper are absorbed into the frequency parameter.}
\subsubsection{Solving the differential equation}\label{sec:solve-cantilever-eqn}
We seek functions $v$ and points $\zeta = \alpha$ for which the Laplace transform $\laplace_{\zeta, \alpha} v$ satisfies the differential equation~\eqref{eqn:triangular_cantilever}. From Proposition~\ref{prop:L-int-op}, we deduce that $\laplace_{\zeta, \alpha} v$ satisfies equation~\eqref{eqn:triangular_cantilever} if and only if $v$ satisfies the integral equation
\begin{equation}\label{eqn:position_cantilever}
\bigl[ \bigl[ \zeta^4 - \omega^2 \bigr] - 2\fracderiv{-1}{\zeta}{\alpha} \circ \zeta^3 \bigr] v = 0.
\end{equation}
Observing that
\[ \frac{\partial}{\partial \zeta} \sqrt{\zeta^4 - \omega^2} = \frac{2\zeta^3}{\sqrt{\zeta^4 - \omega^2}}, \]
we learn that
\[
v_\text{uni} = \frac{1}{\sqrt{\zeta^4 - \omega^2}}
\]
satisfies equation~\eqref{eqn:position_cantilever} whenever $\alpha^4 - \omega^2 = 0$. Thus, a single ``universal solution'' in the position domain leads to four linearly independent solutions $V_\alpha = \laplace_{\zeta, \alpha} v_\text{uni}$ of equation~\eqref{eqn:triangular_cantilever}, indexed by the fourth roots of $\omega^2$.
\subsubsection{Expressing the solutions as thimble integrals}
If we had a thimble integral expression for $V_\alpha$, we could use the thimble projection formula, as written in equation~\eqref{eqn:formula-proof}, to find $v_\text{uni}$. Since we found $v_\text{uni}$ in Section~\ref{sec:solve-cantilever-eqn}, we can work backwards from equation~\eqref{thimble-difference} to guess a thimble integral expression for $V_\alpha$. We hope to express~$V_\alpha$ in the form
\smash{$\int_{\mathcal{C}_a^\theta}{\rm e}^{-z f}\, {\rm d}u$},
where the function $f$ is a holomorphic function on $\C$, the function~$u$ is the standard coordinate on $\C$, the point $a$ is a critical point of $f$ with $f(a) = \alpha$, and the path~$\mathcal{C}_a^\theta$ is the component of $f^{-1}\bigl(\alpha + {\rm e}^{{\rm i}\theta}[0, \infty)\bigr)$ that runs through~$a$. We know from equation~\eqref{thimble-difference} that this thimble integral will evaluate to $V_\alpha$ as long as
\begin{equation}\label{eqn:cantilever-thimble-difference}
\frac{1}{\sqrt{\zeta^4-\omega^2}} = \left[\frac{{\rm d}u}{{\rm d}f}\right]_{\operatorname{start} \mathcal{C}_a^\theta(\zeta)}^{\operatorname{end} \mathcal{C}_a^\theta(\zeta)},
\end{equation}
where $\mathcal{C}_a^\theta(\zeta)$ is the part of $\mathcal{C}_a^\theta$ that goes through $f^{-1}\bigl(\bigl[\alpha,\zeta {\rm e}^{{\rm i}\theta}\bigr]\bigr)$. We can guess a function $f$ with this property by noting that equation~\eqref{eqn:cantilever-thimble-difference} looks like one of the differential equations that characterize the Jacobi elliptic functions.

The Jacobi function
$x = \operatorname{sd}(t, k)$
with elliptic modulus $k \in [0, 1]$ satisfies the differential equation~\cite[equation~(22.13.5)]{dlmf}
\[ \frac{1}{\sqrt{\bigl(1 - k'^2 x^2\bigr)\bigl(1 + k^2 x^2\bigr)}} = \frac{{\rm d}t}{{\rm d}x}, \]
where $k' \in [0, 1]$ is the complementary modulus given by $k^2 + k'^2 = 1$. In particular, the function~\smash{$ x = \operatorname{sd}\big(t, \frac{1}{\sqrt{2}}\big)$}
satisfies the equation
\[
 \frac{1}{\sqrt{1 - \frac{1}{4} x^4}} = \frac{{\rm d}t}{{\rm d}x}.
 \]
Introducing, for convenience, the constant $\varepsilon = (\omega/2)^{1/2}$, we can deduce that the function
\[ f = {\rm i}\varepsilon \operatorname{sd}\left(4\varepsilon u, \frac{1}{\sqrt{2}}\right) \]
satisfies the equation
\[ \frac{1}{2} \left[ \frac{1}{\sqrt{f^4 - \omega^2}} \right] = \frac{{\rm d}u}{{\rm d}f} \]
when we choose the branch of the square root that evaluates to ${\rm i}\omega^2$ when $f = 0$. By definition, $f\bigl(\operatorname{end} \mathcal{C}_a^\theta(\zeta)\bigr) = \zeta$
so we have
\begin{equation}\label{eqn:cantilever-endpoint}
\frac{1}{2}\left[\frac{1}{\sqrt{\zeta^4 - \omega^2}}\right] = \left.\frac{{\rm d}u}{{\rm d}f}\right|_{\operatorname{end} \mathcal{C}_a^\theta(\zeta)},
\end{equation}
with the branch of the square root determined by analytic continuation along $\mathcal{C}_a^\theta$ from $a$.

From our formula for $f$ in terms of $\operatorname{sd}$, we can see that $f$ is odd with respect to $u$, implying that
\[ \left.\frac{{\rm d}u}{{\rm d}f}\right|_{\operatorname{start} \mathcal{C}_a^\theta(\zeta)} = -\left.\frac{{\rm d}u}{{\rm d}f}\right|_{\operatorname{end} \mathcal{C}_a^\theta(\zeta)}. \]
We can therefore rewrite equation~\eqref{eqn:cantilever-endpoint} as
\[ \frac{1}{\sqrt{\zeta^4-\omega^2}} = \left[\frac{{\rm d}u}{{\rm d}f}\right]_{\operatorname{start} \mathcal{C}_a^\theta(\zeta)}^{\operatorname{end} \mathcal{C}_a^\theta(\zeta)}, \]
confirming that our guess for $f$ satisfies equation~\eqref{eqn:cantilever-thimble-difference}. It follows that
\[ V_\alpha = \int_{\mathcal{C}_a^\theta} \exp\left[-z {\rm i}\varepsilon \operatorname{sd}\left(4\varepsilon u, \frac{1}{\sqrt{2}}\right) \right] {\rm d}u, \]
where $\varepsilon = (\omega/2)^{1/2}$ is the convenience constant we introduced earlier. The exponent is doubly periodic, with period lattice~\cite[Table~22.4.1]{dlmf}
\[ L = \varepsilon^{-1} K\left(\frac{1}{\sqrt{2}}\right)\left[\Z + \frac{1}{2}(1 + {\rm i}) \Z\right], \]
where $K$ is the complete elliptic integral~\cite[formula~(19.2.8)]{dlmf}
\[ K(k) = \int_0^1 \frac{{\rm d}x}{\sqrt{\bigl(1 - x^2\bigr)\bigl(1 - k^2 x^2\bigr)}}. \]
We can therefore take the thimble to be a real submanifold of the torus $\C / L$.

\appendix

\section{The Airy equation}\label{airy-appendix}
\subsection{Specializing from the Airy--Lucas example}\label{sec:spec-to-airy}
\subsubsection{Motivation}
Since the Airy equation is a widely used example in the study of Borel summation and resurgence, we find it useful to repeat the Airy--Lucas example calculations from Sections~\ref{sec:airy}--\ref{example_AL} in this special case, with the parameters set to $n = 3$ and $m = 1$. In particular, this makes it easier to compare our approaches and conventions with others from the literature, as we do in Appendix~\ref{airy-comparison}.

\subsubsection{Rewriting as a modified Bessel equation}
Starting from equation~\eqref{eqn:airy-distillation}, we can distill the most interesting part of the Airy function by writing
\smash{$\Ai(y) = -\frac{1}{\pi\sqrt{3}} y^{1/2} K\big(\frac{2}{3} y^{3/2}\big)$},
where
\begin{equation}\label{integral:mod-bessel-airy}
K(z) = -{\rm i}\sqrt{3} \int_{\mathcal{C}^\theta_1} \exp\bigl[z \bigl(4u^3 - 3u\bigr)\bigr]\, {\rm d}u
\end{equation}
and $\mathcal{C}^\theta_1$ is the contour described in Section~\ref{sec:airy}. Saying that $\Ai$ satisfies the Airy equation is equivalent to saying that $K$ satisfies the modified Bessel equation
\[
\left[z^2 \left(\frac{\partial}{\partial z}\right)^2 + z \frac{\partial}{\partial z} - \left[\left(\frac{1}{3}\right)^2 + z^2\right]\right] K = 0.
\]
In fact, $K$ is the modified Bessel function $K_{1/3}$~\cite[equation~(9.6.1)]{dlmf}. Like we did in equation~\eqref{eqn:reg-mod-bessel-AL}, we can rewrite the modified Bessel equation above as
\begin{equation}\label{eqn:reg-mod-bessel}
\left[ \left[ \left(\frac{\partial}{\partial z}\right)^2 - 1 \right] + z^{-1} \frac{\partial}{\partial z} - \left(\frac{1}{3}\right)^2 z^{-2} \right] K = 0.
\end{equation}
\subsubsection{Asymptotic analysis}\label{sec:asympt-airy}
We know from the formal theory of level~$1$ ODEs that equation~\eqref{eqn:reg-mod-bessel} has a frame of formal $1$-Gevrey trans-monomial solutions
$\big\{ {\rm e}^{-\alpha z} z^{-\tau_\alpha} \series{W}_\alpha \mid \alpha^2 - 1 = 0 \big\}$,
where $\tau_\alpha = 1/2$ and $\series{W}_\alpha\in\C\big\llbracket z^{-1} \big\rrbracket_1$. Specializing the solution formulas given in Section~\ref{sec:asympt-AL}, we find that $\smash{K \sim \bigl(\frac{\pi}{2}\bigr)^{1/2}} {\rm e}^{-z} \allowbreak\times z^{-1/2} \series{W}_1$, with
\[
\series{W}_1 = 1 - \frac{\bigl(\frac{1}{6}\bigr)_1 \bigl(\frac{5}{6}\bigr)_1}{2^1 \cdot 1!} z^{-1} + \frac{\bigl(\frac{1}{6}\bigr)_2 \bigl(\frac{5}{6}\bigr)_2}{2^2 \cdot 2!} z^{-2} - \frac{\bigl(\frac{1}{6}\bigr)_3 \bigl(\frac{5}{6}\bigr)_3}{2^3 \cdot 3!} z^{-3} + \cdots.
\]
The analogous holomorphic solutions
$ \big\{ {\rm e}^{-\alpha z} z^{-\tau_\alpha} W_\alpha \mid \alpha^2 - 1 = 0 \big\}$,
which we will describe in Appendix~\ref{big-idea-airy}, can be recovered from the formal solutions using Borel summation. We showed this abstractly in Section~\ref{bessel-regularity-AL}, and we will confirm it concretely in Appendix~\ref{confirmation-borel-regularity-airy}. This confirms Theorem~\ref{thm:summability_ODE}.

\subsubsection{The big idea}\label{big-idea-airy}
As discussed in Section~\ref{big-idea}, a function $\laplace_{\zeta, \alpha} v$ satisfies the differential equation~\eqref{eqn:reg-mod-bessel} if and only if $v$ satisfies the integral equation
\begin{equation}\label{int-eq:spatial-mod-bessel}
\left[ \left[ \zeta^2 - 1 \right] - \fracderiv{-1}{\zeta}{\alpha} \circ \zeta - \left(\frac{1}{3}\right)^2 \fracderiv{-2}{\zeta}{\alpha} \right] v = 0.
\end{equation}
Every solution of equation~\eqref{int-eq:spatial-mod-bessel} will also satisfy the differential equation
\begin{equation}\label{diff-eq:spatial-mod-bessel}
\left[ \left(\frac{\partial}{\partial \zeta}\right)^2 \circ \bigl[ \zeta^2 - 1 \bigr] - \frac{\partial}{\partial \zeta} \circ \zeta - \left(\frac{1}{3}\right)^2 \right] v = 0
\end{equation}
obtained by differentiating twice on both sides. Conversely, as shown in Appendix~\ref{shifting}, a solution of equation~\eqref{diff-eq:spatial-mod-bessel} will satisfy equation~\eqref{int-eq:spatial-mod-bessel} if it belongs to $(\zeta - \alpha)^\sigma \mathcal{O}_{\zeta = \alpha}$ for some $\sigma \in (-1, 0)$. Rewriting equation~\eqref{diff-eq:spatial-mod-bessel} as
\[ \left[ (\zeta^2 - 1) \left(\frac{\partial}{\partial \zeta}\right)^2 + 3\zeta \frac{\partial}{\partial \zeta} + \left[ 1 - \left(\frac{1}{3}\right)^2 \right] \right] v = 0 \]
emphasizes the regular singularities at the roots of $\zeta^2 - 1$. This will lead us, in Appendices~\ref{pos-root} and~\ref{neg-root}, to a solution $v_\alpha \in (\zeta - \alpha)^{1/2} \mathcal{O}_{\zeta = \alpha}$ for each root $\alpha$.

As described in Section~\ref{big-idea}, each solution $v_\alpha$ will extend to a function in $\singexpalg{-1/2}(\Omega_\alpha)$ on any sector $\Omega_\alpha$ that has its tip at $\zeta = \alpha$ and does not touch the other roots. The arguments of Section~\ref{big-idea} go on to show that $\laplace^\theta_{\zeta, \alpha} v_\alpha$ belongs to
\smash{$ cz^{-\tau_\alpha} + \dualsingexp{-3/2}\bigl(\widehat{\Omega}_\alpha^\Lambda\bigr) $}
for some non-zero constant $c$, confirming the existence part of Theorem~\ref{re:thm:exist_uniq_ODE}. Continuing through Section~\ref{big-idea}, we get the decomposition
\smash{$ \laplace^\theta_{\zeta, \alpha} v_\alpha = {\rm e}^{-\alpha z} V_\alpha $}
for $V_\alpha := \laplace^\theta_{\zeta_\alpha, 0} v_\alpha$ and $\zeta = \alpha + \zeta_\alpha$, and the further decomposition
\smash{$\laplace_{\zeta, \alpha} v_\alpha = {\rm e}^{-\alpha z} z^{-1/2} W_\alpha$},
where $W_\alpha$ is a bounded holomorphic function on \smash{$\widehat{\Omega}_\alpha^\Lambda$}.

The argument in Section~\ref{bessel-regularity-AL} shows that the solutions ${\rm e}^{-\alpha z} V_\alpha$ are Borel regular. This argument is just as simple in general as it is in the Airy case, so we will not repeat it.
\subsubsection[Focus on zeta = 1]{Focus on $\boldsymbol{\zeta = 1}$}\label{pos-root}
Switching to the translation coordinate $\zeta_1$ defined by $\zeta = 1 + \zeta_1$, and then to the rescaled coordinate $\xi_1$ defined by $\zeta_1 = -2\xi_1$, we can rewrite equation~\eqref{diff-eq:spatial-mod-bessel} as the hypergeometric equation
\[
\left[\xi_1 (1 - \xi_1) \left(\frac{\partial}{\partial \xi_1}\right)^2 + 3\left(\frac{1}{2} - \xi_1\right) \frac{\partial}{\partial \xi_1} - \left[1 - \left(\frac{1}{3}\right)^2\right]\right] v = 0,
\]
as described in Section~\ref{pos-root-AL}. The solution $v_1 \in (\zeta-1)^{-1/2} \mathcal{O}_{\zeta=1}$ found in that section specifies~to
\begin{gather*}
v_1 = \xi_1^{-1/2} {}_2F_1\left(\frac{1}{6}, \frac{5}{6}; \frac{1}{2}; \xi_1\right)
= -{\rm i}\sqrt{2} \zeta_1^{-1/2} {}_2F_1\left(\frac{1}{6}, \frac{5}{6}; \frac{1}{2}; -\frac{1}{2}\zeta_1\right)
\end{gather*}
in the Airy equation case. As discussed in Section~\ref{pos-root-AL}, $v_1$ is holomorphic throughout the sector~${\Omega_1 = \C \smallsetminus \mathcal{J}^\pi_{\zeta_1, 0}}$, and it belongs to the space
\smash{$ -{\rm i}\sqrt{2} \zeta_1^{-1/2} + \singexp{1/2}{\Lambda}(\Omega_1) $}
for all $\Lambda > 0$. Its Laplace transform can be written as
\smash{$ \laplace^\theta_{\zeta, 1} v_1 = {\rm e}^{-z} z^{-1/2} W_1$},
where $W_1$ is a holomorphic function on the slit plane $z \notin (-\infty, 0]$ which is bounded outside any constant-radius neighborhood of the ray $z \in (-\infty, 0]$.
\subsubsection[Focus on zeta = -1]{Focus on $\boldsymbol{\zeta = -1}$}\label{neg-root}
We now switch to the translation coordinate $\zeta_{-1}$ defined by $\zeta = -1 + \zeta_{-1}$. The solution $v_{-1} \in (\zeta+1)^{-1/2} \mathcal{O}_{\zeta=-1}$ found in Section~\ref{neg-root-AL} specifies to
\begin{gather*}
v_{-1} = (1-\xi_1)^{-1/2} {}_2F_1\left(\frac{1}{6}, \frac{5}{6}; \frac{1}{2}; 1-\xi_1\right)
= \sqrt{2} \zeta_{-1}^{-1/2} {}_2F_1\left(\frac{1}{6}, \frac{5}{6}; \frac{1}{2}; \frac{1}{2}\zeta_{-1}\right)
\end{gather*}
in the Airy equation case. By the same reasoning as before, $v_{-1}$ is holomorphic throughout the sector $\Omega_{-1} = \C \smallsetminus \mathcal{J}^0_{\zeta_{-1}, 0}$, and it belongs to the space
\smash{$ \sqrt{2} \zeta_{-1}^{-1/2} + \singexp{1/2}{\Lambda}(\Omega_1) $}
for all $\Lambda > 0$. Its Laplace transform can be written as
\smash{$ \laplace^\theta_{\zeta, -1} v_1 ={\rm e}^z z^{-1/2} W_{-1}$},
where $W_{-1}$ is a holomorphic function on the slit plane $z \notin [0, \infty)$ which is bounded outside any constant-radius neighborhood of the ray $z \in [0, \infty)$.
\subsubsection{Confirmation of Borel regularity}\label{confirmation-borel-regularity-airy}
We can verify the conclusions of Section~\ref{bessel-regularity-AL} in the Airy case using our explicit expressions for the formal power series $\tilde{W}_\alpha$ and the functions $v_\alpha$. We found in Appendix~\ref{sec:asympt-airy} that
\begin{align*}
\tilde{W}_1 = \sum_{k = 0}^{\infty} \frac{\bigl(\frac{1}{6}\bigr)_k \bigl(\frac{5}{6}\bigr)_k}{k!} \left(-\frac{1}{2}\right)^k z^{-k} ,\qquad
\tilde{W}_{-1} = \sum_{k = 0}^{\infty} \frac{\bigl(\frac{1}{6}\bigr)_k \bigl(\frac{5}{6}\bigr)_k}{k!} \left(\frac{1}{2}\right)^k z^{-k}.
\end{align*}
Repeating the computation in Section~\ref{confirmation-borel-regularity}, we find that
\[ \borel_\zeta \bigl[ {\rm e}^{-z} z^{-1/2} \tilde{W}_1 \bigr] = \frac{\zeta_1^{-\frac{1}{2}}}{\Gamma\bigl(\frac{1}{2}\bigr)} \sum_{k = 0}^{\infty} \frac{\bigl(\frac{1}{6}\bigr)_k \bigl(\frac{5}{6}\bigr)_k}{\bigl(\frac{1}{2}\bigr)_k} \left(-\frac{1}{2}\right)^k \frac{\zeta_1^k}{k!}, \]
making it apparent that $\borel\bigl[ {\rm e}^{-z} z^{-1/2} \tilde{W}_1 \bigr]$ sums to
\[ \frac{1}{\Gamma(1/2)} \zeta_1^{-1/2} {}_2F_1\left(\frac{1}{6}, \frac{5}{6}; \frac{1}{2}; -\frac{1}{2}\zeta_1\right). \]
Looking back at Appendix~\ref{pos-root}, we recognize this as a scalar multiple of $v_1$.

Through a similar calculation, we see that $\borel\bigl[ {\rm e}^z z^{-1/2} \tilde{W}_{-1} \bigr]$ sums to
\[ \frac{1}{\Gamma(1/2)} \zeta_{-1}^{-1/2} {}_2F_1\left(\frac{1}{2}-\frac{m}{n}, \frac{1}{2}+\frac{m}{n}; \frac{1}{2}; \frac{1}{2}\zeta_{-1}\right). \]
Looking back at Appendix~\ref{neg-root}, we recognize this as a scalar multiple of $v_{-1}$.

\subsubsection{Thimble projection reasoning}\label{contour-argument-airy}

In Section~\ref{contour-argument-AL}, we specialized the reasoning behind Lemma~\ref{lem:thimble_proj_formula-proof} to the case of the Airy--Lucas functions. We now specialize even further, down to the case of the Airy function. Like in Section~\ref{contour-argument-AL}, we first recast integral~\eqref{integral:mod-bessel-airy} into the $\zeta$ plane by setting $-\zeta = 4u^3 - 3u$, which implies that $-{\rm d}\zeta = 3\bigl(4u^2 - 1\bigr){\rm d}u$. Recall from Section~\ref{sec:airy} that the integration contour $\mathcal{C}^\theta_1$ runs through the critical point $u = \frac{1}{2}$, with its direction determined by the parameter $\theta$. The critical point splits $\mathcal{C}^\theta_1$ into two pieces: the incoming branch, where the orientation of the thimble runs toward the critical point, and the outgoing branch, where the orientation points away. Recalling that $\mathcal{C}^\theta_1$ is a preimage of $\mathcal{J}^\theta_{\zeta,1}$, we get
\begin{align*}
K(z) & = -{\rm i}\sqrt{3} \int_{\mathcal{C}^\theta_1} \exp\bigl[z \bigl(4u^3 - 3u\bigr)\bigr] \, {\rm d}u \\
& = \frac{{\rm i}}{\sqrt{3}} \left[ \int_{\mathcal{J}^\theta_{\zeta, 1}} {\rm e}^{-z\zeta} \frac{1}{4u_+^2 - 1} \,{\rm d}\zeta - \int_{\mathcal{J}^\theta_{\zeta, 1}} {\rm e}^{-z\zeta} \frac{1}{4u_-^2 - 1} \, {\rm d}\zeta \right],
\end{align*}
where $u_-$ and $u_+$ are the lifts to the incoming and outgoing branches of $\mathcal{C}^\theta_j$, respectively. Recognizing \smash{$4u^3 - 3u$} as the Chebyshev polynomial $T_3(u)$, we introduce a new variable $\phi$ with $u = \cos(\phi)$ and $-\zeta = \cos(3\phi)$. On the $\phi$ plane, which is an infinite branched cover of the $u$ plane, we can lift $\mathcal{C}^\theta_1$ to the path $\frac{\pi}{3} + {\rm i}\R$. The incoming and outgoing branches lift to $\frac{\pi}{3} - {\rm i}[0, \infty)$ and $\frac{\pi}{3} + {\rm i}[0, \infty)$, respectively.

Recognizing $4u^2 - 1$ as the Chebyshev polynomial $U_2(u)$, we can use identity~(15.4.16) from~\cite{dlmf} to write the integrand explicitly in terms of $\zeta$:
\begin{align*}
\frac{1}{4 \cos(\phi)^2 - 1} & = \frac{1}{U_2(\cos(\phi))} = \frac{\sin(\phi)}{\sin(3\phi)} = \frac{1}{3} {}_2F_1\left(\frac{1}{3}, \frac{2}{3}; \frac{3}{2}; \sin(n \phi)^2\right) \\
& = \frac{1}{3} {}_2F_1\left(\frac{1}{3}, \frac{2}{3}; \frac{3}{2}; 1 - \zeta^2\right).
\end{align*}
Just as in Section~\ref{contour-argument-AL}, we simplify the integrand using identities (15.8.4) and (15.8.27)--(15.8.28) from \cite{dlmf}, which tell us that
\begin{align*}
{}_2F_1\left(\frac{1}{3}, \frac{2}{3}; \frac{3}{2}; 1 - \zeta^2\right)
 ={}&\frac{\pi}{\Gamma\left(\frac{1}{3}\right) \Gamma\left(\frac{2}{3}\right)} {}_2F_1\left(\frac{1}{3}, \frac{2}{3}; \frac{1}{2}; \zeta^2\right) - \frac{\pi \zeta}{\Gamma\left(\frac{1}{3}\right) \Gamma\left(\frac{2}{3}\right)} {}_2F_1\left(\frac{1}{3}, \frac{2}{3}; \frac{3}{2}; \zeta^2\right) \\
={}& {}_2F_1\left(\frac{2}{3}, \frac{4}{3}; \frac{3}{2}; \frac{1}{2} - \frac{\zeta}{2}\right) + \hphantom{\frac{0}{0}} {}_2F_1\left(\frac{2}{3}, \frac{4}{3}; \frac{3}{2}; \frac{1}{2} + \frac{\zeta}{2}\right) \\
& + \frac{1}{2} {}_2F_1\left(\frac{2}{3}, \frac{4}{3}; \frac{3}{2}; \frac{1}{2} - \frac{\zeta}{2}\right) - \frac{1}{2} {}_2F_1\left(\frac{2}{3}, \frac{4}{3}; \frac{3}{2}; \frac{1}{2} + \frac{\zeta}{2}\right) \\
={}&\frac{3}{2} {}_2F_1\left(\frac{2}{3}, \frac{4}{3}; \frac{3}{2};\frac{1}{2} - \frac{\zeta}{2}\right) + \frac{1}{2} {}_2F_1\left(\frac{2}{3}, \frac{4}{3}; \frac{3}{2}; \frac{1}{2} + \frac{\zeta}{2}\right)
\end{align*}
away from the line $\Re(\zeta) = 0$ and the rays $\mathcal{J}^0_{\zeta,1}$ and $\mathcal{J}^\pi_{\zeta,-1}$.

In the projected thimble integral
\[ K(z) = \frac{\rm i}{\sqrt{3}} \int_{\mathcal{J}^\theta_{\zeta, 1}} {\rm e}^{-z\zeta}\left[\frac{1}{4u_+^2 - 1} - \frac{1}{4u_-^2 - 1}\right] {\rm d}\zeta, \]
we can now see the integrand as the variation of the function
\[ \frac{1}{4\cos(\phi)^2 - 1} = \frac{1}{2} {}_2F_1\left(\frac{2}{3}, \frac{4}{3}; \frac{3}{2};\frac{1}{2} - \frac{\zeta}{2}\right) + \frac{1}{6} {}_2F_1\left(\frac{2}{3}, \frac{4}{3}; \frac{3}{2}; \frac{1}{2} + \frac{\zeta}{2}\right) \]
across the branch cut $\mathcal{J}^\theta_{\zeta,1}$. The first term is regular at $\zeta = 1$, so only the second term contributes to the jump. We can write the jump explicitly using identity~(15.2.3) from \cite{dlmf}
\[ K(z) = -\frac{1}{2} \int_{\mathcal{J}^\theta_{\zeta, 1}} {\rm e}^{-z\zeta} \left(-\frac{1}{2}+\frac{\zeta}{2}\right)^{-1/2} {}_2F_1\left(\frac{1}{6}, \frac{5}{6}; \frac{1}{2}; \frac{1}{2} - \frac{\zeta}{2}\right) {\rm d}\zeta. \]
Comparing this expression with the expression for $V_1$ in Appendix~\ref{pos-root}, we see that
$K(z) = \frac{\rm i}{2}\laplace_{\zeta, 1} v_1$.

\subsubsection{Thimble projection formula}\label{thimble-proj-airy}

In Appendix~\ref{contour-argument-airy}, we specialized the reasoning behind Lemma~\ref{lem:thimble_proj_formula-proof} to the case of the Airy function. In this appendix, as an alternative, we will simply apply Lemma~\ref{lem:thimble_proj_formula-proof}, using trigonometric substitution. The lemma tells us that $K$ is the Laplace transform of
\[ \kappa_1 = -{\rm i}\sqrt{3} \frac{\partial}{\partial \zeta}\left(\int_{\mathcal{C}_1^\theta(\zeta)}  {\rm d}u\right). \]
Based on the variable $\phi$ from Appendix~\ref{contour-argument-airy}, we introduce a new variable $\eta$ defined by ${\phi = \frac{\pi}{3} + {\rm i}\eta}$, with $u = \cos\bigl(\frac{\pi}{3} + {\rm i}\eta\bigr)$ and $\zeta = \cos\bigl(\pi + 3{\rm i}\eta\bigr)$. We can lift the thimble $\mathcal{C}^\theta_1$ to the path $\eta \in \R$. Its incoming and outgoing branches lift to $(-\infty, 0]$ and $[0, \infty)$, respectively. Using the variable~$\eta$, we calculate
\begin{align*}
\int_{\mathcal{C}_1^\theta(\zeta)} \, {\rm d}u & = u \Big|_{\operatorname{start} \mathcal{C}^\theta_1(\zeta)}^{\operatorname{end} \mathcal{C}^\theta_1(\zeta)} = u_+ - u_-.
\end{align*}
Like in Section~\ref{contour-argument-AL}, the functions $u_-$ and $u_+$ are functions on a neighborhood of $\mathcal{J}_{\zeta, 1}$ in the position domain. They give the values of $u$ after lifting to the incoming and outgoing branches of $\mathcal{C}^\theta_1$, respectively. We define functions $\eta_-$ and $\eta_+$ similarly, by lifting from the position domain to $\eta \in (-\infty, 0]$ and $\eta \in [0, \infty)$. Observing that $\eta_- = -\eta_+$, we can continue the calculation
\begin{align*}
\int_{\mathcal{C}_1^\theta(\zeta)} {\rm d}u = {}& \cos\left(\frac{\pi}{3} + {\rm i}\eta_+\right) - \cos\left(\frac{\pi}{3} + {\rm i}\eta_-\right)
 = \cos\left(\frac{\pi}{3} +{\rm i}\eta_+\right) - \cos\left(\frac{\pi}{3} - {\rm i}\eta_+\right) \\
 = {}& \left[\cos\left(\frac{\pi}{3}\right) \cos({\rm i}\eta_+) - \sin\left(\frac{\pi}{3}\right) \sin({\rm i}\eta_+)\right] \\
& - \left[\cos\left(\frac{\pi}{3}\right) \cos(-{\rm i}\eta_+) - \sin\left(\frac{\pi}{3}\right) \sin(-{\rm i}\eta_+)\right] \\
 = {}&-2 \sin\left(\frac{\pi}{3}\right) \sin(-{\rm i}\eta_+)
 = {\rm i}\sqrt{3} \sinh(\eta_+).
\end{align*}
Identity (15.4.16) from \cite{dlmf} implies that
\[ \sinh(\eta_+) = \frac{2}{3} \sinh\left(\frac{3}{2}\eta_+\right) {}_2 F_1\left(\frac{1}{6}, \frac{5}{6}; \frac{3}{2}; -\sinh\left(\frac{3}{2}\eta_+\right)^2\right). \]
Noticing that $\frac{1}{2}(1 - \zeta) = -\sinh\bigl(\frac{3}{2}\eta_+\bigr)^2$, we conclude that
\[ \int_{\mathcal{C}_1^\theta(\zeta)} {\rm d}u = \frac{2{\rm i}}{\sqrt{3}} \left(-\frac{1}{2} + \frac{\zeta}{2}\right)^{1/2} {}_2 F_1\left(\frac{1}{6}, \frac{5}{6}; \frac{3}{2}; \frac{1}{2} - \frac{\zeta}{2}\right). \]
It follows, through identity~(15.5.4) from \cite{dlmf}, that
\begin{align*}
\kappa_1 & = 2 \frac{\partial}{\partial \zeta} \left[ \left(-\frac{1}{2} + \frac{\zeta}{2}\right)^{1/2} {}_2 F_1\left(\frac{1}{6}, \frac{5}{6}; \frac{3}{2}; \frac{1}{2} - \frac{\zeta}{2}\right) \right] \\
& = -\frac{1}{2} \left(-\frac{1}{2} + \frac{\zeta}{2}\right)^{-1/2} {}_2 F_1\left(\frac{1}{6}, \frac{5}{6}; \frac{1}{2}; \frac{1}{2} - \frac{\zeta}{2}\right),
\end{align*}
matching the conclusion of Appendix~\ref{contour-argument-airy}.
\subsection{Comparison with other treatments of the Airy equation}\label{airy-comparison}
\subsubsection{Other conventions for the Borel transform}
Physicists often use a different version of the Borel transform
\begin{align*}
\borel_{{\rm phys}} \maps \ \C \big\llbracket z^{-1} \big\rrbracket \to \C \llbracket \zeta \rrbracket ,\qquad
z^{-n} \mapsto \frac{\zeta^n}{n!}.
\end{align*}
This version avoids sending $1$ to the convolution unit $\delta$, at the cost of no longer mapping multiplication to convolution or inverting the formal Laplace transform. It is related to the mathematician's Borel transform by the identity $\borel_{\rm phys}(f) = \borel\bigl(z^{-1} f\bigr)$.

For problems involving a small parameter $\hbar$ rather than a large parameter $z$, physicists also define
\begin{align*}
\borel_{{\rm phys}} \maps\ \C \llbracket \hbar \rrbracket \to \C \llbracket \zeta \rrbracket ,\qquad
\hbar^n \mapsto \frac{\zeta^n}{n!}.
\end{align*}
From a combinatorial perspective, this is just the transformation that sends an ordinary generating function to the corresponding exponential generating function.

In \cite{lectures-Marino}, Mari\~{n}o studies the Airy functions as an example of resurgent functions. He starts with the formal trans-monomial solutions of the Airy equation
\begin{align*}
\series{\Phi}_{\mathrm{Ai}}(x)=\frac{1}{2\sqrt{\pi}}x^{-1/4}{\rm e}^{-\frac{2}{3}x^{3/2}}\tilde{W}_1\bigl(x^{-3/2}\bigr),\qquad
\series{\Phi}_{\mathrm{Bi}}(x)=\frac{1}{2\sqrt{\pi}}x^{-1/4}{\rm e}^{\frac{2}{3}x^{3/2}}\tilde{W}_2\bigl(x^{-3/2}\bigr),
\end{align*}
where
\begin{align*}
\tilde{W}_{1,2}(\hbar)=\sum_{n=0}^{\infty}\frac{1}{2\pi}\left(\mp\frac{3}{4}\right)^{n}\frac{\Gamma\bigl(n+\frac{5}{6}\bigr)\Gamma\bigl(n+\frac{1}{6}\bigr)}{n!}\hbar^n.
\end{align*}
By applying the two versions $\borel_{{\rm phys}}$ and $\borel$ of the Borel transform, we get
\begin{align*}
w_{1,2}(\zeta)&:=\borel_{{\rm phys}}\bigl(\tilde{W}_{1,2}\bigr)(\zeta)=\sum_{n=0}^{\infty}\frac{1}{2\pi}\left(\mp\frac{3}{4}\right)^{n}\frac{\Gamma\bigl(n+\frac{5}{6}\bigr)\Gamma\bigl(n+\frac{1}{6}\bigr)}{n!}\frac{\zeta^n}{n!}\\
&={}_2F_1\left(\frac{1}{6},\frac{5}{6};1;\mp\frac{3}{4}\zeta\right)
\end{align*}
and
\begin{align*}
\borel\bigl(\tilde{W}_{1,2}\bigr)(\zeta)&=\frac{1}{2\pi}\delta+\sum_{n=1}^{\infty} \frac{1}{2\pi}\left(\mp\frac{3}{4}\right)^{n}\frac{\Gamma\bigl(n+\frac{5}{6}\bigr)\Gamma\bigl(n+\frac{1}{6}\bigr)}{n!}\frac{\zeta^{n-1}}{(n-1)!}\\
&=\frac{1}{2\pi}\delta+\sum_{n=0}^{\infty} \frac{1}{2\pi}\left(\mp\frac{3}{4}\right)^{n+1}\frac{\Gamma\bigl(n+1+\frac{5}{6}\bigr)\Gamma\bigl(n+1+\frac{1}{6}\bigr)}{(n+1)!}\frac{\zeta^{n}}{n!}\\
&=\frac{1}{2\pi}\delta\mp\frac{3}{4}\sum_{n=0}^{\infty} \frac{1}{2\pi}\left(\mp\frac{3}{4}\right)^{n}\frac{\Gamma(n+\frac{11}{6})\Gamma\bigl(n+\frac{7}{6}\bigr)}{\Gamma(n+2)}\frac{\zeta^{n}}{n!}\\
&=\frac{1}{2\pi}\delta\mp\frac{5}{48} {}_2F_1\left(\frac{7}{6},\frac{11}{6};2;\mp\frac{3}{4}\zeta\right).
\end{align*}
Comparing the two solutions, we notice that up to the factor of $\delta$,
\begin{equation}\label{Borel-W12}
\borel\bigl(\tilde{W}_{1,2}\bigr)(\zeta)-\frac{1}{2\pi}\delta=\mp \frac{5}{48} {}_2F_1\left(\frac{7}{6},\frac{11}{6};2;\mp\frac{3}{4}\zeta\right)=\frac{\rm d}{{\rm d}\zeta}\borel_{{\rm phys}}\bigl(\tilde{W}_{1,2}\bigr)(\zeta).
\end{equation}
More generally, if $\series{\Phi}(z)\in z^{-1}\C \big\llbracket z^{-1} \big\rrbracket$, i.e., it has no constant term, then
$
\frac{\rm d}{{\rm d}\zeta}\circ\borel_{{\rm phys}}\series{\Phi}=\borel\series{\Phi}$.
In particular, \smash{$\frac{\rm d}{{\rm d}\zeta}\circ\borel_{{\rm phys}}\bigl[z^{-1/2}\tilde{W}_{1}\bigr]\left(\frac{2}{3}\zeta\right)=v_1(\zeta)$}.
\subsubsection{Integral formula for hypergeometric functions}
In \cite{diverg-resurg-i}, Mitschi and Sauzin study summability and resurgent properties of solutions of the Airy equation. They consider the formal power series
\[\series{\Phi}_{\pm}(z):= \sum_{n=0}^{\infty}\frac{1}{2\pi}\left(\mp\frac{1}{2}\right)^{n}\frac{\Gamma\bigl(n+\frac{5}{6}\bigr)\Gamma\bigl(n+\frac{1}{6}\bigr)}{n!}z^{-n}\]
such that
\begin{align*}
\series{\Phi}_{\mathrm{Ai}}(y)&=\frac{1}{2\sqrt{\pi}}y^{-1/4}{\rm e}^{-\frac{2}{3}y^{3/2}}\series{\Phi}_{+}\left(\frac{2}{3}y^{3/2}\right),\qquad
\series{\Phi}_{\mathrm{Bi}}(y)=\frac{1}{2\sqrt{\pi}}y^{-1/4}{\rm e}^{\frac{2}{3}y^{3/2}}\series{\Phi}_{-}\left(\frac{2}{3}y^{3/2}\right)
\end{align*}
are formal solutions of the Airy equation. Notice that compared to Mari\~{n}o's formal solutions, Mitschi and Sauzin adopt a different change of coordinates $z=\frac{2}{3}y^{3/2}$.

Seeking solutions of the Borel-transformed equation, Mitschi and Sauzin write the Borel transform of $\series{\Phi}_{\pm}$ as a convolution product
\begin{align*}
\hat{\phi}_+(\zeta):=\borel\series{\Phi}_+=\delta+\frac{\rm d}{{\rm d}\zeta}\hat{\chi}(\zeta), \qquad \hat{\phi}_-(\zeta):=\borel\series{\Phi}_-=\delta-\frac{\rm d}{{\rm d}\zeta}\hat{\chi}(-\zeta),
\end{align*}
 where
\[
\hat{\chi}(\zeta)=\frac{2^{1/6}}{\Gamma(1/6)\Gamma(5/6)}\bigl(2\zeta+\zeta^2\bigr)^{-1/6} \ast_\zeta \zeta^{-5/6}.
\]
Notice that $\hat{\chi}(\zeta)$ is a hypergeometric function
\begin{align*}
\hat{\chi}(\zeta)&=\frac{2^{1/6}}{\Gamma(1/6)\Gamma(5/6)}\bigl(2\zeta+\zeta^2\bigr)^{-1/6}\ast_\zeta \zeta^{-5/6}\\
&=\frac{2^{1/6}}{\Gamma(1/6)\Gamma(5/6)}\int_0^{\zeta}\bigl(2\zeta'+\zeta'^2\bigr)^{-1/6} \bigl(\zeta-\zeta'\bigr)^{-5/6}\,{\rm d}\zeta'\\
&=\frac{2^{1/6}}{\Gamma(1/6)\Gamma(5/6)}\int_0^{1}(\zeta t)^{-1/6}(2+\zeta t)^{-1/6} (\zeta-\zeta t)^{-5/6} \zeta \,{\rm d}t\\
&=\frac{2^{1/6}}{\Gamma(1/6)\Gamma(5/6)}\int_0^{1} t^{-1/6} 2^{-1/6}\left(1+\frac{\zeta}{2} t\right)^{-1/6} (1-t)^{-5/6}\,{\rm d}\zeta'\\
&=\frac{1}{\Gamma(1/6)\Gamma(5/6)}\int_0^{1} t^{-1/6} \left(1+\frac{\zeta}{2} t\right)^{-1/6} (1-t)^{-5/6}\,{\rm d}\zeta'\\
&={}_2F_1\left(\frac{1}{6},\frac{5}{6};1;-\frac{\zeta}{2}\right),
\end{align*}
where in the last step we use the Euler formula for hypergeometric functions.\footnote{The Euler formula is \smash{$
{}_{2}F_1 (a,b;c;x )=\frac{\Gamma(c)}{\Gamma(b)\Gamma(c-b)}\int_0^1 t^{b-1}(1-t)^{c-b-1}(1-xt)^{-a}\,{\rm d}t$}.
}
By taking derivatives, we recover $\hat{\phi}_{\pm}(\zeta)$
\begin{align*}
\hat{\phi}_+(\zeta)&=\delta-\frac{1}{2}\frac{5}{36} {}_2F_1\left(\frac{7}{6},\frac{11}{6};2;-\frac{\zeta}{2}\right)=\delta-\frac{2}{3}\frac{5}{48} {}_2F_1\left(\frac{7}{6},\frac{11}{6};2;-\frac{\zeta}{2}\right),\\
\hat{\phi}_-(\zeta)&=\delta+\frac{1}{2}\frac{5}{36} {}_2F_1\left(\frac{7}{6},\frac{11}{6};2;\frac{\zeta}{2}\right)=\delta+\frac{2}{3}\frac{5}{48} {}_2F_1\left(\frac{7}{6},\frac{11}{6};2;\frac{\zeta}{2}\right)
\end{align*}
and up to a multiplicative constant they match our computations for $\borel\bigl(\tilde{W}_{1,2}\bigr)$ (see equation~\eqref{Borel-W12}).
The main advantage of writing Gauss hypergeometric functions as a convolution product relies on \'Ecalle's singularity theory. Indeed \smash{$\bigl(2\zeta+\zeta^2\bigr)^{-1/6}$} extends analytically to the universal cover of $\C\setminus\lbrace 0,-2\rbrace$ and the convolution with $\zeta^{-5/6}$ does not change the set of singularities \cite[Section~6.14.5\,(c)]{diverg-resurg-i}. Furthermore, Mitschi and Sauzin prove that $\hat{\phi}_{\pm}(\zeta)$ are simple resurgent functions \cite[Lemma 6.106]{diverg-resurg-i}.
\subsubsection{Comparison with exact WKB}
In \cite{kawai-takei}, Kawai and Takei carry out an exact WKB analysis of the Airy-type Schr\"{o}dinger equation%
\begin{equation}
\label{WKB_Airy}
\left[\left(\frac{\rm d}{{\rm d}x}\right)^2 - \eta^2 x \right] \psi(x, \eta) = 0
\end{equation}
in the $\eta \to \infty$ limit. They define $\psi_B(x, y)$ as the inverse Laplace transform of $\psi(x, \eta)$ with respect to $\eta$. In the coordinate $t=\frac{3}{2}yx^{-3/2}$ they find an explicit formula for $\psi_B(x,y)$ in terms of Gauss hypergeometric functions
\begin{align*}
\psi_{+,B}(x,y)&=\frac{1}{x}\phi_+(t)=\frac{\sqrt{3}}{2\sqrt{\pi}}\frac{1}{x}s^{-1/2} {}_2F_1\left(\frac{1}{6},\frac{5}{6};\frac{1}{2};s\right),\\
\psi_{-,B}(x,y)&=\frac{1}{x}\phi_-(t)=\frac{\sqrt{3}}{2\sqrt{\pi}}\frac{1}{x}(1-s)^{-1/2} {}_2F_1\left(\frac{1}{6},\frac{5}{6};\frac{1}{2};1-s\right),
\end{align*}
where $s=t/2+1/2$. We found the same hypergeometric functions in the position domain while solving the Airy equation
\smash{$
\big[\big(\frac{\rm d}{{\rm d}w}\big)^2 - w \big] f(w) = 0
$}
in Appendix~\ref{sec:spec-to-airy}. Although the two equations look closely related, and they are equivalent under the change of coordinates $w=x\eta^{2/3}$, they differ in an important way. The Borel transform of~$\psi$ is computed with respect to $\frac{2}{3}\eta x^{3/2}$, the conjugate variable of $t$, while the Borel transform of~$f(w)$ is computed with respect to $w$. We need to find a different change of coordinates to explain why the Borel transforms of $\psi(x,\eta)$ and~$f(w)$ are given by the same hypergeometric function.

To do this, first notice that if $\eta$ and $y$ are conjugate variables under Borel transform, meaning
\begin{align*}
\sum_{n\geq 0}a_n\eta^{-n-1} \overset{\borel}{\longrightarrow} \sum_{n\geq 0}\frac{a_n}{n!} y^{n}
\end{align*}
then \smash{$t=\frac{3}{2}yx^{-3/2}$} is the conjugate variable of \smash{$q=\frac{2}{3}\eta x^{3/2}$} up to correction by a factor of \smash{$\frac{3}{2}x^{-3/2}$}
\begin{gather*}
\sum_{n\geq 0}a_nq^{-n-1}=\sum_{n\geq 0}a_nx^{-3/2(n+1)}\left(\frac{2}{3}\eta\right)^{-n-1}\\
\qquad \overset{\borel}{\longrightarrow} \sum_{n\geq 0}\frac{a_nx^{-3/2(n+1)}}{n!}\left(\frac{3}{2}\right)^{n+1} y^{n}=\frac{3}{2}x^{-3/2}\sum_{n\geq 0}\frac{a_n}{n!} t^{n}.
\end{gather*}
In addition, $\psi_{B,\pm}(x,y)=\frac{1}{x}\phi_{\pm}(t)$, therefore we expect that $\psi(x,\eta)=x^{1/2}\Phi(q)$. Assume that~$\psi(x,y)$ is a solution of \eqref{WKB_Airy}. Then $\Phi(q)$ satisfies
\begin{equation}
\label{eq_Phi}
\left[\left(\frac{\rm d}{{\rm d}x}\right)^2+x^{-1}\frac{\rm d}{{\rm d}x}-\frac{1}{4}x^{-2} - \eta^2 x \right] \Phi(q) = 0,
\end{equation}
because
\begin{gather*}
\left[\left(\frac{\rm d}{{\rm d}x}\right)^2 - \eta^2 x \right] \psi(x, \eta) = 0,\qquad
\left[\left(\frac{\rm d}{{\rm d}x}\right)^2 - \eta^2 x \right] x^{1/2}\Phi(q) = 0,\\
\frac{\rm d}{{\rm d}x}\left[\frac{1}{2}x^{-1/2}\Phi+x^{1/2}\frac{\rm d}{{\rm d}x}\Phi\right]-\eta^2x^{3/2}\Phi=0,\\
-\frac{1}{4}x^{-3/2}\Phi+\frac{1}{2}x^{-1/2}\frac{\rm d}{{\rm d}x}\Phi+\frac{1}{2}x^{-1/2}\frac{\rm d}{{\rm d}x}\Phi+x^{1/2}\left(\frac{\rm d}{{\rm d}x}\right)^2\Phi-\eta^2x^{3/2}\Phi=0,\\
\left[x^{1/2}\left(\frac{\rm d}{{\rm d}x}\right)^2+x^{-1/2}\frac{\rm d}{{\rm d}x}-\frac{1}{4}x^{-3/2}-\eta^2x^{3/2}\right]\Phi=0,\\
\left[\left(\frac{\rm d}{{\rm d}x}\right)^2+x^{-1}\frac{\rm d}{{\rm d}x}-\frac{1}{4}x^{-2}-\eta^2x\right]\Phi=0.
\end{gather*}
Now, rewrite \eqref{eq_Phi} in the coordinates $q=\frac{2}{3}\eta x^{3/2}$,
\begin{gather*}
\left[\left(\frac{\rm d}{{\rm d}x}\right)^2+x^{-1}\frac{\rm d}{{\rm d}x}-\frac{1}{4}x^{-2}-\eta^2x\right]\Phi=0,\\
\left[\eta^2x\left(\frac{\rm d}{{\rm d}q}\right)^2+\frac{1}{2}\eta x^{-1/2}\frac{\rm d}{{\rm d}q}+x^{-1}\cdot \eta x^{1/2}\frac{\rm d}{{\rm d}q}-\frac{1}{4}x^{-2}-\eta^2x\right]\Phi=0,\\
\left[\eta^2x\left(\frac{\rm d}{{\rm d}q}\right)^2+\frac{3}{2}\eta x^{-1/2}\frac{\rm d}{{\rm d}q}-\frac{1}{4}x^{-2}-\eta^2x\right]\Phi=0,\\
\left[\left(\frac{\rm d}{{\rm d}q}\right)^2+\frac{3}{2}\eta^{-1} x^{-3/2}\frac{\rm d}{{\rm d}q}-\frac{1}{4}\eta^{-2}x^{-3}-1\right]\Phi=0,\\
\left[\left(\frac{\rm d}{{\rm d}q}\right)^2+q^{-1}\frac{\rm d}{{\rm d}q}-\frac{1}{9}q^{-2}-1\right]\Phi=0.
\end{gather*}
This shows that $\Phi(q)$ satisfies the transformed Airy equation \eqref{eqn:reg-mod-bessel}.

\begin{Remark}
The new coordinate $q = \frac{2}{3}\eta x^{3/2}$ is carefully chosen. Recall that the WKB ansatz for a Schr\"{o}dinger equation is
\[
\psi(x,\eta) = \exp\left(\int_{x_0}^xS(\eta,x')\, {\rm d}x'\right),
\]
with $S(\eta,x)=\sum_{k\geq -1} S_k(x) \eta^{-k}$. For the Airy-type Schr\"odinger equation, we have $S_{-1}^2 = x$, so up to a choice of sign for the square root,
\[ q = \frac{2}{3}\eta x^{3/2}=\eta\int_0^x (x')^{1/2} \, {\rm d}x' = \eta\int_{0}^x S_{-1}(x')\, {\rm d}x'. \]
We expect that the change of coordinates $q=\eta\int_0^{x}S_{-1}(x') \,{\rm d}x'$ would explain the analogies between the Borel transform of the WKB solution of a Schr\"{o}dinger equation and the Borel transform of the associated ODE.
\end{Remark}

\section{Order shifting for integro-differential equations}\label{shifting}
For a differential equation with a regular singularity at $\zeta = 0$, the Frobenius method typically gives a convergent power series solution in $\zeta^\rho \C\{\zeta\}$ for some $\rho \in \R$, which sums to a germ in~$\zeta^\rho \mathcal{O}_{\zeta=0}$. When $\rho$ is a non-integer greater than $-1$, we can do a nice trick in this space of germs: we can pass back and forth between differential and integral equations without worrying about boundary conditions.
\begin{Proposition}\label{prop:shifting}
Consider a germ $\phi \in \zeta^\rho \mathcal{O}_{\zeta=0}$. When $\rho$ is a non-integer greater than~$-1$,
\[ \left[ \sum_{k = 0}^n \left(\frac{\partial}{\partial \zeta}\right)^k \circ h_k + \sum_{k = 1}^m \fracderiv{-k}{\zeta}{0} \circ h_{-k} \right] \phi = 0, \]
if and only if
\[ \left[ \sum_{k = 1}^n \left(\frac{\partial}{\partial \zeta}\right)^{k-1} \circ h_k + \sum_{k = 0}^m \fracderiv{-k-1}{\zeta}{0} \circ h_{-k} \right] \phi = 0, \]
assuming that the coefficients $h_n, \dots, h_{-m}$ are in $\mathcal{O}_{\zeta=0}$.
\end{Proposition}
\begin{proof}
The reverse implication holds without any special condition on $\phi$, because $\frac{\partial}{\partial \zeta} \fracderiv{-1}{\zeta}{0}$ acts as the identity on functions which are locally integrable at $\zeta = 0$.

To prove the forward implication, rewrite the first equation in the statement as
\begin{equation}\label{eqn:diff-int-split}
\frac{\partial}{\partial \zeta} \left[ \sum_{k = 1}^n \left(\frac{\partial}{\partial \zeta}\right)^{k-1} \circ h_k \right] \psi = -\left[ h_0 + \sum_{k = 1}^m \fracderiv{-k}{\zeta}{0} \circ h_{-k} \right] \psi.
\end{equation}
The germ
\[ f := \left[ \sum_{k = 1}^n \left(\frac{\partial}{\partial \zeta}\right)^{k-1} \circ h_k \right] \psi \]
belongs to $\zeta^{\rho-(n-1)} \mathcal{O}_{\zeta = 0}$. This tells us, in particular, that $f$ has a shifted Taylor expansion with no constant term. Looking at the right-hand side of equation~\eqref{eqn:diff-int-split}, we can also see that $\frac{\partial}{\partial \zeta} f$ belongs to \smash{$\zeta^\rho \mathcal{O}_{\zeta=0}$}. Combining these observations, we can deduce that $f$ belongs to \smash{$\zeta^{\rho+1} \mathcal{O}_{\zeta=0}$}. Since $\rho > -1$, that means $f$ vanishes at $\zeta = 0$.

Integrating both sides of equation~\eqref{eqn:diff-int-split}, we get
\[ \fracderiv{-1}{\zeta}{0} \frac{\partial}{\partial \zeta} f = -\left[ \fracderiv{-1}{\zeta}{0} \circ h_0 + \sum_{k = 1}^m \fracderiv{-k-1}{\zeta}{0} \circ h_{-k} \right] \phi. \]
Since $\fracderiv{-1}{\zeta}{0} \frac{\partial}{\partial \zeta}$ acts as the identity on functions that vanish at $\zeta = 0$, this simplifies to
\[ f = -\left[ \sum_{k = 0}^m \fracderiv{-k-1}{\zeta}{0} \circ h_{-k} \right] \phi, \]
which rearranges to the second equation in the statement.
\end{proof}

\subsection*{Acknowledgements}
This paper is a result of the ERC-SyG project, Recursive and Exact New Quantum Theory (ReNewQuantum) which received funding from the European Research Council (ERC) under the European Union's Horizon 2020 research and innovation programme under grant agreement No 810573.

We thank Fondation Math\'ematique Jacques Hadamard for supporting the visit of the second author at IH\'{E}S, under the program \textit{Junior Scientific Visibility}.

We thank Fr\'ed\'eric Fauvet, Maxim Kontsevich, Andrew Neitzke, and David Sauzin for fruitful discussions and suggestions. We thank the referees for their careful reading, thoughtful comments, and useful references.

\pdfbookmark[1]{References}{ref}
\LastPageEnding

\end{document}